\documentstyle{amsppt}
\pagewidth{5.26in} \pageheight{7.86in} \NoRunningHeads
\magnification = 1200

\topmatter
\title On $(H,\widetilde{H})$-harmonic Maps between pseudo-Hermitian manifolds*\endtitle
\author Yuxin Dong \endauthor
\thanks {*Supported by NSFC grant No. 11271071, and LMNS, Fudan.}
\endthanks
\abstract {In this paper, we investigate critical maps of the
horizontal energy functional $E_{H,\widetilde{H}}(f)$ for maps
between two pseudo-Hermitian manifolds $(M^{2m+1},H(M),J,\theta )$
and $(N^{2n+1},\widetilde{H}(N), \widetilde{J},\widetilde{\theta
})$. These critical maps are referred to as $
(H,\widetilde{H})$-harmonic maps. We derive a CR Bochner formula
for the horizontal energy density $|df_{H, \widetilde{H}}|^{2}$,
and introduce a Paneitz type operator acting on maps to refine the
Bochner formula. As a result, we are able to establish some
Bochner type theorems for $(H,\widetilde{H})$-harmonic maps. We
also introduce $(H,\widetilde{H})$-pluriharmonic,
$(H,\widetilde{H})$-holomorphic maps between these manifolds,
which provide us examples of $(H,\widetilde{H})$-harmonic maps.
Moreover, a Lichnerowicz type result is established to show that
foliated $(H,\widetilde{ H})$-holomorphic maps are actually
minimizers of $E_{H,\widetilde{H}}(f)$ in their foliated homotopy
classes. We also prove some unique continuation results for
characterizing either horizontally constant maps or foliated
$(H,\widetilde{H})$-holomorphic maps. Furthermore, Eells-Sampson
type existence results for $(H,\widetilde{H})$-harmonic maps are
established if both manifolds are compact Sasakian and the target
is regular with non-positive horizontal sectional curvature.
Finally, we give a foliated rigidity result for
$(H,\widetilde{H})$-harmonic maps and Siu type strong rigidity
results for compact regular Sasakian manifolds with either
strongly negative horizontal curvature or adequately negative
horizontal curvature.}
\endabstract
\endsubjclass
\endkeywords
\endtopmatter
\document
\heading{\bf Introduction}
\endheading
\vskip 0.3 true cm

A smooth map $f$ between two Riemannian manifolds $M$ and $N$ is
called harmonic if it is a critical point of the energy functional
(cf. [EL])
$$
E(f)=\frac{1}{2}\int_{M}|df|^{2}dv_{M}.
$$
Harmonic maps became a useful tool for studying complex structures
of K\"{a}hler manifolds through the fundamental work of Siu
[Si1,2]. In his generalization of Mostow's rigidity theorem for
Hermitian symmetric spaces, Siu proved that a harmonic map of
sufficiently high maximum rank of a compact K\"{a}hler manifold to
a compact K\"{a}hler manifold with strongly negative curvature or
a compact quotient of an irreducible bounded symmetric domain must
be holomorphic or anti-holomorphic. It follows that if a compact
K\"{a}hler manifold is homotopic to such a target K\"{a}hler
manifold, then the homotopy equivalent map is homotopic to a
biholomorphic or anti-biholomorphic map. The latter result is
usually known as Siu's strong rigidity theorem. Although harmonic
maps are successful in the study of geometric and topological
structures of K\"{a}hler manifolds, they are not always effective
in other settings. For example, it is known that harmonic maps are
no longer in force for investigating general Hermitian manifolds,
since holomorphic maps between these manifolds are not necessarily
harmonic. In recent years, some generalized harmonic maps were
introduced and investigated in various geometric backgrounds (cf.
[JY], [Kok], [BD], [BDU], [KW], [Pe], [CZ]).

First, we recall the notion of transversally harmonic maps between
Riemannian foliations.  Let $(M,g,F)$ and $(N,h,\widetilde{F})$ be
two compact Riemann manifolds with Riemannian foliations $F$ and
$\widetilde{F}$ respectively, and $ f:M\rightarrow N$ a smooth
foliated map. Denote by $\upsilon (F)$ and $ \upsilon
(\widetilde{F})$ the normal bundles of the foliations $F$ and $
\widetilde{F}$ respectively. The differential $df$ gives rise
naturally to a smooth section $d_{T}f$ of $Hom(\upsilon
(F),f^{-1}\upsilon (\widetilde{F}))$. Then we may define the
transverse energy
$$
E_{T}(f)=\frac{1}{2}\int_{M}|d_{T}f|^{2}dv_{M}
$$
where $d_{T}f:\upsilon (F)\rightarrow \upsilon (\widetilde{F})$ is
the induced map of the differential map $df$, called the
transversally differential map. A smooth foliated map
$f:(M,g,F)\rightarrow (N,h, \widetilde{F})$ is called
transversally harmonic if it is an extremal of $ E_{T}(\cdot )$
for any variation of $f$ by foliated maps (cf. [BD], [KW]).
Suppose now that the foliations $F$ and $\widetilde{F}$ are two
K\"{a}hlerian foliations with complex structures $J$ and
$\widetilde{J}$ on their normal spaces $\upsilon (F)$ and
$\upsilon (\widetilde{F})$ respectively. A foliated map
$f:M\rightarrow N$ is said to be transversally holomorphic if $
d_{T}f\circ J=\widetilde{J}\circ d_{T}f$. It was proved in [BD]
that any transversally holomorphic map $f:(M,g,F)\rightarrow (N,h,
\widetilde{F})$ of K\"ahlerian foliations with $F$ harmonic is
transversally harmonic.

Recall that a CR structure on an $(2m+1)$-dimensional manifold
$M^{2m+1}$ is an $2m$-dimensional distribution $H(M)$ endowed with
a formally integrable complex structure $J$. The manifold $M$ with
the pair $(H(M),J)$ is called a CR-manifold. A pseudo-Hermitian
manifold, which is an odd-dimensional analogue of Hermitian
manifolds, is a CR manifold $M$ endowed with a pseudo-Hermitian
structure $\theta $. The pseudo-Hermitian structure $\theta $
determines uniquely a global nowhere zero vector field $\xi $ and
a Riemannian metric $g_{\theta }$ on $M$. The integral curves of
$\xi $ forms a foliation $F_{\xi }$, called the Reeb foliation.
From [Ta], [We], we know that each pseudo-Hermitian manifold
admits a unique canonical connection $\nabla $ (the Tanaka-Webster
connection), which is compatible with both the metric $g_{\theta
}$ and the CR structure (see Proposition 1.1). However, this
canonical connection always has nonvanishing torsion $T_{\nabla
}(\cdot ,\cdot )$, whose partial component $T_{\nabla }(\xi ,\cdot
)$ is an important pseudo-Hermitian invariant, called the
pseudo-Hermitian torsion. Pseudo-Hermitian manifolds with
vanishing pseudo-Hermitian torsion are referred to as Sasakian
manifolds which play an important role in  AdS/CFT correspondence
stemming from string theory (cf. [MSY]). An equivalent
characterization for a pseudo-Hermitian manifold to be Sasakian is
that $F_{\xi }$ is a Riemannian foliation. It is known that
Sasakian geometry sits naturally in between two K\"ahler
geometries. On the one hand, Sasakian manifolds are the bases of
metric cones which are K\"ahler. On the other hand, the Reeb
foliation $F_\xi$ of a Sasakian manifold is a K\"ahlerian
foliation (cf. [BG], [BGS]). Due to these two aspects, a Sasakian
manifold can be viewed as an odd-dimensional analogue of a
K\"ahler manifold. Consequently, it is natural to expect similar
Siu type rigidity theorems for Sasakian manifolds. In [CZ], an
interesting Siu type holomorphicity result was asserted for
transversally harmonic maps between Sasakian manifolds when the
target has strongly negative transverse curvature. However, Siu
type strong rigidity theorems have not been established for
Sasakian manifolds yet.

For a map $f:(M^{2m+1},H(M),J,\theta )\rightarrow
(N^{2n+1},\widetilde{H}(N),\widetilde{J},\widetilde{\theta })$
between two pseudo-Hermitian manifolds, Petit [Pe] defined a
natural horizontal energy functional
$$
E_{H,\widetilde{H}}(f)=\frac{1}{2}\int_{M}|df_{H,\widetilde{H}
}|^{2}dv_{\theta}
$$
where $|df_{H,\widetilde{H}}|^{2}$ denotes the horizontal energy
density (cf. \S 3), and he called a critical map of
$E_{H,\widetilde{H}}(f)$ a pseudoharmonic map. The main purpose of
[Pe] is to derive Mok-Siu-Yeung type formulas for horizontal maps
from compact contact locally sub-symmetric spaces into
pseudo-Hermitian manifolds and obtain some rigidity theorem for
the horizontal pseudoharmonic maps. Note that the Euler-Lagrange
equation for $E_{H,\widetilde{H}}(f)$ derived in [Pe] contains an
extra condition on the pull-back torsion (see (3.10)). Besides,
the authors in [BDU] introduced another kind of pseudoharmonic
maps from a pseudo-Hermitian manifold to a Riemannian manifold. To
avoid the extra torsion condition in the Euler-Lagrange equation
of [Pe] and any possible confusion with the notion of
pseudo-harmonic maps in [BDU], we modify Petit's definition
slightly by restricting the variational vector field to be
horizontal and refer to the corresponding critical maps as
$(H,\widetilde{H})$-harmonic maps. Although the energy functional
$ E_{H,\widetilde{H}}(f)$ is defined in a way similar to that of
the transversal energy functional, we would like to point out the
differences between the notions of $(H, \widetilde{H})$-harmonic
maps and transversally harmonic maps. First, the Reeb foliation of
a pseudo-Hermitian manifold is not a Riemannian foliation in
general; secondly, a $(H,\widetilde{H})$-harmonic map is not a
priori required to be a foliated map; thirdly, a horizontal
variational vector field is not necessarily foliated too.
Therefore the starting points for their definitions are different,
though they may coincide for foliated maps between Sasakian
manifolds. Actually we will see that $(H,\widetilde{H})$-harmonic
maps may display some geometric phenomenon which are invisible
from transversally harmonic maps.

In this paper, we will study some basic geometric properties and
problems for $(H,\widetilde{H})$-harmonic maps such as
Bochner-type, Lichnerowicz-type and Eells-Sampson-type, Siu type
holomorphicity results, etc. Our main aim is to utilize
$(H,\widetilde{H})$-harmonic maps to establish Siu type strong
rigidity theorems for Sasakian manifolds. The paper is organized
as follows. Section 1 begins to recall some basic facts and
notions of pseudo-Hermitian geometry, including some properties of
the curvature tensor of a pseudo-Hermitian manifold. Next, we
introduce the notions of strongly negative or strongly
semi-negative horizontal curvature, including adequately negative
horizontal curvature, for Sasakian manifolds. Some model Sasakian
spaces with either strongly negative or adequately negative
horizontal curvature are given. In Section 2, we introduce the
second fundamental form $\beta (\cdot ,\cdot)$ with respect to the
Tanaka-Webster connections for a map $f$ between two
pseudo-Hermitian manifolds, and derive the commutation relations
of its covariant derivatives. In Section 3, we recall the
definition of a $(H,\widetilde{H})$-harmonic map and the
corresponding Euler-Lagrange equation (see (3.4))
$$
\tau _{H,\widetilde{H}}(f)=0,\tag{0.1}
$$
where $\tau _{H,\widetilde{H}}(f)$ is called the horizontal
tension field of $f$. Then some relationship among
$(H,\widetilde{H})$-harmonic maps, pseudo-harmonic maps and
harmonic maps are discussed. It turns out that although
$(H,\widetilde{H})$-harmonic maps have many nice properties
related to the pseudo-Hermitian structures, the PDE system (0.1)
is too degenerate and thus its solutions may not be regular enough
to detect global geometric properties, including the strong
rigidity, of pseudo-Hermitian manifolds. Actually transversally
harmonic maps between Riemannian foliations also have similar
drawbacks. In order to repair these drawbacks, we define a special
kind of $(H,\widetilde{H})$-harmonic maps as follows. For a map
$f:M\rightarrow N$ between two pseudo-Hermitian manifolds, we set
$\tau_{H}(f)=\tau_{H,\widetilde{H}}(f)+\tau_{H,\widetilde{L}}(f)$,
where $\tau_{H,\widetilde{L}}(f)$ is the vertical component of
$tr_{g_{\theta }}(\beta |_{H})$. For our purpose, we introduce a
nonlinear subelliptic system of equations
$$
\tau _{H}(f)=0,\tag{0.2}
$$
imposed on the map $f$. Since (0.2) implies (0.1), a solution of
(0.2) is referred to as a special $(H,\widetilde{H})$-harmonic map
(see Definition 3.2).  Special $(H,\widetilde{H})$-harmonic maps
will play an important role in our studying of the strong rigidity
for Sasakian manifolds. In Section 4, we derive a CR Bochner
formula for the horizontal energy density $|df_{H,\widetilde{H}
}|^{2}$, whose main difficulty in applications comes from a mixed
term consisting of some contractions of $df_{H,\widetilde{H}}$ and
$\beta (\cdot ,\xi )$. In order to deal with this term, we
introduce a Paneitz type operator acting on the map, which enables
us to refine the Bochner formula. As a result, we are able to
establish some Bochner type theorems for
$(H,\widetilde{H})$-harmonic maps. In Section 5, we first define
the notions of $(H,\widetilde{H})$-pluriharmonic maps,
$(H,\widetilde{H})$-holomorphic maps and
$(H,\widetilde{H})$-biholomorphisms. It turns out that foliated
$(H,\widetilde{H})$-holomorphic maps are
$(H,\widetilde{H})$-pluriharmonic, and
$(H,\widetilde{H})$-pluriharmonic maps are foliated
$(H,\widetilde{H})$-harmonic. Next, we give a unique continuation
theorem which asserts that a foliated $(H,\widetilde{H})$-harmonic
map between two Sasakian manifolds must be
$(H,\widetilde{H})$-holomorphic on the whole manifold if it is
$(H,\widetilde{H})$-holomorphic on an open subset. Some examples
of $(H,\widetilde{H})$-holomorphic maps are also given. From
[BGS], we know that for a given Sasakian structure $S=(\xi ,\theta
,J,g_{\theta })$ on $M$, the Reeb vector field $\xi $ polarizes
the Sasakian manifold $(M,S)$, and the space $S(\xi ,J_{\upsilon
})$ of all Sasakian structures with the fixed Reeb vector field
$\xi $ and the fixed transverse holomorphic structure $J_{\upsilon
}$ on $\upsilon (F_{\xi })$ is an affine space. We show that
$id_{M}:(M,S_{1})\rightarrow (M,S_{2})$ for any $S_{1},S_{2}\in
S(\xi ,J_{\upsilon })$ is a foliated
$(H,\widetilde{H})$-biholomorphism. In addition, we discuss the
case when $id_{M}:$ $(M,S_{1})\rightarrow (M,S_{2})$ is a special
$(H,\widetilde{H})$-biholomorphism. In Section 6, we obtain a
Lichnerowicz type result which asserts that the difference of
horizontal partial energies for a foliated map is a smooth
foliated homotopy invariant. As an application, we deduce that a
foliated $(H,\widetilde{H})$-holomorphic map between two
pseudo-Hermitian manifolds is an absolute minimum of the
horizontal energy $E_{H,\widetilde{H}}(f)$. In Section 7, we study
the existence problem for (0.2) by looking at the following
subelliptic heat flow
$$
\cases \frac{\partial f_t}{\partial t}&=\tau_H(f_t)\\
 f|_{t=0}&=h\endcases\tag{0.3}
$$
where $h:M\rightarrow N$ is a smooth map. In order to show that a
solution of this system exists for all $t>0$ and converges to a
solution of (0.2) as $t\rightarrow \infty $, we impose a
non-positivity condition on the horizontal curvature of $N$. The
main result of this section asserts that if $h:M\rightarrow N$ is
a foliated map between two compact Sasakian manifolds and $N$ is
regular with non-positive horizontal sectional curvature, then
there exists a foliated special $(H,\widetilde{H})$-harmonic map
in the same foliated homotopy class as $h$. In Section 8, we first
give a foliated rigidity result which states that if $
f:M\rightarrow N$ is a $(H,\widetilde{H})$-harmonic map between
two compact Sasakian manifolds and the target $N$ has non-positive
horizontal curvature, then $f$ must be foliated. Next, we obtain a
$(H,\widetilde{H})$-holomorphicity result which asserts that a
$(H,\widetilde{H})$-harmonic map of sufficiently high maximum rank
of a compact Sasakian manifold to a compact Sasakian manifold with
either strongly negative horizontal curvature or adequately
negative horizontal curvature must be
$(H,\widetilde{H})$-holomorphic or
$(H,\widetilde{H})$-antiholomorphic.  Besides, we establish some
foliated strong rigidity theorems for Sasakian manifolds with
either strongly negative horizontal curvature or adequately
negative horizontal curvature (see Theorem 8.12 and Corollary
8.13). In Appendix A, we introduce another natural generalized
harmonic maps between pseudo-Hermitian manifolds, called
pseudo-Hermitian harmonic maps. First, we give a continuation
theorem about the foliated property for pseudo-Hermitian harmonic
maps. Next, we obtain a rigidity result which asserts that if
$f:M\rightarrow N$ is a pseudo-Hermitian harmonic map between two
compact Sasakian manifolds, and $N$ has non-positive horizontal
curvature, then $f$ is a foliated special
$(H,\widetilde{H})$-harmonic map. The latter result shows the
rationality for using special $(H,\widetilde{H})$-harmonic maps as
a tool in our study of global geometric and topological properties
of Sasakian manifolds. In Appendix B, we give explicit
formulations for both (0.2) and (0.3) , which are helpful for us
to understand the existence theory in Section 7. The method for
these formulations is possibly useful in studying the existence of
other generalized harmonic maps.

\heading{\bf 1. Pseudo-Hermitian Geometry}
\endheading
\vskip 0.3 true cm

In this section, we  collect some facts and notations concerning
pseudohermitian structures on CR manifolds (cf. [DTo], [BG] for
details).

\definition{Definition 1.1} Let $M^{2m+1}$ be a real $(2m+1)$-dimensional orientable $C^\infty $
manifold. A CR structure on $M$ is a complex rank $m$ subbundle
$H^{1,0}M$ of $TM\otimes C$ satisfying
\newline
(i) $H^{1,0}M\cap H^{0,1}M=\{0\}$
($H^{0,1}M=\overline{H^{1,0}M}$);
\newline
(ii) $[\Gamma (H^{1,0}M),\Gamma (H^{1,0}M)]\subseteq \Gamma
(H^{1,0}M)$.
\newline The pair $(M,H^{1,0}M)$ is called a CR manifold of CR dimension $m$.
\enddefinition

The complex subbundle $H^{1,0}M$ corresponds to a real subbundle
of $TM$:
$$
H(M)=Re\{H^{1,0}M\oplus H^{0,1}M\}\tag{1.1}
$$
which is called the Levi distribution. The Levi distribution
$H(M)$ admits a natural complex structure defined by $
J(V+\overline{V})=i(V-\overline{V})$ for any $V\in H^{1,0}M$.
Equivalently, the CR structure may be described by the pair
$(H(M),J)$.

Let $E$ be the conormal bundle of $H(M)$ in $T^{*}M$, whose fiber
at each point $x\in M$ is given by
$$E_x=\{\omega \in T_x^{*}M:\ker
\omega \supseteq H_x(M)\}. \tag{1.2}
$$
Since $M$ is assumed to be orientable, and the complex structure
$J$ induces an orientation on $H(M)$, it follows that the real
line bundle $E$ is orientable. Thus $E$ admits globally defined
nowhere vanishing sections.

\definition{Definition 1.2} A globally defined nowhere vanishing section $\theta \in \Gamma
(E)$ is called a pseudo-Hermitian structure on $M$. The Levi-form
$L_\theta $ associated with a pseudo-Hermitian structure $\theta$
is defined by
$$
L_\theta (X,Y)=d\theta (X,JY)\tag{1.3}
$$
for any $X,Y\in H(M)$. If $L_\theta $ is positive definite for
some $\theta$, then $(M,H(M),J)$ is said to be strictly
pseudoconvex.
\enddefinition

When $(M,H(M),J)$ is strictly pseudoconvex, it is natural to
orient $E$ by declaring a section $\theta$ to be positive if
$L_\theta $ is positive. Henceforth we will assume that
$(M,H(M),J)$ is a strictly pseudoconvex CR manifold and $\theta$
is a positive pseudo-Hermitian structure. The quadruple
$(M,H(M),J,\theta )$ is called a pseudo-Hermitian manifold.

On a pseudo-Hermitian manifold $(M,H(M),J,\theta )$, one may use
basic linear algebra to derive that $\ker \theta _x=H_x(M)$ for
each point $x\in M$, and there is a unique globally defined vector
field $\xi$ such that
$$
\theta (\xi)=1,\quad d\theta (\xi,\cdot )=0.  \tag{1.4}
$$
The vector field $\xi$ is referred to as the Reeb vector field.
The collection of all its integral curves forms an oriented
one-dimensional foliation $F_\xi$ on $M$, which is called the Reeb
foliation in this paper. Consequently there is a splitting of the
tangent bundle $TM$
$$
TM=H(M)\oplus L_{\xi},\tag{1.5}
$$
where $L_{\xi}$ is the trivial line bundle generated by $\xi$. Let
$\nu (F_\xi )$ be the vector bundle whose fiber at each point
$p\in M$ is the quotient space $T_pM/L_\xi$, and let
$\pi_\nu:TM\rightarrow \nu (F_\xi)$ be the natural projection.
Clearly $\pi_\nu |_{H(M)}:H(M)\rightarrow \nu (F_\xi)$ is a vector
bundle isomorphism.

Let $\pi _H:TM\rightarrow H(M)$ denote the natural projection
morphism. Set $G_\theta =\pi _H^{*}L_\theta$, that is,
$$G_\theta
(X,Y)=L_\theta (\pi _HX,\pi _HY)\tag{1.6}
$$
for any $X$, $Y\in TM$. Let us extend $J$ to a $(1,1)$-tensor
field on $M$ by requiring that
$$
J\xi=0. \tag{1.7}
$$
Then the integrability condition (ii) in Definition 1.1 implies
that $G_\theta$ is $J$-invariant. The Webster metric on
$(M,H(M),J,\theta )$ is a Riemannian metric defined by
$$
g_\theta =\theta \otimes \theta +G_\theta.\tag{1.8}
$$
It follows that
$$
\theta(X)=g_{\theta}(\xi,X),\quad
d\theta(X,Y)=g_{\theta}(JX,Y)\tag{1.9}
$$
for any $X,Y\in TM$. We find that (1.5) is actually an orthogonal
decomposition of $TM$ with respect to $g_{\theta}$. In terms of
terminology from foliation theory, $H(M)$ and $L_{\xi}$ are also
called the horizontal and vertical distributions respectively.
Clearly $\theta\wedge (d\theta )^m$ is, up to a constant, the
volume form of $(M,g_{\theta})$.

On a pseudo-Hermitian manifold, we have the following canonical
linear connection which preserves both the CR and the metric
structures.

\proclaim{Proposition 1.1 ([Ta], [We])} Let $(M,H(M),J,\theta )$
be a pseudo-Hermitian manifold. Then there exists a unique linear
connection $\nabla$ such that
\newline(i) $\nabla _X\Gamma
(H(M))\subset \Gamma (H(M))$ for any $X\in \Gamma (TM)$;
\newline (ii) $\nabla g_\theta =0$, $\nabla J=0$ (hence $\nabla
\xi=\nabla \theta =0$);
\newline (iii) The torsion $T_\nabla $ of
$\nabla $ is pure, that is, for any $X,Y\in H(M)$, $T_\nabla
(X,Y)=d\theta (X,Y)\xi$ and $T_\nabla (\xi,JX)+JT_\nabla
(\xi,X)=0$.
\endproclaim

The connection $\nabla$ in Proposition 1.1 is called the
Tanaka-Webster connection. Note that the torsion of the
Tanaka-Webster connection is always non-zero. The pseudo-Hermitian
torsion, denoted by $\tau$, is the $TM$-valued $1$-form defined by
$\tau (X)=T_\nabla(\xi,X)$. The anti-symmetry of $T_\nabla$
implies that $\tau(\xi)=0$. Using (iii) of Proposition 1.1 and the
definition of $\tau$, the total torsion of the Tanaka-Webster
connection may be expressed as
$$
T_\nabla (X,Y)=(\theta \wedge \tau )(X,Y)+d\theta (X,Y)\xi
\tag{1.10}
$$
for any $X,Y\in TM$. Set
$$ A(X,Y)=g_\theta (\tau X,Y)  \tag{1.11}
$$
for any $X,Y\in TM$. Then the properties of $\nabla$ in
Proposition 1.1 also imply that $\tau(H^{1,0}(M))\subset
H^{0,1}(M)$ and $A$ is a trace-free symmetric tensor field.

\proclaim{Lemma 1.2 (cf. [DTo])} The Levi-Civita connection
$\nabla ^\theta $ of $(M,g_\theta )$ is related to the
Tanaka-Webster connection by
$$
\nabla ^\theta =\nabla -(\frac 12d\theta +A)\xi+\tau \otimes
\theta +\theta \odot J\tag{1.12}
$$
where $(\theta \odot J)(X,Y)=\frac 12(\theta (X)JY+\theta (Y)JX)$
for any $ X,Y\in TM$.
\endproclaim
\remark{Remark 1.1} The Levi form in this paper is $2$ times that
one in [DTo]. Thus the coefficient of the term $\theta \odot J$ in
(1.12) is different from that in Lemma 1.3 of [DTo].
\endremark

\proclaim{Lemma 1.3}Let $(M^{2m+1},H(M),J,\theta)$ be a
pseudo-Hermitian manifold with the associated Tanaka-Webster
connection $\nabla$. Let $X$ and $\rho$ be a vector field and a
$1$-form on M respectively. Then
$$
div(X)=\sum_{A=0}^{2m}g_\theta (\nabla _{e_A}X,e_A)\quad\text{and
}\quad\delta (\rho )=-\sum_{A=0}^{2m}(\nabla _{e_A}\rho )(e_A)
$$
where $\{e_A\}_{A=0,1,...,2m}=\{\xi,e_1,...,e_{2m}\}$ is any
orthonormal frame field on $M$. Here $div(\cdot )$ and $\delta
(\cdot )$ denote the divergence and the codifferential
respectively.
\endproclaim
\demo{Proof} According to (1.12) and using the property $tr(A)=0$,
it is easy to verify that
$$
\nabla _{\xi}^\theta \xi=0,\quad \sum_{A=1}^{2m}\nabla
_{e_A}^\theta e_A=\sum_{A=1}^{2m}\nabla _{e_A}e_A. \tag{1.13}
$$
Since both $\nabla^\theta$ and $\nabla$ are metric connections, we
may employ (1.13) and the definition of $divX$ to find
$$
\aligned divX &=\sum_{A=0}^{2m}g_\theta (\nabla _{e_A}^\theta
X,e_A)=\sum_{A=0}^{2m}\{e_Ag_\theta (X,e_A)-g_\theta (X,\nabla
_{e_A}^\theta e_A)\}\\
&=\sum_{A=0}^{2m}g_\theta (\nabla_{e_A}X,e_A).\endaligned
$$
Similarly the codifferential of the $1$-form $\rho$ can be
computed as follows
$$
\aligned \delta (\rho ) &=-\sum_{A=0}^{2m}(\nabla _{e_A}^\theta
\rho )(e_A)=-\sum_{A=0}^{2m}\{e_A\rho (e_A)-\rho (\nabla
_{e_A}^\theta e_A)\}\\
&=-\sum_{A=0}^{2m}(\nabla_{e_A}\rho )(e_A)
\endaligned
$$
in view of (1.13) again. \qed
\enddemo

For a pseudo-Hermitian manifold $(M,H(M),J,\theta )$, the
sub-Laplacian operator $\bigtriangleup _H$ on a function $u\in
C^2(M)$ is defined by
$$
\bigtriangleup_H u=div(\nabla_Hu),
$$
where $\nabla_Hu=\pi _H(\nabla u)$ is the horizontal gradient of
$u$. In terms of a local orthonormal frame field
$\{e_{A}\}_{A=1}^{2m}$ of $H(M)$ on an open subset $U\subset M$ ,
the sub-Laplacian can be expressed as
$$
\bigtriangleup _{H}u=\sum_{A=1}^{2m}\{e_{A}(e_{A}u)-(\nabla
_{e_{A}}e_{A})u\}.  \tag{1.14}
$$
The non-degeneracy of the Levi form $L_{\theta }$ on $H(M)$
implies that $\{e_{A},[e_{A},e_{B}]\}_{1\leq A,B\leq 2m}$ spans
the tangent space $T_{p}M$ at each point $p\in U$. From [H\"{o}],
we know that $\bigtriangleup _{H}$ is a hypoelliptic operator.

For simplicity, we will denote by $\langle\cdot ,\cdot\rangle$ the
real inner product induced by $g_\theta $ on various tensor
bundles of $M$. Recall that the curvature tensor $R$ of the
Tanaka-Webster connection $\nabla$ is defined by
$$
R(X,Y)Z=\nabla _X\nabla _YZ-\nabla _Y\nabla _XZ-\nabla _{[X,Y]}Z
\tag{1.15}
$$
for any $X,Y,Z\in \Gamma (TM)$. Set $R(X,Y,Z,W)=\langle
R(Z,W)Y,X\rangle$. Then $R$ satisfies
$$
R(X,Y,Z,W)=-R(X,Y,W,Z)=-R(Y,X,Z,W),\tag{1.16}
$$
where the second equality is because of $\nabla g_\theta =0$.
However, the symmetric property $R(X,Y,Z,W)=R(Z,W,X,Y)$ is no
longer true for a general pseudo-Hermitian manifold due to the
failure of the first Bianchi identity.

The curvature tensor of $(M,\nabla)$ induces a morphism $Q:\wedge
^2TM\rightarrow \wedge ^2TM$ which is determined by
$$
\langle Q(X\wedge Y),Z\wedge W\rangle=R(X,Y,Z,W) \tag{1.17}
$$
for any $X,Y,Z,W\in TM$. The complex extension of $Q$ (resp. $R$)
to a morphism from $\wedge ^2TM^C$(resp. $\otimes ^4TM^C$) is
still denoted by the same notation. For a horizontal $2$ plane
$\sigma =span_R\{X,Y\}\subset H(M)$, we define the horizontal
sectional curvature of $\sigma$ by
$$
K^H(\sigma )=\frac{\langle Q(X\wedge Y),X\wedge Y\rangle}{\langle
X\wedge Y,X\wedge Y\rangle}. \tag{1.18}
$$
In particular, the horizontal holomorphic sectional curvature of a
horizontal holomorphic $2$-plane $\sigma =span\{X,JX\}\subset
H(M)$ is given by
$$
K_{hol}^H(\sigma )=\frac{\langle Q(X\wedge JX),X\wedge JX\rangle
}{\langle X\wedge JX,X\wedge JX\rangle}.  \tag{1.19}
$$

\definition{Definition 1.3} A pseudo-Hermitian manifold $(M,H(M),J,\theta)$ is called a
Sasakian manifold if its pseudo-Hermitian torsion $\tau$ is zero.
If $(M,H(M),J,\theta)$ is Sasakian, then the quadruple $(\xi
,\theta ,J,g_\theta )$ is referred to as a Sasakian structure on
$M$ with underlying CR structure $(H(M),J|_{H(M)})$.
\enddefinition

It turns out that if $(M,H(M),J,\theta)$ is a Sasakian manifold,
then its curvature tensor satisfies the Bianchi identities.
Consequently
$$
R(X,Y,Z,W) =R(Z,W,X,Y)\tag{1.20}
$$
for any $X,Y,Z,W\in TM$. Since $\nabla \xi=0$, it follows that if
one of the vectors $X,Y,Z$ and $W$ is vertical, then
$$
R(X,Y,Z,W)=0.\tag{1.21}
$$
Furthermore, in terms of $\nabla J=0$ and the $J$-invariance of
$G_{\theta}$, we have
$$
R(JX,JY,Z,W)=R(X,Y,JZ,JW)=R(X,Y,Z,W)\tag{1.22}
$$
for any $X,Y,Z,W\in H(M)$. When $X,Y,Z$ and $W$ vary in $H(M)$,
$R(X,Y,Z,W)$ may be referred to as the horizontal curvature
tensor, which will be denoted by $R^H$. Hence we discover that all
curvature information of the Sasakian manifold is contained in its
horizontal curvature tensor $R^H$ which enjoys the same properties
as the curvature tensor of a K\"ahler manifold.

For a Sasakian manifold $(M,H(M),J,\theta)$, it is known that its
Reeb foliation $F_{\xi}$ defines a Riemannian foliation (cf. [BG],
[DT]). In addition, $F_\xi $ is transversely holomorphic in the
following sense. There is an open covering $\{U_\alpha \}$ of $M$
together with a family of diffeomorphisms $\Phi _\alpha :U_\alpha
\rightarrow (-1,1)\times W_\alpha \subset R\times C^m$ and
submersions
$$\varphi _\alpha =\pi \circ \Phi _\alpha :U_\alpha \rightarrow
W_\alpha \subset C^m\tag{1.23}
$$
where $\pi:(-1,1)\times W_\alpha\rightarrow W_\alpha$ is the
natural projection, such that $\varphi_\alpha ^{-1}(\varphi_\alpha
(x))$ is just the integral curve of $\xi $ in $U_\alpha $ passing
through the point $x$, and when $U_\alpha \cap U_\beta \neq
\emptyset $ the map
$$
\varphi_\beta \circ \varphi_\alpha ^{-1}:\varphi _\alpha (U_\alpha
\cap U_\beta )\rightarrow \varphi _\beta (U_\alpha \cap U_\beta)
$$
is a biholomorphism. Such a triple
$(U_\alpha,\Phi_\alpha;\varphi_\alpha)$ is called a foliated
coordinate chart. For every point $x\in U_\alpha $ the
differential $d\varphi_\alpha :H_x(M)\rightarrow T_{\varphi
_\alpha (x)}W_\alpha $ is an isomorphism taking the complex
structure $J_x$ on $H_x(M)$ to that on $T_{\varphi_\alpha
(x)}W_\alpha $.

\definition {Definition 1.4}Let $M$ be a compact Sasakian manifold and
let $F_\xi$ be the Reeb foliation defined by $\xi$. Then the
foliation $F_\xi$ is said to be quasi-regular if there is a
positive integer $k$ such that each point has a foliated
coordinate chart $(U,\varphi )$ such that each leaf of $F_\xi $
passes through $U$ at most $k$ times, otherwise it is called
irregular. If $F_\xi $ is quasi-regular with the integer $k=1$,
then the foliation is called regular.
\enddefinition

It is known that the quasi-regular property is equivalent to the
condition that all the leaves of the foliation are compact. In the
quasi-regular case, the leaf space has the structure of a K\"ahler
orbifold. In the regular case, the foliation is simple so that the
Sasakian manifold can be realized as a $S^1$-bundle over a
K\"ahler manifold (the Boothby-Wang fibration [BW]), and the
natural projection of this fibration is actually a Riemannian
submersion. In general, in the irregular case, the leaf space is
not even Hausdorff.

\definition{Definition 1.5} A Sasakian manifold $M^{2m+1}$ with the
Tanaka-Webster connection $\nabla $ is said to have strongly
negative horizontal curvature (resp. strongly seminegative
horizontal curvature) if
$$
\langle Q(\zeta),\overline{\zeta} \rangle < 0
\quad(\text{resp.}\leq 0)
$$
for any $\zeta=(Z\wedge W)^{(1,1)}\neq 0$, $Z,W\in HM^C$. Here
$\overline{\zeta}$ is the complex conjugate of $\zeta$. In
addition, we say that the horizontal curvature tensor $R^H$ of a
Sasakian manifold $M^{2m+1}$ is negative of order $k$ if it is
strongly seminegative and it enjoys the following property. If
$A=(A_{\overline{i}}^\alpha )$, $B=(B_i^\alpha )$ are any two
$m\times k$ matrices ($1\leq \alpha \leq m$, $1\leq i\leq k$) with
$$
rank\pmatrix A& B\\
\overline{B}& \overline{A}\\
\endpmatrix
=2k
$$
and  if
$$
\sum_{\alpha ,\beta ,\gamma ,\delta }R_{\alpha \overline{\beta
}\gamma \overline{\delta }}\xi _{\overline{i}\overline{j}}^{\alpha
\overline{\beta }}\thinspace\overline{\xi
_{\overline{i}\overline{j}}^{\delta \overline{\gamma }}}=0
$$
for all $1\leq i,j\leq k$, where
$$
\xi_{\overline{i}\overline{j}}^{\alpha \overline{\beta
}}=A_{\overline{i}}^\alpha \overline{B_j^\beta
}-A_{\overline{j}}^\alpha \overline{B_i^\beta },$$ then either
$A=0$ or $B=0$. The horizontal curvature tensor $R^H$ is called
adequately negative if it is negative of order $m$.
\enddefinition

By the $J$-invariant property (1.22), we find that the curvature
operator $Q$ annihilates any $2$-vector of type (2,0) or (0,2).
Therefore strongly negativity (resp. strongly semi-negativity) of
the horizontal curvature tensor implies negativity (resp.
semi-negativity) of the horizontal sectional curvature.

\example{Example 1.1}
\newline \text{  }(i) A Sasakian manifolds
$(M^{2m+1},H(M),J,\theta )$ with $K_{hol}^H\ $constant is called a
Sasakian space form. For each real number $\lambda$, Webster [We]
gave a model $M(\lambda )$ for the Sasakian space form with
$K_{hol}^H=\lambda$. The horizontal curvature tensor of $M(\lambda
)$ may be expressed as
$$
R_{\alpha \overline{\beta }\gamma \overline{\delta }}=-\frac
\lambda 2(g_{\alpha \overline{\beta }}g_{\gamma \overline{\delta
}}+g_{\alpha \overline{\delta }}g_{\gamma \overline{\beta
}})\tag{1.24}
$$
with respect to any frame $\{\eta _\alpha \}$ of $H^{1,0}M$ at
every point. Following the method for complex ball in [Si1], one
may verify that if $\lambda <0$, then $M(\lambda )$ has strongly
negative horizontal curvature.
\newline \text{ }(ii) By a theorem of Kobayashi ([Ko]), we know that if $B$ is a compact
Hodge manifold with integral K\"ahler form, there exists a
Sasakian manifold $M$, which is the total space of a Riemannian
submersion over $B$. Suppose $M^{2m+1}$ is a compact Sasakian
manifold with a Riemannian submersion $\pi :M\rightarrow B$ over a
compact K\"ahler manifold. Denote by $\nabla ^\theta $ and $\nabla
^B$ the Levi-Civita connections of $M$ and $B$. Recall that a
basic vector field on $M$ is one that is both horizontal and
projectable. Suppose $X,Y$ are basic vector fields on $M$. Denote
by $\widetilde{X}$, $\widetilde{Y}$ the vector fields on $B$ that
are $\pi $-related to $X$ and $Y$. A basic property on the
connections of a Riemannian submersion (cf.[O'N], [GW]) gives that
$d\pi [\pi _H(\nabla _X^\theta Y)]=(\nabla
_{\widetilde{X}}\widetilde{Y})\circ \pi $, where $\pi
_H:TM\rightarrow H(M)$ is the natural projection. Since $M$ is
Sasakian, we know from Lemma 1.2 that $\pi _H(\nabla _X^\theta
Y)=\nabla _XY$. Thus $d\pi (\nabla _XY)=(\nabla
_{\widetilde{X}}\widetilde{Y})\circ \pi $, which implies that
$$
R(X,Y,Z,W)=R^B(d\pi(X),d\pi(Y),d\pi(Z),d\pi(W))\tag{1.25}
$$
for any $X,Y,Z,W\in H(M)$, where $R^B$ is the curvature tensor of
$\nabla ^B$. Hence $M$ has strongly curvature horizontal curvature
(resp. strongly semi-negative horizontal curvature with negative
order $k$) if and only if $B$ has strongly negative curvature
(resp. strongly semi-negative curvature with negative order $k$)
in the sense of [Si1]. In terms of [Si1], we find that if $B$ is a
compact quotient of an irreducible symmetric bounded domain, then
$M$ has adequately negative horizontal curvature. Since the
foliation of a Sasakian manifold is locally a Riemannian foliation
over a K\"ahler manifold, the above discussion helps us to
understand the general curvature properties of a Sasakian manifold
from those of a K\"ahler manifold. For example, it is proved in
[Si1] that if the curvature tensor of a K\"ahler manifold is
strongly negative, then it is negative of order 2. Therefore we
may conclude that if a Sasakian manifold has strongly negative
horizontal curvature, then its horizontal curvature tensor is
negative of order 2.

\endexample

\heading{\bf 2. Second fundamental forms and their covariant
derivatives}
\endheading
\vskip 0.3 true cm

Let $(M^{2m+1},H(M),J,\theta )$ and
$(N^{2n+1},\widetilde{H}(N),\widetilde{J} ,\widetilde{\theta })$
be two pseudo-Hermitian manifolds. Denote by $\nabla $ and
$\widetilde{\nabla }$ the Tanaka-Webster connections of $M$ and
$N$ respectively. Let $f:M\rightarrow N$ be a smooth map. Then the
bundle $T^{*}M\otimes f^{-1}TN$ has the induced connection $\nabla
\otimes f^{-1}\widetilde{\nabla }$, where $f^{-1}\widetilde{\nabla
}$ is the pull-back connection in $f^{-1}TN$. For simplicity, we
will write $f^{-1} \widetilde{\nabla }$ as $\widetilde{\nabla }$
when the meaning is clear. The second fundamental form of $f$ with
respect to $(\nabla , \widetilde{\nabla })$ is defined by:
$$
\aligned \beta (X,Y)&=\lbrack (\nabla \otimes
f^{-1}\widetilde{\nabla })_Ydf](X)\\
&=\widetilde{\nabla}_Ydf(X)-df(\nabla _YX)\endaligned\tag{2.1}
$$
for $X,Y\in \Gamma(TM)$. In what follows, we shall use the
summation convention for repeated indices.

\proclaim{Lemma 2.1} Let $f:(M,H(M),J,\theta )\rightarrow
(N,\widetilde{H}(N),\widetilde{J}, \widetilde{\theta })$ be a map.
Then
$$
\widetilde{\nabla }_Xdf(Y)-\widetilde{\nabla
}_Ydf(X)-df([X,Y])=\widetilde{T}_{\widetilde{\nabla
}}(df(X),df(Y))
$$
for any $X,Y\in \Gamma (TM)$, where
$\widetilde{T}_{\widetilde{\nabla }}$ denotes the torsion of the
Tanaka-Webster connection $\widetilde{\nabla }$ on $N$.
\endproclaim

\demo{Proof} Set $S(X,Y)=$ $\widetilde{\nabla
}_Xdf(Y)-\widetilde{\nabla } _Ydf(X)-df([X,Y])$. It is easy to
show that $S$ is $C^\infty (M)$-bilinear. Choose a local
coordinate chart $(x^A)$ around $p$ and a local coordinate chart
$(u^{\widetilde{A}})$ around $f(p)$. Then
$$
df(\frac \partial {\partial x^A})=\frac{\partial
f^{\widetilde{C}}}{\partial x^A}\frac \partial {\partial
u^{\widetilde{C}}}\tag{2.2}
$$
where $f^{\widetilde{C}}=u^{\widetilde{C}}\circ f$. By the
definition of $f^{-1} \widetilde{\nabla }$ and using (2.2) , we
deduce that
$$
\aligned S(\frac \partial {\partial x^A},\frac \partial {\partial
x^B})&=\widetilde{ \nabla }_{\frac \partial {\partial
x^A}}df(\frac
\partial {\partial x^B})- \widetilde{\nabla }_{\frac \partial
{\partial x^B}}df(\frac
\partial{\partial x^A}) \\
&=\widetilde{\nabla }_{\frac \partial {\partial
x^A}}(\frac{\partial f\widetilde{^C}}{\partial x^B}\frac \partial
{\partial u^{\widetilde{C}}})- \widetilde{\nabla }_{\frac \partial
{\partial x^B}}(\frac{\partial f^{\widetilde{C}}}{\partial x^A}\frac \partial {\partial u^{\widetilde{C}}}) \\
&=\frac{\partial f\widetilde{^C}}{\partial x^B}\widetilde{\nabla
}_{\frac\partial {\partial x^A}}\frac \partial {\partial
u^{\widetilde{C}}}-\frac{
\partial f^{\widetilde{C}}}{\partial x^A}\widetilde{\nabla }_{\frac \partial
{\partial x^B}}\frac \partial {\partial u^{\widetilde{C}}} \\
&=\frac{\partial f\widetilde{^C}}{\partial x^B}\frac{\partial
f^{\widetilde{D} }}{\partial x^A}[\widetilde{\nabla }_{\frac
\partial {\partial u^{\widetilde{ D}}}}\frac \partial {\partial
u^{\widetilde{C}}}-\widetilde{\nabla }_{\frac
\partial {\partial u^{\widetilde{C}}}}\frac \partial {\partial u^{\widetilde{D}}}] \\
&=\widetilde{T}_{\widetilde{\nabla }}(df(\frac \partial {\partial
x^A}),df(\frac \partial {\partial x^B})).
\endaligned
$$
Hence this lemma is proved. \qed
\enddemo

We will use the moving frame method to perform local computations
on maps between pseudo-Hermitian manifolds. Let us now recall the
structure equations of the Tanaka-Webster connection on a
pseudo-Hermitian manifold. For the pseudo-Hermitian manifold
$(M^{2m+1},H(M),J,\theta )$, we choose a local orthonormal frame
field $ \{e_A\}_{A=0}^{2m}=\{\xi,e_1,...,e_m,e_{m+1},.$
$..,e_{2m}\}$ with respect to $g_\theta$ such that
$$\{e_{m+1},...,e_{2m}\}=\{Je_1,...,Je_m\}.$$ Set
$$
\eta _j=\frac 1{\sqrt{2}}(e_j-iJe_j),\quad \eta
_{\overline{j}}=\frac 1{\sqrt{2}}(e_j+iJe_j)\quad  (j=1,...,m).
\tag{2.3}
$$
Then $\{\xi,\eta _j,\eta _{\overline{j}}\}$ forms a frame field of
$TM\otimes C $. Let $\{\theta ,\theta ^j,\theta ^{\overline{j}}\}$
be the dual frame field of $\{\xi,\eta _j,\eta _{\overline{j}}\}$.
From Proposition 1.1, we have
$$
\nabla _X\xi =0,\quad\nabla _X\eta _j =\theta _j^i(X)\eta
_i,\quad\nabla _X\eta _{\overline{j}}=\theta
_{\overline{j}}^{\overline{i}}(X)\eta _{\overline{i}}\tag{2.4}
$$
for any $X\in TM$, where $\{\theta _0^0=\theta _0^i=\theta_0^{
\overline{i}}=\theta_i^0=\theta_{\overline{i}}^0=0, \theta
_j^i,\theta _{\overline{j}}^{ \overline{i}}\}$ are the connection
$1$-forms of $\nabla $ with respect to the frame field $\{\xi,\eta
_j,\eta _{\overline{j}}\}$. According to (iii) of Proposition 1.1,
the pseudo-Hermitian torsion may be expressed as
$$
\tau =A_k^{\overline{j}}\theta ^k\otimes \eta
_{\overline{j}}+A_{\overline{k} }^j\theta ^{\overline{k}}\otimes
\eta _j.  \tag{2.5}
$$
The symmetry of $A$ implies that $A_{jk}=A_{kj}$
$=A_k^{\overline{j}}$ . Since $\nabla $ preserves $H^{1,0}M$, we
may write
$$
R(\eta _k,\eta _{\overline{l}})\eta _j=R_{jk\overline{l}}^i\eta
_i.\tag{2.6}
$$
From [We], we know that $\{\theta ,\theta ^i,\theta
^{\overline{i}},\theta _j^i,\theta
_{\overline{j}}^{\overline{i}}\}$ satisfies the following
structure equations (cf. also \S 1.4 of [DTo]):
$$
\cases
d\theta &=\sqrt{-1}\theta ^j\wedge \theta ^{\overline{j}} \\
d\theta ^i&=-\theta _j^i\wedge \theta ^j+A_{\overline{j}}^i\theta
\wedge\theta ^{\overline{j}} \\
d\theta _j^i&=-\theta _k^i\wedge \theta _j^k+\Psi _j^i
\endcases\tag{2.7}
$$
where
$$
\aligned \Psi _j^i=&W_{jk}^i\theta ^k\wedge \theta
-W_{j\overline{k}}^i\theta ^{\overline{k}}\wedge \theta+\sqrt{-1}\theta ^i\wedge A_k^{\overline{j}}\theta ^k \\
&-\sqrt{-1}A_{\overline{k}}^i\theta ^{\overline{k}}\wedge \theta
^{\overline{j}}+R_{jk\overline{l} }^i\theta ^k\wedge \theta
^{\overline{l}}
\endaligned\tag{2.8}
$$
and
$$
W_{kl}^j=A_{kl,\overline{j}},\quad
W_{k\overline{l}}^j=A_{\overline{l} \overline{j},k} . \tag{2.9}
$$

Let $\{\widetilde{\xi},\widetilde{\eta }_\alpha ,\widetilde{\eta
}_{\overline{ \alpha }}\}_{\alpha =1,...,n}$ be a local frame
field on the pseudo-Hermitian manifold $N^{2n+1}$, and let
$\{\widetilde{\theta },\widetilde{\theta }^\alpha
,\widetilde{\theta }^{ \overline{\alpha }}\}_{\alpha =1,...,n}$ be
its dual frame field. We will denote the connection $1$-forms,
torsion and curvature, etc., of the Tanaka-Webster connection
$\widetilde{\nabla}$ on $N$ by the same notations as in $M$, but
with $\widetilde{}$ on them. Then similar structure equations for
$\widetilde{\nabla}$ are valid in $N$ too. Henceforth we shall
make use of the following convention on the ranges of indices:
$$
\aligned
&A,B,C,...=0,1,...,m,\overline{1},...,\overline{m}; \\
&i,j,k,...=1,...,m,\
\overline{i},\overline{j},\overline{k},...=\overline{1}
,...,\overline{m}; \\
&\widetilde{A},\widetilde{B},\widetilde{C},...=0,1,...,n,\overline{1},...,
\overline{n}; \\
&\alpha ,\beta ,\gamma ,...=1,...,n,\ \overline{\alpha
},\overline{\beta }, \overline{\gamma
},...=\overline{1},...,\overline{n}. \endaligned
$$
As usual repeated indices are summed over the respective ranges.

For a map $f:M\rightarrow N$, we express its differential as
$$
df=f_B^{\widetilde{A}}\theta ^B\otimes \widetilde{\eta
}_{\widetilde{A }}.\tag{2.10}
$$
Therefore
$$
\cases f^{*}\widetilde{\theta }&=f_0^0\theta +f_j^0\theta
^j+f_{\overline{j}}^0\theta ^{\overline{j}} \\
f^{*}\widetilde{\theta }^\alpha &=f_0^\alpha \theta +f_j^\alpha
\theta ^j+f_{\overline{j}}^\alpha \theta ^{\overline{j}} \\
f^{*}\widetilde{\theta }^{\overline{\alpha
}}&=f_0^{\overline{\alpha }}\theta +f_j^{\overline{\alpha }}\theta
^j+f_{\overline{j}}^{\overline{\alpha } }\theta ^{\overline{j}}.
\endcases\tag{2.11}
$$

By taking the exterior derivative of the first equation in (2.11)
and making use of the structure equations in $M$ and $N$, we get
$$
\aligned &Df_0^0\wedge \theta +Df_j^0\wedge \theta
^j+Df_{\overline{j}}^0\wedge \theta
^{\overline{j}}+if_0^0\theta ^j\wedge \theta ^{\overline{j}} \\
&+f_j^0\theta \wedge A_{\overline{k}}^j\theta
^{\overline{k}}+f_{\overline{j} }^0\theta \wedge
A_k^{\overline{j}}\theta ^k-if^{*}\widetilde{\theta } ^\alpha
\wedge f^{*}\widetilde{\theta }^{\overline{\alpha }}=0 \endaligned
\tag{2.12}
$$
where
$$
\aligned Df_0^0&=df_0^0=f_{00}^0\theta +f_{0j}^0\theta
^j+f_{0\overline{j}}^0\theta ^{
\overline{j}} \\
Df_j^0&=df_j^0-f_k^0\theta _j^k=f_{j0}^0\theta +f_{jl}^0\theta
^l+f_{j
\overline{l}}^0\theta ^{\overline{l}} \\
Df_{\overline{j}}^0&=df_{\overline{j}}^0-f_{\overline{k}}^0\theta
_{\overline{ j}}^{\overline{k}}=f_{\overline{j}0}^0\theta
+f_{\overline{j}l}^0\theta ^l+f_{\overline{j}\overline{l}}^0\theta
^{\overline{l}}. \endaligned \tag{2.13}
$$
Then (2.12) gives
$$
\aligned
f_{j0}^0-f_{0j}^0+f_{\overline{k}}^0A_j^{\overline{k}}&=i(f_0^\alpha
f_j^{\overline{\alpha }}-f_0^{\overline{\alpha }}f_j^\alpha ) \\
f_{\overline{j}0}^0-f_{0\overline{j}}^0+f_k^0A_{\overline{j}}^k&=i(f_0^\alpha
f_{\overline{j}}^{\overline{\alpha }}-f_0^{\overline{\alpha
}}f_{\overline{j}}^\alpha)\\
f_{jl}^0-f_{lj}^0&=i(f_j^{\overline{\alpha }}f_l^\alpha
-f_j^\alpha f_l^{\overline{\alpha }}) \\
f_{\overline{j}\overline{l}}^0-f_{\overline{l}\overline{j}}^0&=i(f_{\overline{
j}}^{\overline{\alpha}}f_{\overline{l}}^\alpha-f_{\overline{j}}^\alpha
f_{\overline{l}}^{\overline{\alpha }}) \\
f_{j\overline{l}}^0-f_{\overline{l}j}^0-if_0^0\delta
_j^l&=i(f_j^{\overline{ \alpha }}f_{\overline{l}}^\alpha
-f_j^\alpha f_{\overline{l}}^{\overline{ \alpha }}).
\endaligned\tag{2.14}
$$

To simplify the notations, we will set $\widehat{\theta }_\beta
^\alpha =f^{*} \widetilde{\theta }_\beta ^\alpha $,
$\widehat{A}_\beta ^\alpha =f^{*} \widetilde{A}_\beta ^\alpha $,
$\widehat{\Psi }_\beta ^\alpha =f^{*} \widetilde{\Psi }_\beta
^\alpha $, etc. Similar computations for the second equation in
(2.11) yield
$$
\aligned &Df_0^\alpha \wedge \theta +Df_j^\alpha \wedge \theta
^j+Df_{\overline{j} }^\alpha
\wedge\theta^{\overline{j}}+if_0^\alpha \theta ^j\wedge \theta ^{
\overline{j}} \\
&+f_j^\alpha A_{\overline{k}}^j\theta \wedge \theta
^{\overline{k}}+f_{ \overline{j}}^\alpha A_k^{\overline{j}}\theta
\wedge \theta ^k=\widehat{A}_{ \overline{\beta }}^\alpha
f^{*}\widetilde{\theta }\wedge f^{*}\widetilde{ \theta
}^{\overline{\beta }} \endaligned\tag{2.15}
$$
where
$$
\aligned
Df_0^\alpha &=df_0^\alpha +f_0^\beta \widehat{\theta
}_\beta ^\alpha =f_{00}^\alpha \theta +f_{0j}^\alpha \theta
^j+f_{0\overline{j}}^\alpha\theta ^{\overline{j}}, \\
Df_j^\alpha &=df_j^\alpha -f_k^\alpha \theta _j^k+f_j^\beta
\widehat{\theta } _\beta ^\alpha =f_{j0}^\alpha \theta
+f_{jl}^\alpha \theta ^l+f_{j\overline{l
}}^\alpha \theta ^{\overline{l}}, \\
Df_{\overline{j}}^\alpha &=df_{\overline{j}}^\alpha
-f_{\overline{k}}^\alpha \theta
_{\overline{j}}^{\overline{k}}+f_{\overline{j}}^\beta
\widehat{\theta }_\beta ^\alpha =f_{\overline{j}0}^\alpha \theta
+f_{\overline{j}l}^\alpha \theta
^l+f_{\overline{j}\overline{l}}^\alpha \theta ^{\overline{l}}.
\endaligned\tag{2.16}
$$
From (2.15), it follows that
$$
\aligned f_{0j}^\alpha -f_{j0}^\alpha -f_{\overline{k}}^\alpha
A_j^{\overline{k}}&= \widehat{A}_{\overline{\beta }}^\alpha
(f_j^0f_0^{\overline{\beta }
}-f_0^0f_j^{\overline{\beta }}) \\
f_{0\overline{j}}^\alpha -f_{\overline{j}0}^\alpha -f_k^\alpha
A_{\overline{j }}^k&=\widehat{A}_{\overline{\beta }}^\alpha
(f_{\overline{j}}^0f_0^{
\overline{\beta }}-f_0^0f_{\overline{j}}^{\overline{\beta }}) \\
f_{jl}^\alpha -f_{lj}^\alpha &=\widehat{A}_{\overline{\beta
}}^\alpha
(f_l^0f_j^{\overline{\beta }}-f_j^0f_l^{\overline{\beta }}) \\
f_{\overline{j}l}^\alpha -f_{l\overline{j}}^\alpha +if_0^\alpha
\delta _l^j&= \widehat{A}_{\overline{\beta }}^\alpha
(f_l^0f_{\overline{j}}^{\overline{
\beta }}-f_l^{\overline{\beta }}f_{\overline{j}}^0)\\
f_{\overline{j}\overline{l}}^\alpha
-f_{\overline{l}\overline{j}}^\alpha &=
\widehat{A}_{\overline{\beta }}^\alpha
(f_{\overline{l}}^0f_{\overline{j}}^{ \overline{\beta
}}-f_{\overline{j}}^0f_{\overline{l}}^{\overline{\beta }}).
\endaligned\tag{2.17}
$$
Likewise we may deduce those commutative relations of $f_{AB}^{
\overline{\alpha }}$ from the third equation of (2.11) by taking
its exterior derivative, or directly from (2.17) by taking bar of
each index. Clearly the second fundament form $\beta $ can be
expressed as
$$
\beta =f_{BC}^{\widetilde{A}}\theta ^B\otimes \theta ^C\otimes
\widetilde{\eta }_{\widetilde{A}}.\tag{2.18}
$$
Due to the torsions of the Tanaka-Webster connections, $\beta$ is
not symmetric. Although the non-symmetry of $\beta$ causes a
little trouble, we will see that it may also lead to some
unexpected geometric consequences.

By taking the exterior derivative of the first equation of (2.13)
and making use of the structure equations, we get
$$
\aligned &Df_{00}^0\wedge \theta +Df_{0j}^0\wedge \theta
^j+Df_{0\overline{j}}^0\wedge \theta
^{\overline{j}}+if_{00}^0\theta ^j\wedge \theta
^{\overline{j}} \\
&+f_{0j}^0A_{\overline{k} }^j\theta \wedge \theta
^{\overline{k}}+f_{0\overline{j}}^0A_k^{\overline{j} }\theta
\wedge \theta ^k=0
\endaligned\tag{2.19}
$$
where
$$
\aligned Df_{00}^0&=df_{00}^0=f_{000}^0\theta +f_{00j}^0\theta
^j+f_{00\overline{j}}^0\theta ^{\overline{j}} \\
Df_{0j}^0&=df_{0j}^0-f_{0k}^0\theta _j^k=f_{0j0}^0\theta
+f_{0jl}^0\theta
^l+f_{0j\overline{l}}^0\theta ^{\overline{l}} \\
Df_{0\overline{j}}^0&=df_{0\overline{j}}^0-f_{0\overline{k}}^0\theta
_{\overline{j}}^{\overline{k}}=f_{0\overline{j}0}^0\theta
+f_{0\overline{j}l}^0\theta^l+f_{0\overline{j}\overline{l}}^0\theta
^{\overline{l}}. \endaligned\tag{2.20}
$$
It follows from (2.19) that
$$
\aligned
&f_{00j}^0-f_{0j0}^0=f_{0\overline{k}}^0A_j^{\overline{k}} \\
&f_{00\overline{j}}^0-f_{0\overline{j}0}^0=f_{0k}^0A_{\overline{j}}^k \\
&f_{0jl}^0=f_{0lj}^0\\
&f_{0\overline{j}l}^0-f_{0l\overline{j}}^0+if_{00}^0\delta _l^j=0 \\
&f_{0\overline{j}\overline{l}}^0=f_{0\overline{l}\overline{j}}^0.
\endaligned\tag{2.21}
$$
Similar computations for the second equation of (2.13) give
$$
\aligned
&Df_{j0}^0\wedge \theta +Df_{jl}^0\wedge \theta
^l+Df_{j\overline{l}}^0\wedge\theta ^{\overline{l}}+if_{j0}^0\theta ^k\wedge \theta ^{\overline{k}} \\
&+f_{jl}^0A_{\overline{k}}^l\theta \wedge \theta
^{\overline{k}}+f_{j\overline{l}}^0A_k^{\overline{l}}\theta \wedge
\theta ^k=-f_k^0\Psi _j^k \endaligned\tag{2.22}
$$
where
$$
\aligned &Df_{j0}^0=df_{j0}^0-f_{k0}^0\theta _j^k=f_{j00}^0\theta
+f_{j0l}^0\theta
^l+f_{j0\overline{l}}^0\theta ^{\overline{l}}\\
&Df_{jl}^0=df_{jl}^0-f_{kl}^0\theta _j^k-f_{jk}^0\theta
_l^k=f_{jl0}^0\theta+f_{jlk}^0\theta ^k+f_{jl\overline{k}}^0\theta ^{\overline{k}}\\
&Df_{j\overline{l}}^0=df_{j\overline{l}}^0-f_{k\overline{l}}^0\theta
_j^k-f_{j\overline{k}}^0\theta_{\overline{l}}^{\overline{k}}=f_{j\overline{l}
0}^0\theta +f_{j\overline{l}k}^0\theta
^k+f_{j\overline{l}\overline{k} }^0\theta ^{\overline{k}}.
\endaligned\tag{2.23}
$$
Then (2.22) implies
$$
\aligned
&f_{j0l}^0-f_{jl0}^0=f_{j\overline{k}}^0A_l^{\overline{k}}-f_k^0W_{jl}^k \\
&f_{j0\overline{l}}^0-f_{j\overline{l}0}^0=f_{jk}^0A_{\overline{l}
}^k+f_k^0W_{j\overline{l}}^k \\
&f_{jkl}^0-f_{jlk}^0=\sqrt{-1}f_k^0A_l^{\overline{j}}-\sqrt{-1}f_l^0A_k^{
\overline{j}} \\
&f_{jk\overline{l}}^0-f_{j\overline{l}k}^0-f_{j0}^0\delta
_k^l=f_t^0R_{jk\overline{l}}^t \\
&f_{j\overline{k}\overline{l}}^0-f_{j\overline{l}\overline{k}}^0=\sqrt{-1}
\delta_{\overline{k}}^{\overline{j}}f_t^0A_{\overline{l}}^t-\sqrt{-1}\delta
_{\overline{l}}^{\overline{j}}f_t^0A_{\overline{k}}^t.
\endaligned\tag{2.24}
$$
The commutative relations for $f_{\overline{j}AB}^0$ may be
derived similarly from the third equation in (2.13) or directly
from (2.24) by taking bar of each index.

By taking the exterior derivative of the first equation of (2.16),
we deduce that
$$
\aligned &Df_{00}^\alpha \wedge \theta +Df_{0j}^\alpha \wedge
\theta ^j+Df_{0\overline{j}}^\alpha \wedge \theta
^{\overline{j}}+if_{00}^\alpha \theta ^j\wedge
\theta ^{\overline{j}} \\
&+f_{0j}^\alpha A_{\overline{k}}^j\theta \wedge \theta
^{\overline{k}}+f_{0 \overline{j}}^\alpha A_k^{\overline{j}}\theta
\wedge \theta ^k=f_0^\beta \widehat{\Psi }_\beta
^\alpha\endaligned\tag{ 2.25}
$$
where
$$
\aligned Df_{00}^\alpha &=df_{00}^\alpha +f_{00}^\beta
\widehat{\theta }_\beta ^\alpha =f_{000}^\alpha \theta
+f_{00j}^\alpha \theta ^j+f_{00\overline{j}}^\alpha
\theta ^{\overline{j}} \\
Df_{0j}^\alpha &=df_{0j}^\alpha -f_{0k}^\alpha \theta
_j^k+f_{0j}^\beta \widehat{\theta }_\beta ^\alpha =f_{0j0}^\alpha
\theta +f_{0jl}^\alpha
\theta ^l+f_{0j\overline{l}}^\alpha \theta ^{\overline{l}} \\
Df_{0\overline{j}}^\alpha &=df_{0\overline{j}}^\alpha
-f_{0\overline{k} }^\alpha \theta
_{\overline{j}}^{\overline{k}}+f_{0\overline{j}}^\alpha
\widehat{\theta }_\beta ^\alpha =f_{0\overline{j}0}^\alpha \theta
+f_{0 \overline{j}l}^\alpha \theta
^l+f_{0\overline{j}\overline{l}}^\alpha \theta ^{\overline{l}}.
\endaligned\tag{2.26}
$$
Consequently
$$
\aligned &f_{00j}^\alpha -f_{0j0}^\alpha -f_{0\overline{k}}^\alpha
A_j^{\overline{k} }=f_0^\beta \widehat{R}_{\beta \gamma
\overline{\delta }}^\alpha (f_j^\gamma f_0^{\overline{\delta
}}-f_j^{\overline{\delta }}f_0^\gamma )+f_0^\beta
\widehat{W}_{\beta \gamma }^\alpha (f_j^\gamma
f_0^0-f_j^0f_0^\gamma ) \\
&-f_0^\beta \widehat{W}_{\beta \overline{\gamma } }^\alpha
(f_j^{\overline{\gamma }}f_0^0-f_j^0f_0^{\overline{\gamma
}})+if_0^\beta \widehat{A}_\delta ^{\overline{\beta }}(f_j^\alpha
f_0^\delta -f_j^\delta f_0^\alpha )-if_0^\beta
\widehat{A}_{\overline{\gamma }}^\alpha (f_j^{\overline{\gamma
}}f_0^{\overline{\beta }}-f_j^{\overline{\beta }}f_0^{
\overline{\gamma }}),
\endaligned\tag{2.27}
$$
$$
\aligned &f_{00\overline{j}}^\alpha -f_{0\overline{j}0}^\alpha
-f_{0k}^\alpha A_{ \overline{j}}^k=f_0^\beta \widehat{R}_{\beta
\gamma \overline{\delta } }^\alpha (f_{\overline{j}}^\gamma
f_0^{\overline{\delta }}-f_{\overline{j}}^{ \overline{\delta
}}f_0^\gamma )+f_0^\beta \widehat{W}_{\beta \gamma }^\alpha
(f_{\overline{j}}^\gamma f_0^0-f_{\overline{j}}^0f_0^\gamma
) \\
&-f_0^\beta \widehat{W}_{\beta \overline{ \gamma }}^\alpha
(f_{\overline{j}}^{\overline{\gamma }}f_0^0-f_{\overline{j}
}^0f_0^{\overline{\gamma }})+if_0^\beta \widehat{A}_\delta
^{\overline{\beta }}(f_{\overline{j}}^\alpha f_0^\delta
-f_{\overline{j}}^\delta f_0^\alpha )-if_0^\beta \widehat{A}_{
\overline{\gamma }}^\alpha (f_{\overline{j}}^{\overline{\gamma
}}f_0^{ \overline{\beta }}-f_{\overline{j}}^{\overline{\beta
}}f_0^{\overline{\gamma }}),
\endaligned \tag{2.28}
$$
$$
\aligned &f_{0jl}^\alpha -f_{0lj}^\alpha =f_0^\beta
\widehat{R}_{\beta \gamma \overline{\delta }}^\alpha (f_l^\gamma
f_j^{\overline{\delta }}-f_j^\gamma f_l^{\overline{\delta
}})+f_0^\beta \widehat{W}_{\beta \gamma }^\alpha (f_l^\gamma
f_j^0-f_j^\gamma f_l^0) \\
&-f_0^\beta \widehat{W}_{\beta \overline{\gamma }}^\alpha (f_l^{
\overline{\gamma }}f_j^0-f_j^{\overline{\gamma }}f_l^0)
+if_0^\beta \widehat{A}_\delta ^{\overline{\beta }}(f_l^\alpha
f_j^\delta -f_j^\alpha f_l^\delta )-if_0^\beta
\widehat{A}_{\overline{\gamma }}^\alpha (f_l^{\overline{\gamma
}}f_j^{\overline{\beta }}-f_j^{\overline{\gamma }
}f_l^{\overline{\beta }}),
\endaligned\tag{2.29}
$$
$$
\aligned &f_{0\overline{j}\overline{l}}^\alpha
-f_{0\overline{l}\overline{j}}^\alpha =f_0^\beta
\widehat{R}_{\beta \gamma \overline{\delta }}^\alpha (f_{
\overline{l}}^\gamma f_{\overline{j}}^{\overline{\delta
}}-f_{\overline{j} }^\gamma f_{\overline{l}}^{\overline{\delta
}})+f_0^\beta \widehat{W}_{\beta \gamma }^\alpha
(f_{\overline{l}}^\gamma f_{
\overline{j}}^0-f_{\overline{j}}^\gamma
f_{\overline{l}}^0) \\
&-f_0^\beta \widehat{W}_{\beta \overline{\gamma }}^\alpha
(f_{\overline{l}}^{\overline{ \gamma
}}f_{\overline{j}}^0-f_{\overline{j}}^{\overline{\gamma }}f_{
\overline{l}}^0) +if_0^\beta \widehat{A}_\delta ^{\overline{\beta
}}(f_{\overline{l}}^\alpha f_{\overline{j}}^\delta
-f_{\overline{j}}^\alpha f_{\overline{l}}^\delta )-if_0^\beta
\widehat{A}_{\overline{\gamma }}^\alpha (f_{\overline{l}}^{
\overline{\gamma }}f_{\overline{j}}^{\overline{\beta
}}-f_{\overline{j}}^{ \overline{\gamma
}}f_{\overline{l}}^{\overline{\beta }}),\endaligned \tag{2.30}
$$
$$
\aligned &f_{0\overline{j}l}^\alpha -f_{0l\overline{j}}^\alpha
+if_{00}^\alpha \delta _l^j=f_0^\beta \widehat{R}_{\beta \gamma
\overline{\delta }}^\alpha (f_l^\gamma
f_{\overline{j}}^{\overline{\delta }}-f_l^{\overline{\delta }}f_{
\overline{j}}^\gamma )+f_0^\beta \widehat{W}_{\beta \gamma
}^\alpha (f_l^\gamma f_{\overline{j}
}^0-f_{\overline{j}}^\gamma f_l^0) \\
&-f_0^\beta \widehat{W}_{\beta \overline{ \gamma }}^\alpha
(f_l^{\overline{\gamma }}f_{\overline{j}}^0-f_{\overline{j}
}^{\overline{\gamma }}f_l^0)+if_0^\beta \widehat{A}_\delta
^{\overline{\beta }}(f_l^\alpha f_{\overline{j }}^\delta
-f_{\overline{j}}^\alpha f_l^\delta )-if_0^\beta \widehat{A}_{
\overline{\gamma }}^\alpha (f_l^{\overline{\gamma
}}f_{\overline{j}}^{ \overline{\beta
}}-f_{\overline{j}}^{\overline{\gamma }}f_l^{\overline{\beta }}).
\endaligned\tag{2.31}
$$
Applying the exterior derivative to the second equation of (2.16),
we obtain
$$
\aligned &Df_{j0}^\alpha \wedge \theta +Df_{jl}^\alpha \wedge
\theta ^l+Df_{j\overline{l}}^\alpha \wedge \theta
^{\overline{l}}+if_{j0}^\alpha \theta
^k\wedge\theta ^{\overline{k}}\\
&+f_{jl}^\alpha A_{\overline{k}}^l\theta \wedge \theta
^{\overline{k}}+f_{j \overline{l}}^\alpha A_k^{\overline{l}}\theta
\wedge \theta ^k=-f_k^\alpha \Psi _j^k+f_j^\beta
\widehat{\Psi}_\beta ^\alpha\endaligned\tag{2.32}
$$
where
$$
\aligned &Df_{j0}^\alpha =df_{j0}^\alpha -f_{k0}^\alpha \theta
_j^k+f_{j0}^\beta \widehat{\theta }_\beta ^\alpha =f_{j00}^\alpha
\theta +f_{j0l}^\alpha
\theta ^l+f_{j0\overline{l}}^\alpha \theta ^{\overline{l}} \\
&Df_{jl}^\alpha =df_{jl}^\alpha -f_{kl}^\alpha \theta
_j^k-f_{jk}^\alpha \theta _l^k+f_{jl}^\beta \widehat{\theta
}_\beta ^\alpha =f_{jl0}^\alpha \theta +f_{jlk}^\alpha \theta
^k+f_{jl\overline{k}}^\alpha \theta ^{
\overline{k}} \\
&Df_{j\overline{l}}^\alpha =df_{j\overline{l}}^\alpha
-f_{k\overline{l} }^\alpha \theta _j^k-f_{j\overline{k}}^\alpha
\theta _{\overline{l}}^{ \overline{k}}+f_{j\overline{l}}^\beta
\widehat{\theta }_\beta ^\alpha =f_{j \overline{l}0}^\alpha \theta
+f_{j\overline{l}k}^\alpha \theta ^k+f_{j
\overline{l}\overline{k}}^\alpha \theta ^{\overline{k}}.
\endaligned\tag{2.33}
$$
Let us substitute (2.33) into (2.32) to get
$$
\aligned &f_{j0l}^\alpha -f_{jl0}^\alpha -f_{j\overline{k}}^\alpha
A_l^{\overline{k} }=-f_t^\alpha W_{jl}^t+f_j^\beta
\widehat{R}_{\beta \gamma \overline{\delta } }^\alpha (f_l^\gamma
f_0^{\overline{\delta }}-f_0^\gamma f_l^{\overline{ \delta
}})\\
&+f_j^\beta \widehat{W}_{\beta \gamma }^\alpha (f_l^\gamma
f_0^0-f_0^\gamma f_l^0) -f_j^\beta \widehat{W}_{\beta
\overline{\gamma }}^\alpha (f_l^{ \overline{\gamma
}}f_0^0-f_0^{\overline{\gamma }}f_l^0) \\
&+if_j^\beta \widehat{A}_\delta ^{\overline{\beta }}(f_l^\alpha
f_0^\delta -f_0^\alpha f_l^\delta )-if_j^\beta
\widehat{A}_{\overline{\gamma }}^\alpha (f_l^{\overline{\gamma
}}f_0^{\overline{\beta }}-f_0^{\overline{\gamma}
}f_l^{\overline{\beta }}), \endaligned \tag{2.34}
$$
$$
\aligned
&f_{j0\overline{l}}^\alpha -f_{j\overline{l}0}^\alpha
-f_{jk}^\alpha A_{ \overline{l}}^k=f_t^\alpha
W_{j\overline{l}}^t+f_j^\beta \widehat{R}_{\beta \gamma
\overline{\delta }}^\alpha (f_{\overline{l}}^\gamma
f_0^{\overline{\delta }}-f_0^\gamma f_{\overline{l}}^{\overline{\delta }}) \\
&+f_j^\beta\widehat{W}_{\beta\gamma}^\alpha(f_{\overline{l}}^\gamma
f_0^0-f_0^\gamma f_{\overline{l}}^0)-f_j^\beta \widehat{W}_{\beta
\overline{ \gamma }}^\alpha (f_{\overline{l}}^{\overline{\gamma
}}f_0^0-f_0^{\overline{\gamma }}f_{\overline{l}}^0) \\
&+if_j^\beta \widehat{A}_\delta ^{\overline{\beta
}}(f_{\overline{l}}^\alpha f_0^\delta -f_0^\alpha
f_{\overline{l}}^\delta )-if_j^\beta \widehat{A}_{
\overline{\gamma }}^\alpha (f_{\overline{l}}^{\overline{\gamma
}}f_0^{ \overline{\beta }}-f_0^{\overline{\gamma
}}f_{\overline{l}}^{\overline{\beta }}),
\endaligned\tag{2.35}
$$
$$
\aligned &f_{jlk}^\alpha -f_{jkl}^\alpha
=i(A_k^{\overline{j}}f_l^\alpha -f_k^\alpha
A_l^{\overline{j}})+f_j^\beta \widehat{R}_{\beta \gamma
\overline{\delta } }^\alpha (f_k^\gamma f_l^{\overline{\delta
}}-f_l^\gamma f_k^{\overline{
\delta }}) \\
&+f_j^\beta \widehat{W}_{\beta \gamma }^\alpha (f_k^\gamma
f_l^0-f_l^\gamma f_k^0)-f_j^\beta \widehat{W}_{\beta
\overline{\gamma }}^\alpha (f_k^{
\overline{\gamma }}f_l^0-f_l^{\overline{\gamma }}f_k^0) \\
&+if_j^\beta \widehat{A}_\delta ^{\overline{\beta }}(f_k^\alpha
f_l^\delta -f_l^\alpha f_k^\delta )-if_j^\beta
\widehat{A}_{\overline{\gamma }}^\alpha (f_k^{\overline{\gamma
}}f_l^{\overline{\beta }}-f_l^{\overline{\gamma }
}f_k^{\overline{\beta }}),
\endaligned\tag{2.36}
$$
$$
\aligned &f_{j\overline{l}k}^\alpha -f_{jk\overline{l}}^\alpha
+if_{j0}^\alpha \delta _{\overline{l}}^{\overline{k}}=-f_t^\alpha
R_{jk\overline{l}}^t+f_j^\beta \widehat{R}_{\beta \gamma
\overline{\delta }}^\alpha (f_k^\gamma f_{
\overline{l}}^{\overline{\delta }}-f_{\overline{l}}^\gamma
f_k^{\overline{
\delta }}) \\
&+f_j^\beta \widehat{W}_{\beta \gamma }^\alpha (f_k^\gamma
f_{\overline{l} }^0-f_{\overline{l}}^\gamma f_k^0)-f_j^\beta
\widehat{W}_{\beta \overline{ \gamma }}^\alpha
(f_k^{\overline{\gamma }}f_{\overline{l}}^0-f_{\overline{l}
}^{\overline{\gamma }}f_k^0) \\
&+if_j^\beta \widehat{A}_\delta ^{\overline{\beta }}(f_k^\alpha
f_{\overline{l }}^\delta -f_{\overline{l}}^\alpha f_k^\delta
)-if_j^\beta \widehat{A}_{ \overline{\gamma }}^\alpha
(f_k^{\overline{\gamma }}f_{\overline{l}}^{ \overline{\beta
}}-f_{\overline{l}}^{\overline{\gamma }}f_k^{\overline{\beta }}),
\endaligned\tag{2.37}
$$
$$
\aligned &f_{j\overline{l}\overline{k}}^\alpha
-f_{j\overline{k}\overline{l}}^\alpha =if_t^\alpha
(A_{\overline{k}}^t\delta _{\overline{l}}^{\overline{j}}-A_{
\overline{l}}^t\delta _{\overline{k}}^{\overline{j}})+f_j^\beta
\widehat{R} _{\beta \gamma \overline{\delta }}^\alpha
(f_{\overline{k}}^\gamma f_{ \overline{l}}^{\overline{\delta
}}-f_{\overline{l}}^\gamma f_{\overline{k}}^{
\overline{\delta }}) \\
&+f_j^\beta \widehat{W}_{\beta \gamma }^\alpha
(f_{\overline{k}}^\gamma f_{
\overline{l}}^0-f_{\overline{l}}^\gamma
f_{\overline{k}}^0)-f_j^\beta \widehat{W}_{\beta \overline{\gamma
}}^\alpha (f_{\overline{k}}^{\overline{ \gamma
}}f_{\overline{l}}^0-f_{\overline{l}}^{\overline{\gamma }}f_{
\overline{k}}^0) \\
&+if_j^\beta \widehat{A}_\delta ^{\overline{\beta
}}(f_{\overline{k}}^\alpha f_{\overline{l}}^\delta
-f_{\overline{l}}^\alpha f_{\overline{k}}^\delta )-if_j^\beta
\widehat{A}_{\overline{\gamma }}^\alpha (f_{\overline{k}}^{
\overline{\gamma }}f_{\overline{l}}^{\overline{\beta
}}-f_{\overline{l}}^{ \overline{\gamma
}}f_{\overline{k}}^{\overline{\beta }}). \endaligned\tag{2.38}
$$

Next, computing the exterior derivative of the third equation of
(2.16), we derive that
$$
\aligned &Df_{\overline{j}0}^\alpha \wedge \theta
+Df_{\overline{j}l}^\alpha \wedge \theta
^l+Df_{\overline{j}\overline{l}}^\alpha \wedge \theta
^{\overline{l}
}+if_{\overline{j}0}^\alpha \theta ^k\wedge \theta ^{\overline{k}} \\
&+f_{\overline{j}l}^\alpha A_{\overline{k}}^l\theta \wedge \theta
^{\overline{ k}}+f_{\overline{j}\overline{l}}^\alpha
A_k^{\overline{l}}\theta \wedge \theta ^k=-f_{\overline{k}}^\alpha
\Psi _{\overline{j}}^{\overline{k}}+f_{ \overline{j}}^\beta
\widehat{\Psi }_\beta ^\alpha
\endaligned\tag{2.39}
$$
where
$$
\aligned &Df_{\overline{j}0}^\alpha =df_{\overline{j}0}^\alpha
-f_{\overline{k} 0}^\alpha \theta
_{\overline{j}}^{\overline{k}}+f_{\overline{j}0}^\beta
\widehat{\theta }_\beta ^\alpha =f_{\overline{j}00}^\alpha \theta
+f_{ \overline{j}0l}^\alpha \theta
^l+f_{\overline{j}0\overline{l}}^\alpha \theta
^{\overline{l}} \\
&Df_{\overline{j}l}^\alpha =df_{\overline{j}l}^\alpha
-f_{\overline{k} l}^\alpha \theta
_{\overline{j}}^{\overline{k}}-f_{\overline{j}k}^\alpha \theta
_l^k+f_{\overline{j}l}^\beta \widehat{\theta }_\beta ^\alpha =f_{
\overline{j}l0}^\alpha \theta +f_{\overline{j}lk}^\alpha \theta
^k+f_{
\overline{j}l\overline{k}}^\alpha \theta ^{\overline{k}} \\
&Df_{\overline{j}\overline{l}}^\alpha
=df_{\overline{j}\overline{l}}^\alpha
-f_{\overline{k}\overline{l}}^\alpha \theta
_{\overline{j}}^{\overline{k}
}-f_{\overline{j}\overline{k}}^\alpha \theta
_{\overline{l}}^{\overline{k} }+f_{\overline{j}\overline{l}}^\beta
\widehat{\theta }_\beta ^\alpha =f_{
\overline{j}\overline{l}0}^\alpha \theta
+f_{\overline{j}\overline{l} k}^\alpha \theta
^k+f_{\overline{j}\overline{l}\overline{k}}^\alpha \theta ^{
\overline{k}}. \endaligned\tag{2.40}
$$
Therefore we have
$$
\aligned
&f_{\overline{j}0l}^\alpha -f_{\overline{j}l0}^\alpha
-f_{\overline{j} \overline{k}}^\alpha
A_l^{\overline{k}}=f_{\overline{t}}^\alpha W_{\overline{
j}l}^{\overline{t}}+f_{\overline{j}}^\beta \widehat{R}_{\beta
\gamma \overline{\delta }}^\alpha (f_l^\gamma
f_0^{\overline{\delta }}-f_0^\gamma
f_l^{\overline{\delta }}) \\
&+f_{\overline{j}}^\beta \widehat{W}_{\beta \gamma }^\alpha
(f_l^\gamma f_0^0-f_0^\gamma f_l^0)-f_{\overline{j}}^\beta
\widehat{W}_{\beta \overline{ \gamma }}^\alpha
(f_l^{\overline{\gamma }}f_0^0-f_0^{\overline{\gamma }
}f_l^0) \\
&+if_{\overline{j}}^\beta \widehat{A}_\delta ^{\overline{\beta
}}(f_l^\alpha f_0^\delta -f_0^\alpha f_l^\delta
)-if_{\overline{j}}^\beta \widehat{A}_{ \overline{\gamma }}^\alpha
(f_l^{\overline{\gamma }}f_0^{\overline{\beta }
}-f_0^{\overline{\gamma }}f_l^{\overline{\beta }}),
\endaligned\tag{2.41}
$$
$$\aligned
&f_{\overline{j}0\overline{l}}^\alpha
-f_{\overline{j}\overline{l}0}^\alpha -f_{\overline{j}k}^\alpha
A_{\overline{l}}^k=-f_{\overline{t}}^\alpha W_{
\overline{j}\overline{l}}^{\overline{t}}+f_{\overline{j}}^\beta
\widehat{R} _{\beta \gamma \overline{\delta }}^\alpha
(f_{\overline{l}}^\gamma f_0^{
\overline{\delta }}-f_0^\gamma f_{\overline{l}}^{\overline{\delta }}) \\
&+f_{\overline{j}}^\beta \widehat{W}_{\beta \gamma }^\alpha
(f_{\overline{l} }^\gamma f_0^0-f_0^\gamma
f_{\overline{l}}^0)-f_{\overline{j}}^\beta \widehat{W}_{\beta
\overline{\gamma }}^\alpha (f_{\overline{l}}^{\overline{
\gamma }}f_0^0-f_0^{\overline{\gamma }}f_{\overline{l}}^0) \\
&+if_{\overline{j}}^\beta \widehat{A}_\delta ^{\overline{\beta
}}(f_{ \overline{l}}^\alpha f_0^\delta -f_0^\alpha
f_{\overline{l}}^\delta )-if_{ \overline{j}}^\beta
\widehat{A}_{\overline{\gamma }}^\alpha (f_{\overline{l}
}^{\overline{\gamma }}f_0^{\overline{\beta
}}-f_0^{\overline{\gamma }}f_{ \overline{l}}^{\overline{\beta }}),
\endaligned \tag{2.42}
$$
$$
\aligned &f_{\overline{j}lk}^\alpha -f_{\overline{j}kl}^\alpha
=-if_{\overline{t} }^\alpha A_k^{\overline{t}}\delta
_l^j+if_{\overline{t}}^\alpha A_l^{ \overline{t}}\delta
_k^j+f_{\overline{j}}^\beta \widehat{R}_{\beta \gamma
\overline{\delta }}^\alpha (f_k^\gamma f_l^{\overline{\delta
}}-f_l^\gamma
f_k^{\overline{\delta }}) \\
&+f_{\overline{j}}^\beta \widehat{W}_{\beta \gamma }^\alpha
(f_k^\gamma f_l^0-f_l^\gamma f_k^0)-f_{\overline{j}}^\beta
\widehat{W}_{\beta \overline{ \gamma }}^\alpha
(f_k^{\overline{\gamma }}f_l^0-f_l^{\overline{\gamma }
}f_k^0) \\
&+if_{\overline{j}}^\beta \widehat{A}_\delta ^{\overline{\beta
}}(f_k^\alpha f_l^\delta -f_l^\alpha f_k^\delta
)-if_{\overline{j}}^\beta \widehat{A}_{ \overline{\gamma }}^\alpha
(f_k^{\overline{\gamma }}f_l^{\overline{\beta }
}-f_l^{\overline{\gamma }}f_k^{\overline{\beta }}),
\endaligned\tag{2.43}
$$
$$
\aligned &f_{\overline{j}\overline{l}k}^\alpha
-f_{\overline{j}k\overline{l}}^\alpha +if_{\overline{j}0}^\alpha
\delta _{\overline{l}}^{\overline{k}}=-f_{ \overline{t}}^\alpha
R_{\overline{j}k\overline{l}}^{\overline{t}}+f_{
\overline{j}}^\beta \widehat{R}_{\beta \gamma \overline{\delta
}}^\alpha (f_k^\gamma f_{\overline{l}}^{\overline{\delta
}}-f_{\overline{l}}^\gamma
f_k^{\overline{\delta }}) \\
&+f_{\overline{j}}^\beta \widehat{W}_{\beta \gamma }^\alpha
(f_k^\gamma f_{ \overline{l}}^0-f_{\overline{l}}^\gamma
f_k^0)-f_{\overline{j}}^\beta \widehat{W}_{\beta \overline{\gamma
}}^\alpha (f_k^{\overline{\gamma }}f_{
\overline{l}}^0-f_{\overline{l}}^{\overline{\gamma }}f_k^0) \\
&+if_{\overline{j}}^\beta \widehat{A}_\delta ^{\overline{\beta
}}(f_k^\alpha f_{\overline{l}}^\delta -f_{\overline{l}}^\alpha
f_k^\delta )-if_{\overline{j }}^\beta
\widehat{A}_{\overline{\gamma }}^\alpha (f_k^{\overline{\gamma
}}f_{ \overline{l}}^{\overline{\beta
}}-f_{\overline{l}}^{\overline{\gamma }}f_k^{ \overline{\beta }}),
\endaligned\tag{2.44}
$$
$$\aligned
&f_{\overline{j}\overline{l}\overline{k}}^\alpha
-f_{\overline{j}\overline{k} \overline{l}}^\alpha
=if_{\overline{k}}^\alpha A_{\overline{l}}^j-if_{
\overline{l}}^\alpha A_{\overline{k}}^j+f_{\overline{j}}^\beta
\widehat{R} _{\beta \gamma \overline{\delta }}^\alpha
(f_{\overline{k}}^\gamma f_{ \overline{l}}^{\overline{\delta
}}-f_{\overline{l}}^\gamma f_{\overline{k}}^{
\overline{\delta }}) \\
&+f_{\overline{j}}^\beta \widehat{W}_{\beta \gamma }^\alpha
(f_{\overline{k} }^\gamma
f_{\overline{l}}^0-f_{\overline{l}}^\gamma f_{\overline{k}}^0)-f_{
\overline{j}}^\beta \widehat{W}_{\beta \overline{\gamma }}^\alpha
(f_{ \overline{k}}^{\overline{\gamma
}}f_{\overline{l}}^0-f_{\overline{l}}^{
\overline{\gamma }}f_{\overline{k}}^0) \\
&+if_{\overline{j}}^\beta \widehat{A}_\delta ^{\overline{\beta
}}(f_{ \overline{k}}^\alpha f_{\overline{l}}^\delta
-f_{\overline{l}}^\alpha f_{ \overline{k}}^\delta
)-if_{\overline{j}}^\beta \widehat{A}_{\overline{\gamma }}^\alpha
(f_{\overline{k}}^{\overline{\gamma }}f_{\overline{l}}^{\overline{
\beta }}-f_{\overline{l}}^{\overline{\gamma
}}f_{\overline{k}}^{\overline{ \beta }}).
\endaligned\tag{2.45}
$$
Similarly the commutative relations of $f_{ABC}^{\overline{\alpha
}}$ can be deduced from (2.27)-(2.31), (2.34)-(2.38) and
(2.41)-(2.45) by taking bar of each index.

\heading{\bf 3. $(H,\widetilde{H})$-harmonic maps}
\endheading
\vskip 0.3 true cm

Let $(M^{2m+1},H(M),J,\theta )$ and
$(N^{2n+1},\widetilde{H}(N),\widetilde{J},\widetilde{ \theta })$
be two pseudo-Hermitian manifolds endowed with the Tanaka-Webster
connections $\nabla $ and $\widetilde{\nabla }$ respectively.
Suppose $\Psi $ is a section of $Hom(\otimes ^kTM,f^{-1}TN)$. Let
$\Psi _{H,\widetilde{H}}$ be a section of $Hom(\otimes
^kH(M),f^{-1} \widetilde{H}(N))$ defined by
$$
\Psi _{H,\widetilde{H}}(X_1,..,X_k)=\pi _{\widetilde{H}}\circ \Psi
(X_1,...,X_k)\tag{3.1}
$$
for any $X_1,...,X_k\in HM$, where
$\pi_{\widetilde{H}}:TN\rightarrow \widetilde{H}(N)$ is the
natural projection morphism. For convenience, one may extend $\Psi
_{H,\widetilde{H}}$ to a section of $Hom(\otimes ^kTM,f^{-1}TN)$
by requiring that $\Psi _{H,\widetilde{H}}(Z_1,..,Z_k)=0$ if one
of $Z_1,...,Z_k\in TM$ is vertical.

For any smooth map $f:M\rightarrow N$ between the pseudo-Hermitian
manifolds, Petit ([Pe]) introduced the following horizontal energy
functional
$$
E_{H,\widetilde{H}}(f)=\frac
12\int_M|df_{H,\widetilde{H}}|^2dv_{\theta } \tag{3.2}
$$
where $dv_{\theta}=\theta \wedge (d\theta )^m$ and then he derived
its first variational formula. Since our notations are slightly
different from those in [Pe], we will derive the first variational
formula of $E_{H,\widetilde{H}}$ again for the convenience of the
readers.

\proclaim{Proposition 3.1 ([Pe])} Let $\{f_t\}_{|t|<\varepsilon}$
be a family of maps with $f_0= f$ and $\frac{\partial
f_t}{\partial t}|_{t=0}=v\in \Gamma (f^{-1}TN)$. Then
$$ \frac d{dt}E_{H,\widetilde{H}}(f_t)|_{t=0}=-\int_M\langle v,\tau _{H,
\widetilde{H}}(f)-tr_{G_\theta
}(f^{*}\widetilde{A})_H\widetilde{\xi} \rangle dv_{\theta }
\tag{3.3}
$$
where $\tau _{H,\widetilde{H}}(f)$ is the horizontal tension field
given by
$$
\tau _{H,\widetilde{H}}(f)=tr_{G_\theta }\big(\beta
_{H,\widetilde{H}}+[(f^{*} \widetilde{\theta })\otimes
(f^{*}\widetilde{\tau })]_{H,\widetilde{H}}\big). \tag{3.4}
$$
\endproclaim
\demo{Proof} Let $\Phi :M\times (-\varepsilon ,\varepsilon
)\rightarrow N$ be the map defined by $\Phi (x,t)=f_t(x)$. A
vector $X\in T_xM$ may be identified with a vector $(X,0)\in
T_{(x,t)}(M\times (-\varepsilon ,\varepsilon))$. This
identification gives the following distribution $$
H_{(x,t)}(M\times (-\varepsilon ,\varepsilon ))=span\{(X,0):X\in
H(M)\}$$ on $M\times(-\varepsilon, \varepsilon)$. Set $d\Phi
_{H,\widetilde{H}}=\pi _{\widetilde{H}}(d\Phi |_{H(M\times
(-\varepsilon ,\varepsilon ))})$. When the meaning is clear, we
often write $(X,0)$ as $X$ for simplicity.

We denote by $\widetilde{\nabla }$ the pull-back connection in
$\Phi ^{-1}TN$ . Choose a local orthonormal frame field
$\{e_A\}_{A=1}^{2m}$ of $H(M)$. Applying Lemma 2.1 and (1.10), we
obtain that
$$
\aligned &\pi _{\widetilde{H}}\big(\widetilde{\nabla }_{\frac
\partial {\partial t}}d\Phi (e_A)-\widetilde{\nabla }_{e_A}d\Phi
(\frac
\partial {\partial t})\big)_{t=0}\\
&=\pi _{\widetilde{H}}\big(\widetilde{T}_{\widetilde{\nabla
}}(d\Phi (\frac \partial
{\partial t}),d\Phi (e_A))\big)_{t=0} \\
&=\pi _{\widetilde{H}}\big(\widetilde{\theta }(v)\widetilde{\tau
}(d\Phi (e_A))_{t=0}-\widetilde{\theta }(d\Phi
(e_A))_{t=0}\widetilde{\tau }(v)+d
\widetilde{\theta }(v,d\Phi (e_A)_{t=0})\widetilde{\xi}\big) \\
&=\widetilde{\theta }(v)\widetilde{\tau
}(df(e_A))-\widetilde{\theta } (df(e_A))\widetilde{\tau }(v)
\endaligned
$$
that is,
$$
\lbrack \widetilde{\nabla }_{\frac \partial {\partial t}}d\Phi
_{H, \widetilde{H}}(e_A)]_{t=0}=\pi
_{\widetilde{H}}(\widetilde{\nabla }_{e_A}v)+ \widetilde{\theta
}(v)\widetilde{\tau }(df(e_A))-\widetilde{\theta }(df(e_A))
\widetilde{\tau }(v).  \tag{3.5}
$$
By (3.5), we compute
$$
\aligned \frac 12&\frac d{dt}|df_{tH,\widetilde{H}}|^2\mid
_{t=0}=\sum_{A=1}^{2m}g_{ \widetilde{\theta
}}\big(\widetilde{\nabla }_{\frac \partial {\partial t}}d\Phi
_{H,\widetilde{H}}(e_A),d\Phi _{H,\widetilde{H}}(e_A)\big)_{t=0} \\
=&\sum_{A=1}^{2m}\{g_{\widetilde{\theta }}\big(\widetilde{\nabla
}_{e_A}v,df_{H, \widetilde{H}}(e_A)\big)+\widetilde{\theta
}(v)g_{\widetilde{\theta }}\big(\widetilde{\tau }(df(e_A)),df_{H,\widetilde{H}}(e_A)\big) \\
&-\widetilde{\theta }(df(e_A))g_{\widetilde{\theta
}}\big(\widetilde{\tau }(v),df_{H,\widetilde{H}}(e_A)\big)\} \\
=&\sum_{A=1}^{2m}\{ e_Ag_{\widetilde{\theta
}}\big(v,df_{H,\widetilde{H}}(e_A)\big)-g_{\widetilde{\theta }}\big(v,df_{H,\widetilde{H}}(\nabla _{e_A}e_A)\big)\}\\
&-\sum_{A=1}^{2m}\{g_{\widetilde{\theta
}}\big(v,(\widetilde{\nabla }_{e_A}df_{H,
\widetilde{H}})(e_A)\big)+(f^{*}\widetilde{\theta
})(e_A)g_{\widetilde{\theta }}\big(\widetilde{\tau }(df_{H,\widetilde{H}}(e_A)),v\big) \\
&-(f^{*}\widetilde{A})(e_A,e_A)g_{\widetilde{\theta}}(\widetilde{\xi},v)\}.
\endaligned\tag{3.6}
$$
Set $\alpha (X)=g_{\widetilde{\theta
}}(v,df_{H,\widetilde{H}}(X))$ for any $ X\in TM$. Using Lemma
1.3, we deduce that
$$
\aligned \delta(\alpha) &=-(\nabla _T\alpha
)(T)-\sum_{A=1}^{2m}(\nabla _{e_A}\alpha
)(e_A) \\
&=-\sum_{A=1}^{2m}[e_Ag_{\widetilde{\theta
}}(v,df_{H,\widetilde{H}}(e_A))-g_{ \widetilde{\theta
}}(v,df_{H,\widetilde{H}}(\nabla _{e_A}e_A))].
\endaligned\tag{3.7}
$$
It follows from (3.6), (3.7) and the divergence theorem that
$$
\aligned \frac
d{dt}E_{H,\widetilde{H}}(f_t)|_{t=0}=&-\sum_{A=1}^{2m}\int_M\{g_{
\widetilde{\theta }}\big(v,(\widetilde{\nabla
}_{e_A}df_{H,\widetilde{H} })(e_A)+(f^{*}\widetilde{\theta
})(e_A)\widetilde{\tau }(df_{H,\widetilde{H}
}(e_A))\big) \\
&-(f^{*}\widetilde{A})(e_A,e_A)g_{\widetilde{\theta
}}(\widetilde{\xi} ,v)\}dv_{\theta }. \endaligned\tag{3.8}
$$
Note that $\widetilde{\tau }$ is a $\widetilde{H}(N)$-valued
$1$-form and $ \widetilde{\tau }(\widetilde{\xi})=0$. Thus
$$
\aligned \sum_{A=1}^{2m}(f^{*}\widetilde{\theta
})(e_A)\widetilde{\tau }(df_{H, \widetilde{H}}(e_A))&=tr_{G_\theta
}(f^{*}\widetilde{\theta }\otimes f^{*}
\widetilde{\tau })_H \\
&=tr_{G_\theta }(f^{*}\widetilde{\theta }\otimes
f^{*}\widetilde{\tau } )_{H,\widetilde{H}}.
\endaligned\tag{3.9}
$$
Therefore (3.8) and (3.9) complete the proof of this
proposition.\qed
\enddemo
According to Proposition 3.1, $f:M\rightarrow N$ is a critical
point of $ E_{H,\widetilde{H}}$ if and only if
$$
\tau _{H,\widetilde{H}}(f)=0\ \text{ and \ }tr_{G_\theta
}(f^{*}\widetilde{A} )_H=0.  \tag{3.10}
$$
The critical point is referred to as a pseudoharmonic map in [Pe].
However, there is another kind of critical maps, which is also
called a pseudoharmonic map (see [BDU]). To avoid any possible
confusion, we modify Petit's definition slightly to introduce the
following
\definition {Definition 3.1}
A map $f:(M,H(M),J,\theta )\rightarrow
(N,\widetilde{H}(N),\widetilde{J}, \widetilde{\theta })$ is said
to be $(H,\widetilde{H})$-harmonic if $D_vE_{H,
\widetilde{H}}(f)=0$ for any $v\in \Gamma
(f^{-1}\widetilde{H}(N))$.
\enddefinition

By Proposition 3.1, we have
\proclaim{Corollary 3.2} Let
$f:(M,H(M),J,\theta )\rightarrow
(N,\widetilde{H}(N),\widetilde{J}, \widetilde{\theta })$ be a map.
Then $f$ is $(H,\widetilde{H})$-harmonic if and only if $ \tau
_{H,\widetilde{H}}(f)=0$, that is,
$$
tr_{G_\theta }\{\beta
_{H,\widetilde{H}}(f)+(f^{*}\widetilde{\theta } \otimes
f^{*}\widetilde{\tau })_{H,\widetilde{H}}\}=0. \tag{3.11}
$$
\endproclaim
\remark{Remark 3.1} We see from (3.10) that pseudoharmonic maps in
the sense of [Pe2] require an extra condition on the pull-back
pseudo-Hermitian torsion. If the target manifold is Sasakian, then
(3.3) implies that $D_vE_{H,\widetilde{H}}(f)=0$ automatically for
any vertical variation field $v$ along $f$. Consequently, Petit's
pseudoharmonic maps coincide with ours in this special case.
\endremark

In terms of the notations in \S 2, $\beta _{H,\widetilde{H}}$ and
$\big(f^{*}\widetilde{\theta }\otimes f^{*}\widetilde{\tau
}\big)_{H,\widetilde{H}}$ may be expressed as follows
$$
\aligned \beta _{H,\widetilde{H}}&=f_{ij}^\alpha \theta ^i\otimes
\theta ^j\otimes \widetilde{\eta }_\alpha
+f_{i\overline{j}}^\alpha \theta ^i\otimes \theta ^{
\overline{j}}\otimes \widetilde{\eta }_\alpha
+f_{\overline{i}j}^\alpha
\theta ^{\overline{i}}\otimes \theta ^j\otimes \widetilde{\eta }_\alpha \\
&+f_{\overline{i}\overline{j}}^\alpha \theta
^{\overline{i}}\otimes \theta ^{ \overline{j}}\otimes
\widetilde{\eta }_\alpha +f_{ij}^{\overline{\alpha } }\theta
^i\otimes \theta ^j\otimes \widetilde{\eta }_{\overline{\alpha
}}+f_{ \overline{i}j}^{\overline{\alpha }}\theta
^{\overline{i}}\otimes \theta
^j\otimes \widetilde{\eta }_{\overline{\alpha }} \\
&+f_{i\overline{j}}^{\overline{\alpha }}\theta ^i\otimes \theta
^{\overline{j} }\otimes \widetilde{\eta }_{\overline{\alpha
}}+f_{\overline{i}\overline{j} }^{\overline{\alpha }}\theta
\overline{^i}\otimes \theta ^{\overline{j} }\otimes
\widetilde{\eta }_{\overline{\alpha }}
\endaligned\tag{3.12}
$$
and
$$
\aligned \big(f^{*}\widetilde{\theta }\otimes f^{*}\widetilde{\tau
}\big)_{H, \widetilde{H}}&=\widehat{A}_{\overline{\beta }}^\alpha
(f_j^0f_k^{\overline{ \beta }}\theta ^j\otimes \theta
^k+f_j^0f_{\overline{k}}^{\overline{\beta }
}\theta ^j\otimes \theta ^{\overline{k}})\otimes \widetilde{\eta }_\alpha \\
&+\widehat{A}_{\overline{\beta }}^\alpha
(f_{\overline{j}}^0f_k^{\overline{ \beta }}\theta
^{\overline{j}}\otimes \theta ^k+f_{\overline{j}}^0f_{
\overline{k}}^{\overline{\beta }}\theta ^{\overline{j}}\otimes
\theta ^{\overline{k}})\otimes \widetilde{\eta }_\alpha \\
&+\widehat{A}_\beta ^{\overline{\alpha }}(f_j^0f_k^\beta \theta
^j\otimes \theta ^k+f_{\ j}^0f_{\overline{k}}^\beta \theta
^j\otimes \theta ^{\overline{k}})\otimes \widetilde{\eta }_{\overline{\alpha }} \\
&+\widehat{A}_\beta^{\overline{\alpha}}(f_{\overline{j}}^0f_k^\beta
\theta ^{\overline{j}}\otimes \theta ^k+f_{\
\overline{j}}^0f_{\overline{k}}^\beta \theta
^{\overline{j}}\otimes \theta ^{\overline{k}})\otimes
\widetilde{\eta }_{\overline{\alpha }}.
\endaligned\tag{3.13}
$$
Hence (3.11) is equivalent to
$$
f_{k\overline{k}}^\alpha +f_{\overline{k}k}^\alpha +\widehat{A}_{
\overline{\beta }}^\alpha f_k^0f_{\overline{k}}^{\overline{\beta
}}+\widehat{ A}_{\overline{\beta }}^\alpha
f_{\overline{k}}^0f_k^{\overline{\beta } }=0.\tag{3.14}
$$
Recall that a map $f:(M,H(M),J,\theta )\rightarrow
(N,\widetilde{H}(N), \widetilde{J},\widetilde{\theta })$ is said
to be horizontal if $ df(H(M))\subset\widetilde{H}(N)$ (cf.
[Pe2]), or equivalently, $f^{*}\widetilde{\theta }=u\theta $ for
some $u\in C^\infty (M)$. It follows from (3.14) that

\proclaim{Corollary 3.3} Let $f:(M,H(M),J,\theta )\rightarrow
(N,\widetilde{H}(N),\widetilde{J}, \widetilde{\theta })$ be a map.
Suppose that either $f$ is horizontal or $N$ is Sasakian. Then $f$
is $(H,\widetilde{H})$-harmonic if and only if
$$
f_{k\overline{k}}^\alpha +f_{\overline{k}k}^\alpha =0.
$$
\endproclaim

We introduce the following special kind of $(H,\widetilde{H}
)$-harmonic map, which will be an important tool for establishing
rigidity results in this paper.

\definition{Definition 3.2} Let $f:(M,H(M),J,\theta )\rightarrow
(N,\widetilde{H}(N),\widetilde{J}, \widetilde{\theta })$ be a map
between two pseudo-Hermitian manifolds. We say that $f$ is a
special $(H,\widetilde{H})$-harmonic map if it is a
$(H,\widetilde{H} )$-harmonic map with the following additional
property
$$f_{k\overline{k}}^0+f_{\overline{k}k}^0=0.\tag{3.15}$$
\enddefinition

Note that if $f$ is horizontal, then $f_k^0=f_{\overline{k}}^0=0$,
and thus $f_{k\overline{k}}^0+f_{\overline{k}k}^0=0$. As a result,
the map $f$ is a special $(H,\widetilde{H})$-harmonic map if it is
both horizontal and $(H,\widetilde{H} )$-harmonic. Nevertheless, a
$(H,\widetilde{H})$-harmonic map is not necessarily horizontal
(see Example 5.2). We will see that the special condition (3.15)
can not only enhance the regularity of a
$(H,\widetilde{H})$-harmonic map, but also remove superfluous data
for parameterizing all foliated $(H,\widetilde{H}
)$-biholomorphisms between two pseudo-Hermitian manifolds.

Let $(N,g)$ be a Riemannian manifold and let $\nabla ^g$ denote
its Levi-Civita connection. For a map $f:(M,H(M),J,\theta
)\rightarrow (N,g)$ from a compact pseudo-Hermitian manifold to
the Riemannian manifold $(N,g)$, we may define a horizontal energy
for $f$ by
$$
E_H(f)=\frac 12\int_M\sum_{A=1}^{2m}\langle df(e_A),df(e_A)\rangle
dv_\theta\tag{3.16}
$$
where $\{e_A\}_{A=1}^{2m}$ is any orthonormal basis in $H(M)$.
According to [BDU], a critical map of the energy $E_H$ is called
pseudoharmonic. Let us define the following second fundamental
form (with respect to the data $(\nabla ,\nabla ^g)$)
$$
\beta ^g(f)(X,Y)=\nabla _Y^gdf(X)-df(\nabla _YX)\tag{3.17}$$ for
$X,Y\in TM$. Set
$$
\tau_H ^g(f)=tr_{G_\theta }\beta ^g(f)=\sum_{A=1}^{2m}\beta
^g(f)(e_A,e_A).\tag{3.18}
$$
For any variation $f_t$ of $f$, we have (cf. [BDU], [DT])
$$
\frac d{dt}E_H(f_t)|_{t=0}=-\int_M\langle v,\tau_H^g(f)\rangle
dv_\theta\tag{3.19}
$$
where $v=(\partial f_t/\partial t)|_{t=0}$. Hence $f$ is
pseudoharmonic if and only if $\tau_H^g(f)=0$. Set
$$
\tau ^g(f)=tr_{g_\theta }\beta ^g(f)=\sum_{A=0}^{2m}\beta
^g(f)(e_A,e_A).\tag{3.20}
$$
From (1.13), we see that $\tau ^g(f)$ is the usual tension field.
Thus $f:(M,g_\theta )\rightarrow (N,g)$ is harmonic if and only if
$\tau ^g(f)=0$ (cf. [EL]).

Let $f:(M,H(M),J,\theta )\rightarrow
(N,H(N),\widetilde{J},\widetilde{\theta })$ be a map between two
pseudo-Hermitian manifolds, and let $\nabla ^{\widetilde{\theta
}}$ denote the Levi-Civita connection of the Webster metric
$g_{\widetilde{\theta }}$. Using Lemma 1.2, we deduce that
$$
\aligned \tau^{g_{\widetilde{\theta
}}}_H(f)=&\sum_{A=1}^{2m}[\nabla _{e_A}^{\widetilde{\theta
}}df(e_A)-df(\nabla
_{e_A}e_A)]\\
=&\sum_{A=1}^{2m}[\widetilde{\nabla }_{e_A}df(e_A)-df(\nabla
_{e_A}e_A)]-\sum_{A=1}^{2m}\widetilde{A}(df(e_A),df(e_A))\widetilde{\xi}\\
&+\sum_{A=1}^{2m}\widetilde{\theta }(df(e_A))\widetilde{\tau
}(df(e_A))+\sum_{A=1}^{2m}\widetilde{\theta
}(df(e_A))\widetilde{J}df(e_A)\endaligned\tag{3.21}
$$
From (3.21), we immediately get

\proclaim{Proposition 3.4} Let $f:(M,H(M),J,\theta )\rightarrow
(N,H(N),\widetilde{J},\widetilde{\theta })$ be a horizontal map.
Suppose $N$ is Sasakian. Then $f$ is $(H,\widetilde{H} )$-harmonic
if and only if $f$ is pseudoharmonic.
\endproclaim

\definition{Definition 3.3}
A map $f:(M^{2m+1},H(M),J,\theta )\rightarrow
(N^{2n+1},\widetilde{H}(N), \widetilde{J},\widetilde{\theta })$
between two pseudo-Hermitian manifolds is called a foliated map if
it preserves the leaves of the Reeb foliations.
\enddefinition

Clearly a map $f:(M^{2m+1},H(M),J,\theta )\rightarrow
(N^{2n+1},\widetilde{H}(N), \widetilde{J},\widetilde{\theta })$ is
foliated if and only if $f_0^\alpha =0 $ for $\alpha=1,...m$ or
equivalently, $f_0^{\overline{\alpha }}=0$ for $\alpha=1,...m$.

\proclaim{Lemma 3.5} Suppose a map $f:(M^{2m+1},J,\theta
)\rightarrow $ $(N^{2n+1},\widetilde{J}, \widetilde{\theta })$ is
foliated and horizontal. Then $df(\xi)=\lambda \widetilde{\xi}$
and $f^{*}\widetilde{\theta }=\lambda \theta $ for some constant
$\lambda $.
\endproclaim
\demo{Proof} Since $f$ is both foliated and horizontal, there
exists a function $\lambda$ such that $df(\xi)=\lambda
\widetilde{\xi}$ and $f^{*}\widetilde{\theta }=\lambda \theta $.
Consequently,
$$
f_0^0=\lambda, \quad f_0^\alpha =f_0^{\overline{\alpha }}=0,
\tag{3.22}
$$
and
$$
f_j^0=f_{\overline{j}}^0=0.\tag{3.23}
$$
From the first and second equations in (2.14), (3.22) and (3.23),
we find
$$
f_{0j}^0=f_{0\overline{j}}^0=0,\tag{3.24}
$$
that is, $e_j(\lambda)=e_{\overline{j}}(\lambda)=0$ ($j=1,...,m$)
in view of the first equation of (2.13). It follows that $\lambda$
is constant.\qed
\enddemo

The following result gives some relationship between
$(H,\widetilde{H})$-harmonic maps and harmonic maps of
pseudo-Hermitian manifolds.

\proclaim{Proposition 3.6} Let $f:(M^{2m+1},H(M),J,\theta
)\rightarrow (N^{2n+1},\widetilde{H}(N),
\widetilde{J},\widetilde{\theta })$ be a foliated and horizontal
map between two pseudo-Hermitian manifolds. Suppose $N$ is
Sasakian. Then $f$ is $(H,\widetilde{H})$-harmonic if and only if
$f$ is harmonic.
\endproclaim
\demo{Proof} Choose a local orthonormal frame field
$\{\xi,e_A\}_{A=1,...,2m}$ on $M$. According to (1.13) and Lemma
3.5, we have
$$
\widetilde{\nabla }_\xi^\theta df(\xi)=\lambda ^2\widetilde{\nabla
}_{\widetilde{\xi}}^\theta \widetilde{\xi}=0.\tag{3.25}
$$
Under the assumptions that $f$ is horizontal and $N$ is Sasakian,
we apply Lemma 1.2 to deduce that
$$
\sum_{A=1}^{2m}\widetilde{\nabla }_{e_A}^\theta
df(e_A)=\sum_{A=1}^{2m} \widetilde{\nabla }_{e_A}df(e_A)
=\sum_{A=1}^{2m}\widetilde{\nabla}_{e_A}df_{H,\widetilde{H}}(e_A).
\tag{3.26}
$$
From (1.13), (3.25) and (3.26), we find that
$$
\tau ^{g_{\widetilde{\theta }}}(f)=\tau _{H,\widetilde{H}}(f)
$$
Therefore $f$ is $(H,\widetilde{H})$-harmonic if and only if $f$
is harmonic.\qed
\enddemo

\remark{Remark 3.2} From Lemma 3.5 and Proposition 3.6, we realize
that preserving both horizontal and vertical distributions is a
too restrictive condition for a map between two pseudo-Hermitian
manifolds (see also Example 5.2).
\endremark

We know that a pseudo-Hermitian manifold $N$ is a compact regular
Sasakian manifold if and only if the foliation of $N$ induces a
Riemannian submersion $\pi :(N,g_{\widetilde{\theta }})\rightarrow
(B,g_B)$ over a compact K\"ahler manifold $B$.

\proclaim{Proposition 3.7}Let $(M^{2m+1},H(M),J,\theta)$ be a
compact pseudo-Hermitian manifold and let $N$ be a Sasakian
manifold which can be realized as a Riemannian submersion $\pi
:N\rightarrow B$ over a K\"ahler manifold $B$. Suppose
$f:M\rightarrow N$ is a map from $M$ to $N$ and $\varphi =\pi
\circ f$. Then $E_{H,\widetilde{H}}(f)=E_H(\varphi )$ and $d\pi
(\tau _{H,\widetilde{H}}(f))=\tau^{g_B}_H(\varphi )$. In
particular, if $f$ is foliated, then
$E_{H,\widetilde{H}}(f)=E(\varphi )$ and $ d\pi (\tau
_{H,\widetilde{H}}(f))=\tau ^{g_B}(\varphi )$, where $E(\cdot)$
and $\tau^{g_B}(\cdot)$ denote the usual energy functional and
tension field for maps between the Riemannian manifolds
$(M,g_\theta )$ and $(B,g_B)$.
\endproclaim
\demo{Proof} Since $\pi :N\rightarrow B$ is a Riemannian
submersion, we have
$$
\aligned
E_{H,\widetilde{H}}(f)&=\frac12\int_M\sum_{A=1}^{2m}\langle
df_{H,\widetilde{H}}(e_A),df_{H,\widetilde{H}}(e_A)\rangle
dv_{\theta}\\
&=\frac 12\int_M\sum_{A=1}^{2m}\langle d\pi \circ
df(e_A),d\pi \circ df(e_A)\rangle dv_{\theta }\\
&=E_H(\varphi).\endaligned\tag{3.27}
$$
Let $f_t$ ($|t|<\varepsilon $) be any variation of $f$ with
$f_0=f$ and $v=\frac{\partial f_t}{\partial t}|_{t=0}$. Set
$\varphi _t=\pi (f_t)$ and $w=d\pi (v)$. Clearly $\varphi_t$ is a
variation of $\varphi $ with $w=\frac{\partial \varphi
_t}{\partial t}|_{t=0}$. Then (3.27) yields
$$
\frac d{dt}E_{H,\widetilde{H}}(f_t)|_{t=0}=\frac d{dt}E_H(\varphi
_t)|_{t=0}.\tag{3.28}
$$
Applying Proposition 3.1 and (3.19) to the left hand and right
hand sides of (3.28) respectively, we get
$$
\int_M\langle w,d\pi (\tau _{H,\widetilde{H}}(f))-\tau^{g_B}
_H(\varphi)\rangle dv_{\theta}=0.
$$
Since $w$ can be arbitrary vector field on $B$, we have $d\pi
(\tau _{H,\widetilde{H}}(f))=\tau^{g_B}_H(\varphi)$.

Now we assume that $f$ is foliated, that is, $df(\xi )=0$. Hence
$E_{H,\widetilde{H}}(f)=E(\varphi).$ For any vector field $w$ on
$B$, we may lift it as a basic vector field $v$ on $N$. Let $\psi
_t$ be the one parameter of transformations generated by $v$.
Since $\psi_t:N\rightarrow N$ is foliated for each $t$, we see
that $\{f_t=\psi _t\circ f\}$ is a foliated variation of $f$, that
is, each $f_t$ is a foliated map and $f|_{t=0}=f$. Then (3.27)
implies $E_{H,\widetilde{H}}(f_t)=E(\varphi_t)$, where
$\varphi_t=\pi\circ f_t$. Consequently
$$
\int_M\langle w,d\pi (\tau _{H,\widetilde{H}}(f))-\tau
^{g_B}(\varphi)\rangle dv_{\theta} =0\tag{3.29}
$$
since the gradient of the energy functional $E$ is $-\tau
^{g_B}(\varphi)$. Since $w$ is arbitrary, we deduce from (3.29)
that $d\pi (\tau _{H,\widetilde{H}}(f))=\tau ^{g_B}(\varphi)$.
\qed
\enddemo

\proclaim{Corollary 3.8} Let $M$, $N$, $B$,  $f$ and $\varphi $ be
as in Proposition 3.8. Then $f$ is $(H,\widetilde{H})$-harmonic if
and only if $\varphi $ is pseudoharmonic. In particular, if $f$ is
foliated, then $f$ is $(H,\widetilde{H})$-harmonic if and only if
$\varphi $ is harmonic.
\endproclaim

\heading{\bf 4. Bochner formula for the horizontal energy density
$e_{H,\widetilde{H}}$}
\endheading

In this section, we will derive the formula of $\bigtriangleup
_be_{H, \widetilde{H}}$ for a map $f$ between two pseudo-Hermitian
manifolds. According to the notations in \S 2,
$$
df_{H,\widetilde{H}}(\eta _j) =f_j^\alpha \widetilde{\eta }_\alpha
+f_j^{ \overline{\alpha }}\widetilde{\eta }_{\overline{\alpha }},
\quad df_{H,\widetilde{H}}(\eta _{\overline{j}})
=f_{\overline{j}}^\alpha \widetilde{\eta }_\alpha
+f_{\overline{j}}^{\overline{\alpha }}\widetilde{ \eta
}_{\overline{\alpha }}
$$
and thus
$$
e_{H,\widetilde{H}}=\frac 12|df_{H,\widetilde{H}}|^2 =f_j^\alpha
f_{\overline{j}}^{\overline{\alpha }}+f_j^{\overline{\alpha }}f_{
\overline{j}}^\alpha. \tag{4.1}
$$
The horizontal differential of $e_{H,\widetilde{H}}$ is given by
$$
\aligned &d_He_{H,\widetilde{H}}\\
&=(f_j^\alpha f_{\overline{j}}^{\overline{\alpha }
}+f_j^{\overline{\alpha }}f_{\overline{j}}^\alpha )_k\theta
^k+(f_j^\alpha f_{\overline{j}}^{\overline{\alpha
}}+f_j^{\overline{\alpha }}f_{\overline{j}
}^\alpha )_{\overline{k}}\theta ^{\overline{k}} \\
&=(f_{jk}^\alpha f_{\overline{j}}^{\overline{\alpha }}+f_j^\alpha
f_{ \overline{j}k}^{\overline{\alpha }}+f_{jk}^{\overline{\alpha
}}f_{\overline{j}}^\alpha +f_j^{\overline{\alpha
}}f_{\overline{j}k}^\alpha )\theta ^k +(f_{j\overline{k}}^\alpha
f_{\overline{j}}^{\overline{\alpha } }+f_j^\alpha
f_{\overline{j}\overline{k}}^{\overline{\alpha }}+f_{j\overline{
k}}^{\overline{\alpha }}f_{\overline{j}}^\alpha
+f_j^{\overline{\alpha }}f_{ \overline{j}\overline{k}}^\alpha
)\theta ^{\overline{k}}.
\endaligned\tag{4.2}
$$
Consequently
$$
\aligned\bigtriangleup _be_{H,\widetilde{H}}=&|\beta
_{H,\widetilde{H}}|^2+f_{ \overline{j}}^{\overline{\alpha
}}f_{jk\overline{k}}^\alpha +f_j^\alpha f_{
\overline{j}k\overline{k}}^{\overline{\alpha
}}+f_{\overline{j}}^\alpha f_{jk
\overline{k}}^{\overline{\alpha }} \\
&+f_j^{\overline{\alpha }}f_{\overline{j}k\overline{k}}^\alpha
+f_{\overline{j }}^{\overline{\alpha }}f_{j\overline{k}k}^\alpha
+f_j^\alpha f_{\overline{j} \overline{k}k}^{\overline{\alpha
}}+f_{\overline{j}}^\alpha f_{j\overline{k} k}^{\overline{\alpha
}}+f_j^{\overline{\alpha }}f_{\overline{j}\overline{k} k}^\alpha.
\endaligned \tag{4.3}
$$
Using (2.17), (2.37), (2.38), (2.43) and (2.44), we perform the
following computations
$$
\aligned f_{jk\overline{k}}^\alpha =&[f_{kj}^\alpha
+\widehat{A}_{\overline{\beta } }^\alpha
(f_k^0f_j^{\overline{\beta }}-f_k^{\overline{\beta }}f_j^0)]_{
\overline{k}} \\
=&f_{kj\overline{k}}^\alpha +\widehat{A}_{\overline{\beta
},\overline{k} }^\alpha (f_k^0f_j^{\overline{\beta
}}-f_k^{\overline{\beta }}f_j^0)+ \widehat{A}_{\overline{\beta
}}^\alpha (f_k^0f_j^{\overline{\beta }}-f_k^{
\overline{\beta }}f_j^0)_{\overline{k}} \\
=&f_{k\overline{k}j}^\alpha +if_{k0}^\alpha \delta
_{\overline{k}}^{\overline{ j}}+f_t^\alpha
R_{kj\overline{k}}^t-f_k^\beta \widehat{R}_{\beta \gamma
\overline{\delta }}^\alpha (f_j^\gamma
f_{\overline{k}}^{\overline{\delta }
}-f_{\overline{k}}^\gamma f_j^{\overline{\delta }}) \\
&-f_k^\beta \widehat{W}_{\beta \gamma }^\alpha (f_j^\gamma
f_{\overline{k} }^0-f_{\overline{k}}^\gamma f_j^0)+f_k^\beta
\widehat{W}_{\beta \overline{ \gamma }}^\alpha
(f_j^{\overline{\gamma }}f_{\overline{k}}^0-f_{\overline{k}
}^{\overline{\gamma }}f_j^0) \\
&-if_k^\beta \widehat{A}_\delta ^{\overline{\beta }}(f_j^\alpha
f_{\overline{k }}^\delta -f_{\overline{k}}^\alpha f_j^\delta
)+if_k^\beta \widehat{A}_{ \overline{\gamma }}^\alpha
(f_j^{\overline{\gamma }}f_{\overline{k}}^{ \overline{\beta
}}-f_{\overline{k}}^{\overline{\gamma }}f_j^{\overline{\beta
}}) \\
&+\widehat{A}_{\overline{\beta },\overline{k}}^\alpha
(f_k^0f_j^{\overline{ \beta }}-f_j^0f_k^{\overline{\beta
}})+\widehat{A}_{\overline{\beta } }^\alpha
(f_{k\overline{k}}^0f_j^{\overline{\beta
}}+f_k^0f_{j\overline{k}}^{ \overline{\beta
}}-f_{k\overline{k}}^{\overline{\beta }}f_j^0-f_k^{\overline{
\beta }}f_{j\overline{k}}^0),
\endaligned \tag{4.4}
$$
$$
\aligned f_{j\overline{k}k}^\alpha =&[f_{\overline{k}j}^\alpha
+if_0^\alpha \delta _j^k+\widehat{A}_{\overline{\beta }}^\alpha
(f_j^{\overline{\beta }}f_{
\overline{k}}^0-f_j^0f_{\overline{k}}^{\overline{\beta }})]_k \\
=&f_{\overline{k}kj}^\alpha -if_{\overline{t}}^\alpha
A_k^{\overline{t} }\delta _j^k+if_{\overline{t}}^\alpha
A_j^{\overline{t}}\delta _k^k+f_{ \overline{k}}^\beta
\widehat{R}_{\beta \gamma \overline{\delta }}^\alpha
(f_k^\gamma f_j^{\overline{\delta }}-f_j^\gamma f_k^{\overline{\delta }}) \\
&+f_{\overline{k}}^\beta \widehat{W}_{\beta \gamma }^\alpha
(f_k^\gamma f_j^0-f_j^\gamma f_k^0)-f_{\overline{k}}^\beta
\widehat{W}_{\beta \overline{ \gamma }}^\alpha
(f_k^{\overline{\gamma }}f_j^0-f_j^{\overline{\gamma }
}f_k^0) \\
&+if_{\overline{k}}^\beta \widehat{A}_\delta ^{\overline{\beta
}}(f_k^\alpha f_j^\delta -f_j^\alpha f_k^\delta
)-if_{\overline{k}}^\beta \widehat{A}_{ \overline{\gamma }}^\alpha
(f_k^{\overline{\gamma }}f_j^{\overline{\beta }
}-f_j^{\overline{\gamma }}f_k^{\overline{\beta }})+if_{0k}^\alpha
\delta _j^k
\\
&+\widehat{A}_{\overline{\beta },k}^\alpha (f_j^{\overline{\beta
}}f_{ \overline{k}}^0-f_j^0f_{\overline{k}}^{\overline{\beta
}})+\widehat{A}_{ \overline{\beta }}^\alpha
(f_{jk}^{\overline{\beta }}f_{\overline{k}}^0+f_j^{
\overline{\beta
}}f_{\overline{k}k}^0-f_{jk}^0f_{\overline{k}}^{\overline{ \beta
}}-f_j^0f_{\overline{k}k}^{\overline{\beta }}),
\endaligned\tag{4.5}
$$
$$
\aligned f_{\overline{j}k\overline{k}}^\alpha
=&[f_{k\overline{j}}^\alpha -if_0^\alpha \delta
_k^j+\widehat{A}_{\overline{\beta }}^\alpha
(f_k^0f_{\overline{j}}^{ \overline{\beta }}-f_k^{\overline{\beta
}}f_{\overline{j}}^0)]_{\overline{k}}
\\
=&f_{k\overline{k}\overline{j}}^\alpha +if_t^\alpha
(A_{\overline{k}}^t\delta
_{\overline{j}}^{\overline{k}}-A_{\overline{j}}^t\delta
_{\overline{k}}^{ \overline{k}})+f_k^\beta \widehat{R}_{\beta
\gamma \overline{\delta } }^\alpha (f_{\overline{k}}^\gamma
f_{\overline{j}}^{\overline{\delta }}-f_{
\overline{j}}^\gamma f_{\overline{k}}^{\overline{\delta }}) \\
&+f_k^\beta \widehat{W}_{\beta \gamma }^\alpha
(f_{\overline{k}}^\gamma f_{
\overline{j}}^0-f_{\overline{j}}^\gamma
f_{\overline{k}}^0)-f_k^\beta \widehat{W}_{\beta \overline{\gamma
}}^\alpha (f_{\overline{k}}^{\overline{ \gamma
}}f_{\overline{j}}^0-f_{\overline{j}}^{\overline{\gamma }}f_{
\overline{k}}^0) \\
&+if_k^\beta \widehat{A}_\delta ^{\overline{\beta
}}(f_{\overline{k}}^\alpha f_{\overline{j}}^\delta
-f_{\overline{j}}^\alpha f_{\overline{k}}^\delta )-if_k^\beta
\widehat{A}_{\overline{\gamma }}^\alpha (f_{\overline{k}}^{
\overline{\gamma }}f_{\overline{j}}^{\overline{\beta
}}-f_{\overline{j}}^{ \overline{\gamma
}}f_{\overline{k}}^{\overline{\beta }})-if_{0\overline{k}
}^\alpha \delta _k^j \\
&+\widehat{A}_{\overline{\beta },k}^\alpha
(f_k^0f_{\overline{j}}^{\overline{ \beta }}-f_k^{\overline{\beta
}}f_{\overline{j}}^0)+\widehat{A}_{\overline{ \beta }}^\alpha
(f_{k\overline{k}}^0f_{\overline{j}}^{\overline{\beta }
}+f_k^0f_{\overline{j}\overline{k}}^{\overline{\beta
}}-f_{k\overline{k}}^{ \overline{\beta
}}f_{\overline{j}}^0-f_k^{\overline{\beta }}f_{\overline{j}
\overline{k}}^0)
\endaligned\tag{4.6}
$$
and
$$
\aligned f_{\overline{j}\overline{k}k}^\alpha
=&[f_{\overline{k}\overline{j}}^\alpha +
\widehat{A}_{\overline{\beta }}^\alpha
(f_{\overline{k}}^0f_{\overline{j}}^{ \overline{\beta
}}-f_{\overline{j}}^0f_{\overline{k}}^{\overline{\beta }})]_k
\\
=&f_{\overline{k}k\overline{j}}^\alpha -if_{\overline{k}0}^\alpha
\delta _{ \overline{j}}^{\overline{k}}-f_{\overline{t}}^\alpha
R_{\overline{k}k
\overline{j}}^{\overline{t}}+f_{\overline{k}}^\beta
\widehat{R}_{\beta \gamma \overline{\delta }}^\alpha (f_k^\gamma
f_{\overline{j}}^{\overline{
\delta }}-f_{\overline{j}}^\gamma f_k^{\overline{\delta }}) \\
&+f_{\overline{k}}^\beta \widehat{W}_{\beta \gamma }^\alpha
(f_k^\gamma f_{ \overline{j}}^0-f_{\overline{j}}^\gamma
f_k^0)-f_{\overline{k}}^\beta \widehat{W}_{\beta \overline{\gamma
}}^\alpha (f_k^{\overline{\gamma }}f_{
\overline{j}}^0-f_{\overline{j}}^{\overline{\gamma }}f_k^0) \\
&+if_{\overline{k}}^\beta \widehat{A}_\delta ^{\overline{\beta
}}(f_k^\alpha f_{\overline{j}}^\delta -f_{\overline{j}}^\alpha
f_k^\delta )-if_{\overline{k }}^\beta
\widehat{A}_{\overline{\gamma }}^\alpha (f_k^{\overline{\gamma
}}f_{ \overline{j}}^{\overline{\beta
}}-f_{\overline{j}}^{\overline{\gamma }}f_k^{
\overline{\beta }}) \\
&+\widehat{A}_{\overline{\beta },k}^\alpha
(f_{\overline{k}}^0f_{\overline{j} }^{\overline{\beta
}}-f_{\overline{j}}^0f_{\overline{k}}^{\overline{\beta }
})+\widehat{A}_{\overline{\beta }}^\alpha
(f_{\overline{k}k}^0f_{\overline{j} }^{\overline{\beta
}}+f_{\overline{k}}^0f_{\overline{j}k}^{\overline{\beta }
}-f_{\overline{j}k}^0f_{\overline{k}}^{\overline{\beta
}}-f_{\overline{j} }^0f_{\overline{k}k}^{\overline{\beta }}).
\endaligned\tag{4.7}
$$
Clearly the conjugates of (4.4), (4.5), (4.6) and (4.7) yield the
expressions of $f_{jk\overline{k}}^{\overline{\alpha }}$,
$f_{j\overline{k} k}^{\overline{\alpha }}$,
$f_{\overline{j}k\overline{k}}^{\overline{\alpha } } $ and
$f_{\overline{j}\overline{k}k}^{\overline{\alpha }}$.

In the remaining of this section, we assume that N is Sasakian. It
follows from (4.3), (4.4), (4.5), (4.6) and (4.7) that
$$
\aligned \bigtriangleup_He_{H,\widetilde{H}}=&|\beta
_{H,\widetilde{H}}|^2+\langle \widetilde{\nabla }\tau
_{H,\widetilde{H}},df_{H,\widetilde{H}}\rangle
+2i(f_{\overline{j}}^{\overline{\alpha }}f_{j0}^\alpha -f_j^\alpha
f_{\overline{j}0}^{\overline{\alpha }}) \\
&+2i(f_{\overline{j}}^\alpha f_{j0}^{\overline{\alpha
}}-f_j^{\overline{ \alpha }}f_{\overline{j}0}^\alpha
)+2mi(f_{\overline{j}}^\alpha f_{\overline{ k}}^{\overline{\alpha
}}A_{jk}-f_j^\alpha f_k^{\overline{\alpha }}A_{
\overline{j}\overline{k}}) \\
&+f_t^\alpha f_{\overline{j}}^{\overline{\alpha
}}R_{kj\overline{k}}^t+f_{ \overline{t}}^{\overline{\alpha
}}f_j^\alpha R_{\overline{k}\overline{j}k}^{
\overline{t}}+f_t^{\overline{\alpha }}f_{\overline{j}}^\alpha
R_{kj\overline{ k}}^t+f_{\overline{t}}^\alpha
f_j^{\overline{\alpha }}R_{\overline{k}
\overline{j}k}^{\overline{t}} \\
&-f_{\overline{j}}^{\overline{\alpha }}f_k^\beta
\widehat{R}_{\beta \gamma \overline{\delta }}^\alpha (f_j^\gamma
f_{\overline{k}}^{\overline{\delta } }-f_{\overline{k}}^\gamma
f_j^{\overline{\delta }})-f_j^\alpha f_{\overline{k
}}^{\overline{\beta }}\widehat{R}_{\overline{\beta
}\overline{\gamma }\delta }^{\overline{\alpha
}}(f_{\overline{j}}^{\overline{\gamma }}f_k^\delta -f_k^{
\overline{\gamma }}f_{\overline{j}}^\delta ) \\
&-f_{\overline{j}}^{\overline{\alpha }}f_{\overline{k}}^\beta
\widehat{R} _{\beta \gamma \overline{\delta }}^\alpha (f_j^\gamma
f_k^{\overline{\delta } }-f_k^\gamma f_j^{\overline{\delta
}})-f_j^\alpha f_k^{\overline{\beta }}
\widehat{R}_{\overline{\beta }\overline{\gamma }\delta
}^{\overline{\alpha } }(f_{\overline{j}}^{\overline{\gamma
}}f_{\overline{k}}^\delta -f_{\overline{
k}}^{\overline{\gamma }}f_{\overline{j}}^\delta ) \\
&-f_j^{\overline{\alpha }}f_k^\beta \widehat{R}_{\beta \gamma
\overline{ \delta }}^\alpha (f_{\overline{j}}^\gamma
f_{\overline{k}}^{\overline{\delta }}-f_{\overline{k}}^\gamma
f_{\overline{j}}^{\overline{\delta }})-f_{ \overline{j}}^\alpha
f_{\overline{k}}^{\overline{\beta }}\widehat{R}_{ \overline{\beta
}\overline{\gamma }\delta }^{\overline{\alpha }}(f_j^{
\overline{\gamma }}f_k^\delta -f_k^{\overline{\gamma }}f_j^\delta ) \\
&-f_j^{\overline{\alpha }}f_{\overline{k}}^\beta
\widehat{R}_{\beta \gamma \overline{\delta }}^\alpha
(f_{\overline{j}}^\gamma f_k^{\overline{\delta } }-f_k^\gamma
f_{\overline{j}}^{\overline{\delta }})-f_{\overline{j}}^\alpha
f_k^{\overline{\beta }}\widehat{R}_{\overline{\beta
}\overline{\gamma } \delta }^\alpha (f_j^{\overline{\gamma
}}f_{\overline{k}}^\delta -f_{\overline{k}}^{\overline{\gamma
}}f_j^\delta )
\endaligned \tag{4.8}
$$
where
$$|\beta _{H,\widetilde{H}}|^2=2(f_{jk}^\alpha f_{\overline{j}\overline{
k}}^{\overline{\alpha }}+f_{j\overline{k}}^\alpha
f_{\overline{j}k}^{ \overline{\alpha }}+f_{jk}^{\overline{\alpha
}}f_{\overline{j}\overline{k}}^\alpha
+f_{j\overline{k}}^{\overline{\alpha }}f_{\overline{j}k}^\alpha ),
\tag{4.9}
$$
$$\aligned
\langle\widetilde{\nabla }\tau
_{H,\widetilde{H}},df_{H,\widetilde{H}}\rangle=&f_{
\overline{j}}^{\overline{\alpha }}(f_{k\overline{k}j}^\alpha
+f_{\overline{k} kj}^\alpha )+f_j^{\overline{\alpha
}}(f_{k\overline{k}\overline{j}}^\alpha
+f_{\overline{k}k\overline{j}}^\alpha ) \\
&+f_j^\alpha (f_{k\overline{k}\overline{j}}^{\overline{\alpha
}}+f_{\overline{ k}k\overline{j}}^{\overline{\alpha
}})+f_{\overline{j}}^\alpha (f_{k \overline{k}j}^{\overline{\alpha
}}+f_{\overline{k}kj}^{\overline{\alpha }}).
\endaligned\tag{4.10}
$$

The pseudo-Hermitian Ricci curvature is given by (cf. [We], [DTo])
$$
R_{j\overline{k}}=R_{jt\overline{k}}^t=R_{\overline{t}jt\overline{k}}
\tag{4.11}
$$
which has the property
$$
R_{\overline{j}k}=R_{k\overline{j}}. \tag{4.12}
$$
One may define the pseudo-Hermitian Ricci transformation
$Ric_H:HM^C\rightarrow HM^C$ by (cf. [Ta], [DTo, page 57])
$$
Ric_H(\eta _j)=R_{j\overline{k}}\eta _k, \quad Ric_H(\eta
_{\overline{j }})=R_{\overline{j}k}\eta _{\overline{k}}.
\tag{4.13}
$$
From (4.11), (4.12) and (4.13), we find
$$
\aligned &f_t^\alpha f_{\overline{j}}^{\overline{\alpha
}}R_{kj\overline{k}}^t+f_{ \overline{t}}^{\overline{\alpha
}}f_j^\alpha R_{\overline{k}\overline{j}k}^{
\overline{t}}+f_t^{\overline{\alpha }}f_{\overline{j}}^\alpha
R_{kj\overline{ k}}^t+f_{\overline{t}}^\alpha
f_j^{\overline{\alpha }}R_{\overline{k}
\overline{j}k}^{\overline{t}} \\
&=f_t^\alpha f_{\overline{j}}^{\overline{\alpha
}}R_{j\overline{t}}+f_{ \overline{t}}^{\overline{\alpha
}}f_j^\alpha R_{t\overline{j}}+f_t^{ \overline{\alpha
}}f_{\overline{j}}^\alpha R_{j\overline{t}}+f_{\overline{t}
}^\alpha f_j^{\overline{\alpha }}R_{t\overline{j}} \\
&=2\langle df_{H,\widetilde{H}}(Ric_H(\eta
_j)),df_{H,\widetilde{H}}(\eta _{ \overline{j}})\rangle.
\endaligned\tag{4.14}
$$
In terms of (1.20) and (1.22), the curvature terms of $N$
appearing in (4.8) can be expressed as
$$
\aligned -&f_{\overline{j}}^{\overline{\alpha }}f_k^\beta
\widehat{R}_{\beta \gamma \overline{\delta }}^\alpha (f_j^\gamma
f_{\overline{k}}^{\overline{\delta } }-f_{\overline{k}}^\gamma
f_j^{\overline{\delta }})-f_j^\alpha f_{\overline{k
}}^{\overline{\beta }}\widehat{R}_{\overline{\beta
}\overline{\gamma }\delta }^{\overline{\alpha
}}(f_{\overline{j}}^{\overline{\gamma }}f_k^\delta -f_k^{
\overline{\gamma }}f_{\overline{j}}^\delta ) \\
-&f_{\overline{j}}^{\overline{\alpha }}f_{\overline{k}}^\beta
\widehat{R} _{\beta \gamma \overline{\delta }}^\alpha (f_j^\gamma
f_k^{\overline{\delta } }-f_k^\gamma f_j^{\overline{\delta
}})-f_j^\alpha f_k^{\overline{\beta }}
\widehat{R}_{\overline{\beta }\overline{\gamma }\delta
}^{\overline{\alpha } }(f_{\overline{j}}^{\overline{\gamma
}}f_{\overline{k}}^\delta -f_{\overline{
k}}^{\overline{\gamma }}f_{\overline{j}}^\delta ) \\
-&f_j^{\overline{\alpha }}f_k^\beta \widehat{R}_{\beta \gamma
\overline{ \delta }}^\alpha (f_{\overline{j}}^\gamma
f_{\overline{k}}^{\overline{\delta }}-f_{\overline{k}}^\gamma
f_{\overline{j}}^{\overline{\delta }})-f_{ \overline{j}}^\alpha
f_{\overline{k}}^{\overline{\beta }}\widehat{R}_{ \overline{\beta
}\overline{\gamma }\delta }^{\overline{\alpha }}(f_j^{
\overline{\gamma }}f_k^\delta -f_k^{\overline{\gamma }}f_j^\delta )\\
-&f_j^{\overline{\alpha }}f_{\overline{k}}^\beta
\widehat{R}_{\beta \gamma \overline{\delta }}^\alpha
(f_{\overline{j}}^\gamma f_k^{\overline{\delta } }-f_k^\gamma
f_{\overline{j}}^{\overline{\delta }})-f_{\overline{j}}^\alpha
f_k^{\overline{\beta }}\widehat{R}_{\overline{\beta
}\overline{\gamma } \delta }^\alpha (f_j^{\overline{\gamma
}}f_{\overline{k}}^\delta -f_{
\overline{k}}^{\overline{\gamma }}f_j^\delta ) \\
=&-\widetilde{R}\big(df_{H,\widetilde{H}}(\eta
_j),df_{H,\widetilde{H}}(\eta _k),df_{H,\widetilde{H}}(\eta
_{\overline{j}}),df_{H,\widetilde{H}}(\eta _{
\overline{k}})\big) \\
&-\widetilde{R}\big(df_{H,\widetilde{H}}(\eta
_{\overline{j}}),df_{H,\widetilde{H} }(\eta
_{\overline{k}}),df_{H,\widetilde{H}}(\eta _j),df_{H,\widetilde{H}
}(\eta _k)\big) \\
&-\widetilde{R}\big(df_{H,\widetilde{H}}(\eta
_j),df_{H,\widetilde{H}}(\eta _{
\overline{k}}),df_{H,\widetilde{H}}(\eta
_{\overline{j}}),df_{H,\widetilde{H}
}(\eta _k)\big) \\
&-\widetilde{R}\big(df_{H,\widetilde{H}}(\eta
_{\overline{j}}),df_{H,\widetilde{H} }(\eta
_k),df_{H,\widetilde{H}}(\eta _j),df_{H,\widetilde{H}}(\eta _{
\overline{k}})\big) \\
=&-2\widetilde{R}\big(df_{H,\widetilde{H}}(\eta
_j),df_{H,\widetilde{H}}(\eta _k),df_{H,\widetilde{H}}(\eta
_{\overline{j}}),df_{H,\widetilde{H}}(\eta _{
\overline{k}})\big) \\
&-2\widetilde{R}\big(df_{H,\widetilde{H}}(\eta
_j),df_{H,\widetilde{H}}(\eta_{\overline{k}}),df_{H,\widetilde{H}}(\eta
_{\overline{j}}),df_{H,\widetilde{H} }(\eta _k)\big).
\endaligned\tag{4.15}
$$
From (4.8), (4.13), (4.14) and (4.15), we obtain
$$
\aligned \bigtriangleup_He_{H,\widetilde{H}}=&|\beta
_{H,\widetilde{H}}|^2+\langle \widetilde{\nabla }\tau
_{H,\widetilde{H}},df_{H,\widetilde{H}}\rangle
-2i(f_j^{\overline{\alpha }}f_{\overline{j}0}^\alpha +f_j^\alpha
f_{ \overline{j}0}^{\overline{\alpha }}-f_{\overline{j}}^\alpha
f_{j0}^{
\overline{\alpha }}-f_{\overline{j}}^{\overline{\alpha }}f_{j0}^\alpha ) \\
&+2mi(f_{\overline{j}}^\alpha f_{\overline{k}}^{\overline{\alpha }
}A_{jk}-f_j^\alpha f_k^{\overline{\alpha
}}A_{\overline{j}\overline{k} })+2\langle
df_{H,\widetilde{H}}(Ric_H(\eta _j)),df_{H,\widetilde{H}}(\eta _{
\overline{j}})\rangle \\
&-2\widetilde{R}(df_{H,\widetilde{H}}(\eta
_j),df_{H,\widetilde{H}}(\eta _k),df_{H,\widetilde{H}}(\eta
_{\overline{j}}),df_{H,\widetilde{H}}(\eta _{ \overline{k}}))\\
&-2\widetilde{R}(df_{H,\widetilde{H}}(\eta
_j),df_{H,\widetilde{H}}(\eta _{
\overline{k}}),df_{H,\widetilde{H}}(\eta
_{\overline{j}}),df_{H,\widetilde{H} }(\eta _k)).
\endaligned\tag{4.16}
$$
The main difficulty in applications of (4.16) comes from the mixed
term, that is, the third term on the right hand side of (4.16). It
is known that the CR Paneitz operator is a useful tool to deal
with such a term. The usual Paneitz operator is a fourth order
differential operator defined as the divergence of a third order
differential operator $P$ acting on functions. One property of the
Paneitz operator is its nonnegativity, which plays an important
role in some rigidity problems in pseudo-Hermitian geometry. We
will generalize the operator $P$ to a differential operator, still
denoted by $P$, acting on maps between pseudo-Hermitian manifolds
as follows.

\definition{Definition 4.1}
Let $f:(M^{2m+1},J,\theta )\rightarrow (N^{2n+1},\widetilde{J},
\widetilde{\theta })$ be a map between two pseudo-Hermitian
manifolds. The primitive Paneitz operator $P(f)$ is a third order
differential operator given by
$$
\aligned
P(f)&=(f_{\overline{k}kj}^\alpha
+imA_{jk}f_{\overline{k}}^\alpha )\theta ^j\otimes \widetilde{\eta
}_\alpha +(f_{\overline{k}kj}^{\overline{\alpha }
}+imA_{jk}f_{\overline{k}}^{\overline{\alpha }})\theta ^j\otimes
\widetilde{\eta }_{\overline{\alpha }} \\
&=P_j^\alpha (f)\theta ^j\otimes \widetilde{\eta }_\alpha
+P_j^{\overline{ \alpha }}(f)\theta ^j\otimes \widetilde{\eta
}_{\overline{\alpha }}\endaligned
$$
where
$$P_j^\alpha (f)=f_{\overline{k}kj}^\alpha +imA_{jk}f_{\overline{k}
}^\alpha, \quad
P_j^{\overline{\alpha}}(f)=f_{\overline{k}kj}^{\overline{\alpha
}}+imA_{jk}f_{\overline{k}}^{\overline{\alpha }}.
$$
We use $\overline{P}(f)$ to denote the conjugate of $P$, that is,
$$
\overline{P}(f)=\overline{P_j^\alpha (f)}\theta
^{\overline{j}}\otimes \widetilde{\eta }_{\overline{\alpha
}}+\overline{P_j^{\overline{\alpha }}(f)} \theta
\overline{^j}\otimes \widetilde{\eta }_\alpha.
$$
\enddefinition

Note that $N$ is assumed to be Sasakian in this section. Using the
first and second equations of (2.17), we deduce that
$$
\aligned i(&f_j^{\overline{\alpha }}f_{\overline{j}0}^\alpha
+f_j^\alpha f_{\overline{j }0}^{\overline{\alpha
}}-f_{\overline{j}}^\alpha f_{j0}^{\overline{\alpha }
}-f_{\overline{j}}^{\overline{\alpha }}f_{j0}^\alpha )
\\
=&i(f_j^{\overline{ \alpha }}f_{0\overline{j}}^\alpha +f_j^\alpha
f_{0\overline{j}}^{\overline{ \alpha }}-f_{\overline{j}}^\alpha
f_{0j}^{\overline{\alpha }}-f_{\overline{j} }^{\overline{\alpha
}}f_{0j}^\alpha ) -iA_{\overline{k}\overline{j}}(f_k^\alpha
f_j^{\overline{\alpha }}+f_k^{ \overline{\alpha }}f_j^\alpha
)\\
&+iA_{kj}(f_{\overline{k}}^\alpha f_{\overline{j
}}^{\overline{\alpha }}+f_{\overline{k}}^{\overline{\alpha
}}f_{\overline{j}}^\alpha )\\
=&i(f_j^{\overline{\alpha }}f_{0\overline{j}}^\alpha +f_j^\alpha
f_{0 \overline{j}}^{\overline{\alpha }}-f_{\overline{j}}^\alpha
f_{0j}^{\overline{ \alpha }}-f_{\overline{j}}^{\overline{\alpha
}}f_{0j}^\alpha )-2i(A_{ \overline{k}\overline{j}}f_k^\alpha
f_j^{\overline{\alpha }}-A_{kj}f_{ \overline{k}}^\alpha
f_{\overline{j}}^{\overline{\alpha }})
\endaligned \tag{4.17}
$$
The fourth equation in (2.17) yields that
$$(f_{k\overline{k}}^\alpha -f_{\overline{k}k}^\alpha
)=mif_0^\alpha ,\ (f_{k\overline{k}}^{\overline{\alpha
}}-f_{\overline{k}k}^{\overline{\alpha }
})=mif_0^{\overline{\alpha }}.\tag{4.18}
$$
It follows from (4.17)and (4.18) that
$$
\aligned  -&2i(f_j^{\overline{\alpha }}f_{\overline{j}0}^\alpha
+f_j^\alpha f_{ \overline{j}0}^{\overline{\alpha
}}-f_{\overline{j}}^\alpha f_{j0}^{ \overline{\alpha
}}-f_{\overline{j}}^{\overline{\alpha
}}f_{j0}^\alpha )\\
=&4i(A_{\overline{k}\overline{j}}f_k^\alpha f_j^{\overline{\alpha
} }-A_{kj}f_{\overline{k}}^\alpha
f_{\overline{j}}^{\overline{\alpha }}) +\frac
2m(f_j^{\overline{\alpha }}f_{\overline{k}k\overline{j}}^\alpha
+f_{ \overline{j}}^\alpha f_{k\overline{k}j}^{\overline{\alpha
}}+f_j^\alpha f_{ \overline{k}k\overline{j}}^{\overline{\alpha
}}+f_{\overline{j}}^{\overline{
\alpha }}f_{k\overline{k}j}^\alpha ) \\
&-\frac 2m(f_j^{\overline{\alpha
}}f_{k\overline{k}\overline{j}}^\alpha +f_{ \overline{j}}^\alpha
f_{\overline{k}kj}^{\overline{\alpha }}+f_j^\alpha f_{k
\overline{k}\overline{j}}^{\overline{\alpha
}}+f_{\overline{j}}^{\overline{ \alpha }}f_{\overline{k}kj}^\alpha
).\endaligned\tag{4.19}
$$
By the definitions of $P(f)$ and $df_{H,\widetilde{H}}$, one has
$$
\aligned \langle P(f)+\overline{P(f)},df_{H,\widetilde{H}}\rangle
=&(f_j^{\overline{ \alpha }}f_{k\overline{k}\overline{j}}^\alpha
+f_{\overline{j}}^\alpha f_{ \overline{k}kj}^{\overline{\alpha
}}+f_j^\alpha f_{k\overline{k}\overline{j} }^{\overline{\alpha
}}+f_{\overline{j}}^{\overline{\alpha }}f_{\overline{k}
kj}^\alpha ) \\
&-2mi(A_{\overline{k}\overline{j}}f_k^\alpha f_j^{\overline{\alpha
} }-A_{kj}f_{\overline{k}}^\alpha
f_{\overline{j}}^{\overline{\alpha }}).
\endaligned\tag{4.20}
$$
Then (4.19) and (4.20) imply that
$$
\aligned &-2i(f_j^{\overline{\alpha }}f_{\overline{j}0}^\alpha
+f_j^\alpha f_{ \overline{j}0}^{\overline{\alpha
}}-f_{\overline{j}}^\alpha f_{j0}^{ \overline{\alpha
}}-f_{\overline{j}}^{\overline{\alpha }}f_{j0}^\alpha )\\
&=-\frac 2m\langle
P(f)+\overline{P(f)},df_{H,\widetilde{H}}\rangle+\frac
2m(f_j^{\overline{\alpha }}f_{\overline{k}k\overline{j}}^\alpha
+f_{ \overline{j}}^\alpha f_{k\overline{k}j}^{\overline{\alpha
}}+f_j^\alpha f_{ \overline{k}k\overline{j}}^{\overline{\alpha
}}+f_{\overline{j}}^{\overline{ \alpha }}f_{k\overline{k}j}^\alpha
).\endaligned\tag{4.21}
$$
On the other hand, using (4.10) and (4.20), one has
$$
\aligned &\langle \widetilde{\nabla }\tau
_{H,\widetilde{H}},df_{H,\widetilde{H} }\rangle
=(f_j^{\overline{\alpha }}f_{\overline{k}k\overline{j}}^\alpha
+f_{ \overline{j}}^\alpha f_{k\overline{k}j}^{\overline{\alpha
}}+f_j^\alpha f_{ \overline{k}k\overline{j}}^{\overline{\alpha
}}+f_{\overline{j}}^{\overline{
\alpha }}f_{k\overline{k}j}^\alpha ) \\
&+\langle P(f)+\overline{P(f)},df_{H,\widetilde{H}}\rangle
+2mi(A_{\overline{k }\overline{j}}f_k^\alpha f_j^{\overline{\alpha
}}-A_{kj}f_{\overline{k} }^\alpha
f_{\overline{j}}^{\overline{\alpha }}). \endaligned \tag{4.22}
$$
It follows from (4.21) and (4.22) that
$$
\aligned &-2i(f_j^{\overline{\alpha }}f_{\overline{j}0}^\alpha
+f_j^\alpha f_{ \overline{j}0}^{\overline{\alpha
}}-f_{\overline{j}}^\alpha f_{j0}^{ \overline{\alpha
}}-f_{\overline{j}}^{\overline{\alpha }}f_{j0}^\alpha )=\frac
2m\langle \widetilde{\nabla }\tau _{H,\widetilde{H}},df_{H,
\widetilde{H}}\rangle \\
&-\frac 4m\langle P(f)+\overline{P(f)},df_{H,\widetilde{H}}\rangle
-4i(A_{ \overline{k}\overline{j}}f_k^\alpha f_j^{\overline{\alpha
}}-A_{kj}f_{ \overline{k}}^\alpha
f_{\overline{j}}^{\overline{\alpha }}).
\endaligned \tag{4.23}
$$
From (4.16), (4.23), we conclude that

\proclaim{Theorem 4.1} Let $f:(M^{2m+1},J,\theta )\rightarrow
(N^{2n+1},\widetilde{J},\widetilde{ \theta })$ be a map between
two pseudo-Hermitian manifolds. If $N$ is Sasakian, then
$$
\aligned \bigtriangleup _be_{H,\widetilde{H}}=&|\beta
_{H,\widetilde{H}}|^2+(1+\frac 2m)\langle \widetilde{\nabla }\tau
_{H,\widetilde{H}},df_{H,\widetilde{H}}\rangle -\frac 4m\langle
P(f)+\overline{P(f)},df_{H,\widetilde{H}}\rangle
\\
&-(2m+4)i\sum (f_j^\alpha f_k^{\overline{\alpha
}}A_{\overline{j}\overline{k} }-f_{\overline{j}}^\alpha
f_{\overline{k}}^{\overline{\alpha } }A_{jk})+2\langle
df_{H,\widetilde{H}}(Ric_H(\eta _j)),df_{H,\widetilde{H}
}(\eta _{\overline{j}})\rangle \\
&-2\widetilde{R}(df_{H,\widetilde{H}}(\eta
_j),df_{H,\widetilde{H}}(\eta _k),df_{H,\widetilde{H}}(\eta
_{\overline{j}}),df_{H,\widetilde{H}}(\eta _{
\overline{k}})) \\
&-2\widetilde{R}(df_{H,\widetilde{H}}(\eta
_j),df_{H,\widetilde{H}}(\eta _{
\overline{k}}),df_{H,\widetilde{H}}(\eta
_{\overline{j}}),df_{H,\widetilde{H} }(\eta _k))
\endaligned
$$
\endproclaim

To apply the above Bochner formula, we want to investigate the
sign of the integral $\int_M \langle
P(f)+\overline{P(f)},df_{H,\widetilde{H}}\rangle$.

\proclaim{Lemma 4.2} Let $f:(M^{2m+1},J,\theta )\rightarrow
(N^{2n+1},\widetilde{J},\widetilde{ \theta })$ be a map from a
compact pseudo-Hermitian manifold to a Sasakian manifold. Then
$$
\aligned &i\int_M(f_j^{\overline{\alpha }}f_{\overline{j}0}^\alpha
+f_j^\alpha f_{ \overline{j}0}^{\overline{\alpha
}}-f_{\overline{j}}^\alpha f_{j0}^{ \overline{\alpha
}}-f_{\overline{j}}^{\overline{\alpha }}f_{j0}^\alpha )dv_{\theta
} \\
&=2m\int_Mf_0^\alpha f_0^{\overline{\alpha }}dv_{\theta }
-2i\int_M(A_{\overline{k} \overline{j}}f_k^\alpha
f_j^{\overline{\alpha }}-A_{kj}f_{\overline{k} }^\alpha
f_{\overline{j}}^{\overline{\alpha }})dv_{\theta }. \endaligned
\tag{4.24}
$$
\endproclaim

\demo{Proof} We introduce a global $1$-form on $M$ as follows
$$
\psi=i\big((f_j^{\overline{\alpha }}f_0^\alpha +f_j^\alpha
f_0^{\overline{ \alpha }})\theta ^j-(f_{\overline{j}}^\alpha
f_0^{\overline{\alpha }}+f_{ \overline{j}}^{\overline{\alpha
}}f_0^\alpha )\theta ^{\overline{j}}\big).
$$
Using (2.17), we compute
$$
\aligned \delta\psi=&-i\big\{(f_j^{\overline{\alpha }}f_0^\alpha
+f_j^\alpha f_0^{ \overline{\alpha
}})_{\overline{j}}-(f_{\overline{j}}^\alpha f_0^{\overline{
\alpha }}+f_{\overline{j}}^{\overline{\alpha }}f_0^\alpha )_j\big\} \\
=&-i\big\{(f_{j\overline{j}}^{\overline{\alpha }}f_0^\alpha
+f_{j\overline{j} }^\alpha f_0^{\overline{\alpha
}})-(f_{\overline{j}j}^\alpha f_0^{\overline{ \alpha
}}+f_{\overline{j}j}^{\overline{\alpha }}f_0^\alpha )\big\}
+i\big\{(f_j^{\overline{\alpha }}f_{0\overline{j}}^\alpha
+f_j^\alpha f_{0 \overline{j}}^{\overline{\alpha
}})-(f_{\overline{j}}^\alpha f_{0j}^{
\overline{\alpha }}+f_{\overline{j}}^{\overline{\alpha }}f_{0j}^\alpha )\big\} \\
=&-i\big\{(f_{j\overline{j}}^{\overline{\alpha
}}-f_{\overline{j}j}^{\overline{ \alpha }})f_0^\alpha
+(f_{j\overline{j}}^\alpha -f_{\overline{j}j}^\alpha
)f_0^{\overline{\alpha }}\big\}+i\big\{(f_j^{\overline{\alpha
}}f_{0\overline{j}}^\alpha +f_j^\alpha f_{0
\overline{j}}^{\overline{\alpha }})-(f_{\overline{j}}^\alpha
f_{0j}^{
\overline{\alpha }}+f_{\overline{j}}^{\overline{\alpha }}f_{0j}^\alpha )\big\} \\
=&2mf_0^\alpha f_0^{\overline{\alpha
}}-i\big\{(f_j^{\overline{\alpha }}f_{0 \overline{j}}^\alpha
+f_j^\alpha f_{0\overline{j}}^{\overline{\alpha }})-(f_{
\overline{j}}^\alpha f_{0j}^{\overline{\alpha
}}+f_{\overline{j}}^{\overline{ \alpha }}f_{0j}^\alpha )\big\}.
\endaligned
$$
The divergence theorem implies that
$$
i\int_M\big\{(f_j^{\overline{\alpha }}f_{0\overline{j}}^\alpha
+f_j^\alpha f_{0 \overline{j}}^{\overline{\alpha
}})-(f_{\overline{j}}^\alpha f_{0j}^{ \overline{\alpha
}}+f_{\overline{j}}^{\overline{\alpha }}f_{0j}^\alpha
)\big\}dv_{\theta } =2m\int_Mf_0^\alpha f_0^{\overline{\alpha
}}dv_{\theta } . \tag{4.25}
$$
Then (4.24) is a consequence of (4.17) and (4.25). \qed
\enddemo
\proclaim {Lemma 4.3} Let $f:(M^{2m+1},J,\theta )\rightarrow
(N^{2n+1},\widetilde{J},\widetilde{ \theta })$ be a map from a
compact pseudo-Hermitian manifold to a Sasakian manifold. Then we
have
$$
\aligned 2&\int_M\langle df_{H,\widetilde{H}}(Ric_H(\eta
_j)),df_{H,\widetilde{H} }(\eta _{\overline{j}})\rangle
dv_{\theta } \\
=&\int_M\big\{mi(f_j^{\overline{\alpha }}f_{ \overline{j}0}^\alpha
+f_j^\alpha f_{\overline{j}0}^{\overline{\alpha }}-f_{
\overline{j}}^\alpha f_{j0}^{\overline{\alpha
}}-f_{\overline{j}}^{\overline{ \alpha }}f_{j0}^\alpha )
-2(f_{jk}^\alpha f_{\overline{j}\overline{k}}^{\overline{\alpha
}}+f_{ \overline{j}\overline{k}}^\alpha f_{jk}^{\overline{\alpha
}}-f_{j\overline{k} }^\alpha f_{\overline{j}k}^{\overline{\alpha
}}-f_{\overline{j}k}^\alpha f_{j
\overline{k}}^{\overline{\alpha }})\big\}dv_{\theta }  \\
&+2\int_M\widetilde{R}(df_{H,\widetilde{H}}(\eta
_j),df_{H,\widetilde{H} }(\eta
_{\overline{k}}),df_{H,\widetilde{H}}(\eta _{\overline{j}}),df_{H,
\widetilde{H}}(\eta_k))dv_{\theta } \\
&-2\int_M\widetilde{R}(df_{H,\widetilde{H}}(\eta
_j),df_{H,\widetilde{H} }(\eta _k),df_{H,\widetilde{H}}(\eta
_{\overline{j}}),df_{H,\widetilde{H} }(\eta
_{\overline{k}}))dv_{\theta } .
\endaligned\tag{4.26}
$$
\endproclaim
\demo{Proof} Taking $k=l$ in (2.44) and (2.37) respectively and
summing over $k$ from $1$ to $m$, we get
$$
imf_j^{\overline{\alpha }}f_{\overline{j}0}^\alpha
=-f_j^{\overline{\alpha } }f_{\overline{j}\overline{k}k}^\alpha
+f_j^{\overline{\alpha }}f_{\overline{j }k\overline{k}}^\alpha
+f_j^{\overline{\alpha }}f_{\overline{t}}^\alpha R_{t
\overline{j}\overline{k}k}+\widetilde{R}_{\overline{\alpha }\beta
\gamma \overline{\delta }}f_j^{\overline{\alpha
}}f_{\overline{j}}^\beta f_k^\gamma
f_{\overline{k}}^{\overline{\delta
}}+\widetilde{R}_{\overline{\alpha }\beta \overline{\delta }\gamma
}f_j^{\overline{\alpha }}f_{\overline{j}}^\beta
f_k^{\overline{\delta }}f_{\overline{k}}^\gamma  \tag{4.27}
$$
and
$$
imf_j^\alpha f_{\overline{j}0}^{\overline{\alpha }}=-f_j^\alpha
f_{\overline{ j}\overline{k}k}^{\overline{\alpha }}+f_j^\alpha
f_{\overline{j}k\overline{k} }^{\overline{\alpha }}+f_j^\alpha
f_{\overline{t}}^{\overline{\alpha }}R_{t
\overline{j}\overline{k}k}+\widetilde{R}_{\alpha \overline{\beta
}\delta \overline{\gamma }}f_j^\alpha
f_{\overline{j}}^{\overline{\beta }}f_k^\delta
f_{\overline{k}}^{\overline{\gamma }}+\widetilde{R}_{\alpha
\overline{\beta } \overline{\gamma }\delta }f_j^\alpha
f_{\overline{j}}^{\overline{\beta } }f_k^{\overline{\gamma
}}f_{\overline{k}}^\delta .  \tag{4.28}
$$
Consequently
$$
\aligned &im(f_j^{\overline{\alpha }}f_{\overline{j}0}^\alpha
+f_j^\alpha f_{\overline{ j}0}^{\overline{\alpha
}}-f_{\overline{j}}^\alpha f_{j0}^{\overline{\alpha }
}-f_{\overline{j}}^{\overline{\alpha
}}f_{j0}^\alpha )\\
&=-f_j^\alpha f_{\overline{j}\overline{k}k}^{\overline{\alpha
}}-f_j^{ \overline{\alpha }}f_{\overline{j}\overline{k}k}^\alpha
-f_{\overline{j} }^\alpha f_{jk\overline{k}}^{\overline{\alpha
}}-f_{\overline{j}}^{\overline{ \alpha }}f_{jk\overline{k}}^\alpha
+f_j^\alpha f_{\overline{j}k\overline{k} }^{\overline{\alpha
}}+f_j^{\overline{\alpha }}f_{\overline{j}k\overline{k} }^\alpha
+f_{\overline{j}}^\alpha f_{j\overline{k}k}^{\overline{\alpha
}}+f_{
\overline{j}}^{\overline{\alpha }}f_{j\overline{k}k}^\alpha \\
&+2\langle df_{H,\widetilde{H}}(Ric_H(\eta
_j)),df_{H,\widetilde{H}}(\eta _{ \overline{j}})\rangle
+2\widetilde{R}(df_{H,\widetilde{H}}(\eta _j),df_{H,
\widetilde{H}}(\eta _{\overline{j}}),df_{H,\widetilde{H}}(\eta
_k),df_{H, \widetilde{H}}(\eta _{\overline{k}})).
\endaligned\tag{4.29}
$$
By Bianchi identity, we find
$$
\aligned \widetilde{R}&(df_{H,\widetilde{H}}(\eta
_j),df_{H,\widetilde{H}}(\eta
_{\overline{j}}),df_{H,\widetilde{H}}(\eta
_k),df_{H,\widetilde{H}}(\eta _{\overline{k}})) \\
=&-\widetilde{R}(df_{H,\widetilde{H}}(\eta
_j),df_{H,\widetilde{H}}(\eta _{
\overline{k}}),df_{H,\widetilde{H}}(\eta
_{\overline{j}}),df_{H,\widetilde{H}
}(\eta _k)) \\
&+\widetilde{R}(df_{H,\widetilde{H}}(\eta
_j),df_{H,\widetilde{H}}(\eta _k),df_{H,\widetilde{H}}(\eta
_{\overline{j}}),df_{H,\widetilde{H}}(\eta _{ \overline{k}})).
\endaligned \tag{4.30}
$$
Thus we obtain (4.26) by applying integration by parts to (4.29)
and using (4.30).\qed
\enddemo

The following result generalizes the non-negativity of the usual
Paneitz operator.

\proclaim{Theorem 4.4} Let $f:(M^{2m+1},H(M),J,\theta )\rightarrow
(N^{2n+1},H(N),\widetilde{J},\widetilde{\theta })$ be a map from a
compact pseudo-Hermitian manifold with $m\geq 2$ to a Sasakian
manifold with strongly seminegative horizontal curvature. Then
$$
-\int_M\langle P(f)+\overline{P(f)},df_{H,\widetilde{H}}\rangle
dv_{\theta }\geq 0.
$$
In particular, if $M$ is also Sasakian, then the above integral is
always nonnegative for $m\geq 1$ and every target Sasakian
manifold without any curvature condition.
\endproclaim
\demo{Proof} First we assume that $m\geq 2$ and $N$ is a Sasakian
manifold with strongly seminegative horizontal curvature.
Integrating (4.16) and using Lemma 4.3, we discover
$$
\aligned &\int_M|\beta _{H,\widetilde{H}}|^2dv_{\theta }
-\int_M|\tau _{H,\widetilde{H} }|^2dv_{\theta }
+(m-2)i\int_M(f_j^{\overline{\alpha }}f_{\overline{j}0}^\alpha
+f_j^\alpha f_{\overline{j}0}^{\overline{\alpha
}}-f_{\overline{j}}^\alpha f_{j0}^{\overline{\alpha
}}-f_{\overline{j}}^{\overline{\alpha }
}f_{j0}^\alpha )dv_{\theta }  \\
&+2mi\int_M(A_{jk}f_{\overline{j}}^\alpha
f_{\overline{k}}^{\overline{\alpha }
}-A_{\overline{j}\overline{k}}f_j^\alpha f_k^{\overline{\alpha }
})dv_{\theta } -2\int_M(f_{jk}^\alpha
f_{\overline{j}\overline{k}}^{\overline{\alpha }
}+f_{\overline{j}\overline{k}}^\alpha f_{jk}^{\overline{\alpha
}}-f_{j \overline{k}}^\alpha f_{\overline{j}k}^{\overline{\alpha
}}-f_{\overline{j}
k}^\alpha f_{j\overline{k}}^{\overline{\alpha }})dv_{\theta }  \\
&-4\int_M\widetilde{R}(df_{H,\widetilde{H}}(\eta
_j),df_{H,\widetilde{H} }(\eta _k),df_{H,\widetilde{H}}(\eta
_{\overline{j}}),df_{H,\widetilde{H} }(\eta
_{\overline{k}}))dv_{\theta } =0
\endaligned\tag{4.31}
$$
It follows from (4.9) and (4.31) that
$$
\aligned &4\int_M(f_{j\overline{k}}^\alpha
f_{\overline{j}k}^{\overline{\alpha }}+f_{ \overline{j}k}^\alpha
f_{j\overline{k}}^{\overline{\alpha }})dv_{\theta } -\int_M|\tau
_{H,\widetilde{H}}|^2dv_{\theta }
+2mi\int_M(A_{jk}f_{\overline{j}}^\alpha f_{\overline{k}
}^{\overline{\alpha }}-A_{\overline{j}\overline{k}}f_j^\alpha
f_k^{\overline{
\alpha }})dv_{\theta } \\
&+(m-2)i\int_M(f_j^{\overline{\alpha }}f_{\overline{j}0}^\alpha
+f_j^\alpha f_{\overline{j}0}^{\overline{\alpha
}}-f_{\overline{j}}^\alpha f_{j0}^{
\overline{\alpha }}-f_{\overline{j}}^{\overline{\alpha }}f_{j0}^\alpha )dv_{\theta } \\
&-4\int_M\widetilde{R}(df_{H,\widetilde{H}}(\eta
_j),df_{H,\widetilde{H} }(\eta _k),df_{H,\widetilde{H}}(\eta
_{\overline{j}}),df_{H,\widetilde{H} }(\eta
_{\overline{k}}))dv_{\theta } =0.
\endaligned\tag{4.32}
$$
The integral of (4.23) yields
$$
\aligned &i\int_M(f_j^{\overline{\alpha }}f_{\overline{j}0}^\alpha
+f_j^\alpha f_{ \overline{j}0}^{\overline{\alpha
}}-f_{\overline{j}}^\alpha f_{j0}^{ \overline{\alpha
}}-f_{\overline{j}}^{\overline{\alpha }}f_{j0}^\alpha
)dv_{\theta } \\
&=\frac 1m\int_M\big\{|\tau _{H,\widetilde{H}}|^2
+2<P(f)+\overline{P(f)},df_{H,\widetilde{H}}>-2mi(A_{jk}f_{
\overline{j}}^\alpha f_{\overline{k}}^{\overline{\alpha
}}-A_{\overline{j} \overline{k}}f_j^\alpha f_k^{\overline{\alpha
}})\big\}dv_{\theta } .
\endaligned\tag{4.33}
$$
We multiply (4.33) by $(m-1)$ and minus (4.24) to get
$$
\aligned (m&-2)i\int_M(f_j^{\overline{\alpha
}}f_{\overline{j}0}^\alpha +f_j^\alpha f_{
\overline{j}0}^{\overline{\alpha }}-f_{\overline{j}}^\alpha
f_{j0}^{ \overline{\alpha }}-f_{\overline{j}}^{\overline{\alpha
}}f_{j0}^\alpha )dv_{\theta }\\
=&\frac{m-1}m\int_M|\tau _{H,\widetilde{H}}|^2dv_{\theta }
-2m\int_Mf_0^\alpha f_0^{ \overline{\alpha }} dv_{\theta }
+\frac{2(m-1)}m\int_M\langle
P(f)+\overline{P(f)},df_{H,\widetilde{H} }\rangle dv_{\theta }\\
&-2mi\int_M(A_{jk}f_{\overline{j}}^\alpha f_{\overline{k}}^{
\overline{\alpha }}-A_{\overline{j}\overline{k}}f_j^\alpha
f_k^{\overline{ \alpha }})dv_{\theta }
\endaligned
$$
that is,
$$
\aligned &\int_M\big\{(m-2)i(f_j^{\overline{\alpha
}}f_{\overline{j}0}^\alpha +f_j^\alpha f_{
\overline{j}0}^{\overline{\alpha }}-f_{\overline{j}}^\alpha
f_{j0}^{ \overline{\alpha }}-f_{\overline{j}}^{\overline{\alpha
}}f_{j0}^\alpha )+2mi(A_{jk}f_{\overline{j}}^\alpha
f_{\overline{k}}^{\overline{\alpha
}}-A_{\overline{j}\overline{k}}f_j^\alpha f_k^{\overline{\alpha }})\big\}dv_{\theta }  \\
&=\frac{m-1}m\int_M|\tau _{H,\widetilde{H}}|^2dv_{\theta }
-2m\int_Mf_0^\alpha f_0^{ \overline{\alpha }}dv_{\theta }
+\frac{2(m-1)}m\int_M\langle P(f)+\overline{P(f)},df_{H,
\widetilde{H}}\rangle dv_{\theta }.
\endaligned\tag{4.34}
$$
Substituting (4.34) into (4.32) gives
$$
\aligned &\int_M\big\{4(f_{j\overline{k}}^\alpha
f_{\overline{j}k}^{\overline{\alpha }}+f_{ \overline{j}k}^\alpha
f_{j\overline{k}}^{\overline{\alpha }})-\frac 1m|\tau
_{H,\widetilde{H}}|^2+\frac{2(m-1)}m\langle
P(f)+\overline{P(f)},df_{H,\widetilde{H}}\rangle\big\}dv_{\theta } \\
&-2m\int_Mf_0^\alpha f_0^{\overline{\alpha }}dv_{\theta }
-4\int_M\widetilde{R}(df_{H,\widetilde{H}}(\eta
_j),df_{H,\widetilde{H} }(\eta _k),df_{H,\widetilde{H}}(\eta
_{\overline{j}}),df_{H,\widetilde{H} }(\eta
_{\overline{k}}))dv_{\theta } =0.
\endaligned\tag{4.35}
$$
Using the Cauchy-Schwarz inequality, the parallelogram law and the
fourth equation of (2.17), we discover
$$
\aligned f_{j\overline{k}}^\alpha
f_{\overline{j}k}^{\overline{\alpha }}+f_{\overline{ j}k}^\alpha
f_{j\overline{k}}^{\overline{\alpha }}&\geq f_{k\overline{k}
}^\alpha f_{\overline{k}k}^{\overline{\alpha
}}+f_{\overline{k}k}^\alpha f_{k
\overline{k}}^{\overline{\alpha }} \\
&\geq \frac 1m\sum_\alpha \big(|\sum_kf_{k\overline{k}}^\alpha
|^2+|\sum_kf_{
\overline{k}k}^\alpha |^2\big) \\
&=\frac 1{2m}\sum_\alpha \big(|\sum_k(f_{k\overline{k}}^\alpha
+f_{\overline{k} k}^\alpha )|^2+|\sum_k(f_{k\overline{k}}^\alpha
-f_{\overline{k}k}^\alpha
)|^2\big) \\
&=\frac 1{4m}|\tau _{H,\widetilde{H}}|^2+\frac m2f_0^\alpha
f_0^{\overline{ \alpha }}
\endaligned\tag{4.36}
$$
where $|\tau _{H,\widetilde{H}}|^2=2\sum_\alpha
|\sum_k(f_{k\overline{k} }^\alpha +f_{\overline{k}k}^\alpha )|^2$.
It follows from (4.35), (4.36) and the curvature assumption on $N$
that
$$
\aligned &-\int_M\langle
P(f)+\overline{P(f)}, df_{H,\widetilde{H}}\rangle dv_{\theta } \\
&\geq-\frac{2m}{m-1}\int_M\widetilde{R}(df_{H,\widetilde{H}}(\eta
_j),df_{H, \widetilde{H}}(\eta _k),df_{H,\widetilde{H}}(\eta
_{\overline{j}}),df_{H,\widetilde{H}}(\eta _{\overline{k}}))dv_{\theta }  \\
&\geq 0.\endaligned
$$
The first claim is proved.

Assume now that $M$ is a Sasakian manifold with $m\geq 1$ and $N$
is an arbitrary Sasakian manifold . Using (4.10) and the
integration by parts, we get from (4.22) that
$$
\aligned -\int_M\langle
P(f)+\overline{P(f)},df_{H,\widetilde{H}}\rangle dv_{\theta }
&=-\int_M\left( f_{\overline{j}}^{\overline{\alpha
}}f_{k\overline{k} j}^\alpha +f_j^{\overline{\alpha
}}f_{\overline{k}k\overline{j}}^\alpha +f_j^\alpha
f_{\overline{k}k\overline{j}}^{\overline{\alpha }}+f_{\overline{j
}}^\alpha f_{k\overline{k}j}^{\overline{\alpha }}\right)dv_{\theta }   \\
&=\int_M\left( f_{\overline{j}j}^{\overline{\alpha
}}f_{k\overline{k} }^\alpha +f_{j\overline{j}}^{\overline{\alpha
}}f_{\overline{k}k}^\alpha +f_{j\overline{j}}^\alpha
f_{\overline{k}k}^{\overline{\alpha }}+f_{
\overline{j}j}^\alpha f_{k\overline{k}}^{\overline{\alpha }}\right)dv_{\theta }  \\
&\geq 0.
\endaligned
$$
This gives the second claim. \qed
\enddemo

\definition{Definition 4.2}
A map $f:(M,H(M),J,\theta )\rightarrow
(N,\widetilde{H}(N),\widetilde{J}, \widetilde{\theta })$ is called
horizontally constant if it maps the domain manifold into a single
leaf of the pseudo-Hermitian foliation on $N$.
\enddefinition

\proclaim{Lemma 4.5} Let $f:(M,H(M),J,\theta )\rightarrow
(N,\widetilde{H}(N),\widetilde{J}, \widetilde{\theta })\ $be a map
between two pseudo-Hermitian manifolds. Then $f$ is horizontally
constant if and only if $df_{H,\widetilde{H}}=0$.
\endproclaim
\demo{Proof} If $f$ is horizontally constant, then $df(X)$ is
tangent to the fiber of $N$ for any $X\in TM$. Clearly we have
$df_{H,\widetilde{H} }=0$.

Conversely, we assume that $df_{H,\widetilde{H}}=0$, which is
equivalent to $ f_j^\alpha =f_{\overline{j}}^\alpha
=f_{\overline{j}}^{\overline{\alpha }
}=f_{\overline{j}}^{\overline{\alpha }}=0. $ Then the fourth
equation in (2.17) yields $f_0^\alpha =f_0^{\overline{\alpha
}}=0$. Hence $ \pi _{\widetilde{H}}df(X)=0$ for any $X\in TM$.
Suppose that $f(p)=q$ and $ \widetilde{C}_q$ is the integral curve
of $\widetilde{\xi}$ passing through the point $q$. For any point
$p^{\prime }\in M$, let $c(t)$ be a smooth curve joining $p$ and
$p^{\prime }$. Obviously $f(c(t))$ is a smooth curve passing
through $p$ and $df(c^{\prime }(t))=\lambda (t)\widetilde{\xi}$
for some function $\lambda(t)$, which means that $f(c(t))$ is the
reparametrization of the integral curve of $\widetilde{\xi}$.
Therefore $f(c(t))\subset \widetilde{C}_q$. In particular,
$f(p^{\prime })\in \widetilde{C}_q$. Since $p^{\prime}$ is
arbitrary, we conclude that $f(M)\subset \widetilde{C}_q$ . \qed
\enddemo

It is easy to see that a horizontally constant map is foliated and
$(H,\widetilde{H})$-harmonic. The following result is also obvious
by intuition.

\proclaim{Lemma 4.6} Suppose $f:(M^{2m+1},J,\theta )\rightarrow $
$(N^{2n+1},\widetilde{J}, \widetilde{\theta })$ is a horizontal
map. If $f$ is horizontally constant, then $f$ is constant.
\endproclaim

\demo{Proof} Since $f$ is horizontal and horizontally constant, we
have $f_j^0=f_{ \overline{j}}^0=0$ and $f_j^\alpha
=f_{\overline{j}}^\alpha =f_{\overline{j} }^\alpha
=f_{\overline{j}}^{\overline{\alpha }}=0$. Hence $df\circ i_H=0$.
Then the fifth equation of (2.14) implies $f_0^0=0$. This shows
that $ df(\xi)=0 $, since $f$ is foliated. Therefore we conclude
that $f$ is constant. \qed
\enddemo

\proclaim{Lemma 4.7} If $N$ has non-positive horizontal sectional
curvature, then
$$
\aligned &\widetilde{R}(df_{H,\widetilde{H}}(\eta
_j),df_{H,\widetilde{H}}(\eta _k),df_{H,\widetilde{H}}(\eta
_{\overline{j}}),df_{H,\widetilde{H}}(\eta _{ \overline{k}}))\\
&+\widetilde{R}(df_{H,\widetilde{H}}(\eta
_j),df_{H,\widetilde{H}}(\eta _{
\overline{k}}),df_{H,\widetilde{H}}(\eta
_{\overline{j}}),df_{H,\widetilde{H} }(\eta _k))\leq 0.
\endaligned \tag{4.37}
$$
\endproclaim

\demo{Proof} Write $df_{H,\widetilde{H}}(\eta _j)=X_j+iY_j$
($j=1,...,m$). We find that
$$
\aligned
&\text{the l.h.s. of (4.37)} \\
&=\widetilde{R}(X_j,X_k,X_j,X_k)-\widetilde{R}(X_j,X_k,Y_j,Y_k)+\widetilde{R}
(X_j,Y_k,X_j,Y_k)+\widetilde{R}(X_j,Y_k,Y_j,X_k) \\
&+\widetilde{R}(Y_j,X_k,X_j,Y_k)+\widetilde{R}(Y_j,X_k,Y_j,X_k)-\widetilde{R}
(Y_j,Y_k,X_j,X_k)+\widetilde{R}(Y_j,Y_k,Y_j,Y_k) \\
&+\widetilde{R}(X_j,X_k,X_j,X_k)+\widetilde{R}(X_j,X_k,Y_j,Y_k)+\widetilde{R}
(X_j,Y_k,X_j,Y_k)-\widetilde{R}(X_j,Y_k,Y_j,X_k) \\
&-\widetilde{R}(Y_j,X_k,X_j,Y_k)+\widetilde{R}(Y_j,X_k,Y_j,X_k)+\widetilde{R}
(Y_j,Y_k,X_j,X_k)+\widetilde{R}(Y_j,Y_k,Y_j,Y_k) \\
&=2\{\widetilde{R}(X_j,X_k,X_j,X_k)+\widetilde{R}(X_j,Y_k,X_j,Y_k)+\widetilde{R
}(Y_j,X_k,Y_j,X_k)+\widetilde{R}(Y_j,Y_k,Y_j,Y_k)\}, \\
\endaligned\tag{4.38}
$$
which is nonpositive by the assumption that $\widetilde{K}^H\leq
0$. \qed
\enddemo

Now we want to give some consequences of the Bochner formula in
Theorem 4.1.

\proclaim{Theorem 4.8} Let $f:(M^{2m+1},J,\theta )\rightarrow $
$(N^{2n+1},\widetilde{J},\widetilde{ \theta })$ be a
$(H,\widetilde{H})$-harmonic map from a compact pseudo-Hermitian
manifold with CR dimension $m\geq 2$ to a Sasakian manifold with
strongly semi-negative horizontal curvature. Let $\sigma_0(x)$ be
the maximal eigenvalue of the symmetric matrix
$(|A_{jk}|_x)_{m\times m}$ at $x\in M$. Suppose that
$$
Ric_H-(m+2)\sigma_0L_\theta \geq 0,  \tag{4.39}
$$
where $L_\theta$ is the Levi-form (see Definition 1.2). Then

(i) $\beta _{H,\widetilde{H}}=0$;

(ii) If $Ric_H-(m+2)\sigma_0L_\theta >0$ at a point in $M$, then
$f$ is horizontally constant;

(iii) If $N$ has negative horizontal sectional curvature, then $f$
is either horizontally constant or of horizontal rank one.
\endproclaim

\demo {Proof} At each point, let $\lambda $ be the minimal
eigenvalue of the Hermitian matrix $(R_{j\overline{k}})$.
Therefore
$$
\aligned \langle df_{H,\widetilde{H}}(Ric_H(\eta
_j)),df_{H,\widetilde{H}}(\eta _{\overline{j }})\rangle
&=R_{j\overline{k}}f_k^\alpha f_{\overline{j}}^{\overline{\alpha
}}+f_k^{
\overline{\alpha }}f_{\overline{j}}^\alpha R_{j\overline{k}} \\
&\geq \lambda (f_k^\alpha f_{\overline{k}}^{\overline{\alpha
}}+f_k^{
\overline{\alpha }}f_k^\alpha ) \\
&=\lambda \sum_{\alpha ,k}\big(|f_k^\alpha
|^2+|f_k^{\overline{\alpha }}|^2\big).
\endaligned\tag{4.40}
$$
By the definition of $\sigma_0$, one has
$$
\aligned |i(f_j^\alpha f_k^{\overline{\alpha
}}A_{\overline{j}\overline{k}}-f_{ \overline{j}}^\alpha
f_{\overline{k}}^{\overline{\alpha }}A_{jk})|&\leq 2\tau
_0\sum_{\alpha ,k}|f_k^\alpha f_k^{\overline{\alpha }}| \\
&\leq \tau _0\sum_{\alpha ,k}\big(|f_k^\alpha
|^2+|f_k^{\overline{\alpha }}|\big).
\endaligned\tag{4.41}
$$
From Theorems 4.1, 4.4, (4.39), (4.40) and (4.41), we immediately
get (i). Clearly, $\beta _{H,\widetilde{H}}=0$ implies that
$e_{H,\widetilde{H}}=$const. Besides, we have
$$ \aligned
&\widetilde{R}(df_{H,\widetilde{H}}(\eta
_j),df_{H,\widetilde{H}}(\eta _k),df_{H,\widetilde{H}}(\eta
_{\overline{j}}),df_{H,\widetilde{H}}(\eta _{ \overline{k}}))\\
&+\widetilde{R}(df_{H,\widetilde{H}}(\eta
_j),df_{H,\widetilde{H}}(\eta _{
\overline{k}}),df_{H,\widetilde{H}}(\eta
_{\overline{j}}),df_{H,\widetilde{H} }(\eta _k))= 0.
\endaligned \tag{4.42}
$$

Next assume that $Ric_H-(m+2)\tau _0L_\theta >0$ at a point $p$ in
$M$. This additional condition clearly implies that $f_k^\alpha
(p)=f_k^{\overline{ \alpha }}(p)=0$. Therefore
$e_{H,\widetilde{H}}=0$, that is, $df_{H, \widetilde{H}}=0$. It
follows from Lemma 4.5 that $f$ is horizontally constant. This
proves (ii).

Now we consider the claim (iii). If $K^H<0$, then (4.42) implies
that $rank(df_{H,\widetilde{H}})$ is zero or one in view of
(4.38). Since $e_{H,\widetilde{H}}$ is constant, the rank is
constant. In the first case, $f$ is horizontally constant; and in
the second case, we say that $f$ is of horizontal rank one. \qed
\enddemo
\remark{Remark 4.1} \newline (a) On a pseudo-Hermitian manifold,
the interchange of two covariant derivatives with respect to the
Tanaka-Webster connection yields not only the curvature terms, but
also the pseudo-Hermitian torsion term. Hence it seems natural
that the conditions for Bochner-type results include both
ingredients. \newline (b) We have already known that a
horizontally constant map is foliated. So the maps in cases (ii)
and (iii) are foliated. Using the Sasakian assumption on the
target manifold, the fourth equation of (2.17) implies that $f$ is
also foliated for the case (i) of Theorem 4.8.
\endremark

The following two results show that if the domain manifold in
Theorem 4.8 is also Sasakian, then the condition $m\geq 2$ is not
necessary, and the curvature condition on the target manifold may
be slightly weakened.

\proclaim{Corollary 4.9} Let $f:(M^{2m+1},J,\theta )\rightarrow $
$(N^{2n+1},\widetilde{J},\widetilde{ \theta })$ be a
$(H,\widetilde{H})$-harmonic map from a compact Sasakian manifold
to a Sasakian manifold with non-positive horizontal sectional
curvature. Suppose $ Ric_H\geq 0$. Then

(i) $\beta _{H,\overline{H}}=0$;

(ii) If $Ric_H>0$ at a point $p$ in $M$, then $f$ is horizontally
constant;

(iii) If $N$ has negative horizontal sectional curvature, then $f$
is either horizontally constant or of horizontal rank one.
\endproclaim
\demo{Proof} Using Lemma 4.7 and the second claim in Theorem 4.4,
the remaining arguments are similar to that for Theorem 4.8.\qed
\enddemo

\proclaim{Corollary 4.10} Let $M$, $N$ and $f$ be as in Corollary
4.9. If $f$ is horizontal, then

(i) $\beta =0$ (This property is called totally geodesic);

(ii) If $Ric_H>0$ at a point, then $f$ is constant;

(iii) If N has negative horizontal sectional curvature, then $f$
is either constant or of horizontal rank one.
\endproclaim
\demo{Proof} From Corollary 4.9 and Remark 4.1 (b), we know that
$f$ is a foliated map with $\beta _{H,\widetilde{H}}=0$, and thus
Lemma 3.6 implies that
$$
df(\xi)=\lambda \widetilde{\xi},\quad f^{*}\widetilde{\theta
}=\lambda \theta\tag{4.43}
$$
for some constant $\lambda $. Clearly the first equation of (4.43)
yields that
$$
\beta (\xi,X)=0  \tag{4.44}
$$
for any $X\in TM$. Since $f$ is horizontal, we have
$f_{jl}^0=f_{j\overline{l}}^0=0$. Therefore we conclude that
$\beta=0$.

The results for cases (ii) and (iii) follow immediately from Lemma
4.6 and Corollary 4.9.\qed
\enddemo
\remark{Remark 4.2} Corollary 4.10 improves a similar theorem of
[Pe2] in two aspects. It not only slightly strengthens the
corresponding results, but also weakens the curvature condition
for the target manifold.
\endremark

\heading{\bf 5. $(H,\widetilde{H})$-pluriharmonic and
$(H,\widetilde{H})$-holomorphic maps}
\endheading

In this section, we first introduce two special kinds of
$(H,\widetilde{H})$-harmonic maps:
$(H,\widetilde{H})$-pluriharmonic maps and foliated
$(H,\widetilde{H})$-holomorphic maps. Secondly, we give a unique
continuation theorem which ensures that a
$(H,\widetilde{H})$-harmonic map must be
$(H,\widetilde{H})$-holomorphic on the whole manifold if it is
$(H,\widetilde{H})$-holomorphic on an open subset. Clearly a
similar unique continuation result holds true for
$(H,\widetilde{H})$-antiholomorphicity. As a result, we easily
deduce a unique continuation theorem for horizontally constant
maps.

The $(H,\widetilde{H})$-harmonicity equation (3.11) suggests us to
introduce the following
\definition{Definition 5.1}
A map $f:(M,H(M),J,\theta )\rightarrow
(N,\widetilde{H}(N),\widetilde{J},\widetilde{ \theta })$ between
two pseudo-Hermitian manifolds is called a
$(H,\widetilde{H})$-pluriharmonic map if it satisfies
$$
(\beta _{H,\widetilde{H}}+f^{*}\theta \otimes f^{*}\widetilde{\tau
} )^{(1,1)}=0,\tag{5.1}
$$
where the left hand side term of (5.1) denotes the restriction of
$\beta _{H,\widetilde{H}}+f^{*}\theta \otimes f^{*}\widetilde{\tau
}$ to $H^{1,1}(M)$. Here $H^{1,1}(M)$ denotes the $(1,1)$-part of
$H(M)^C\otimes H(M)^C$ with respect to the complex structure on
$H(M)$.
\enddefinition

Clearly (5.1) implies that $$tr_{G_\theta }(\beta
_{H,\widetilde{H} }+f^{*}\theta \otimes f^{*}\widetilde{\tau
})=0,$$ that is, a $(H,\widetilde{H})$-pluriharmonic map is
automatically $(H,\widetilde{H})$-harmonic. It follows from
(3.12), (3.13) and (5.1) that a map $f$ is
$(H,\widetilde{H})$-pluriharmonic if and only if
$$
f_{j\overline{k}}^\alpha +\widehat{A}_{\overline{\beta }}^\alpha
f_j^0f_{ \overline{k}}^{\overline{\beta }}=0  \tag{5.2}
$$
and
$$
f_{\overline{k}j}^\alpha +\widehat{A}_{\overline{\beta }}^\alpha
f_{ \overline{k}}^0f_j^{\overline{\beta }}=0  \tag{5.3}
$$
for $1\leq \alpha \leq n$, $1\leq j,k\leq m$  or equivalently,
$f_{j \overline{k}}^{\overline{\alpha }}+\widehat{A}_\beta
^{\overline{\alpha } }f_j^0f_{\overline{k}}^\beta =0$ and
$f_{\overline{k}j}^{\overline{\alpha }}+ \widehat{A}_\beta
^{\overline{\alpha }}f_{\overline{k}}^0f_j^\beta =0$ for $1\leq
\alpha \leq n$, $1\leq j,k\leq m$.

\proclaim{Proposition 5.1} Suppose that
$f:(M,H(M),\theta,J)\rightarrow
(N,\widetilde{H}(N),\widetilde{\theta},\widetilde{J})$ is a
$(H,\widetilde{H})$-pluriharmonic map. Then $f$ is a foliated
$(H,\widetilde{H})$-harmonic map.
\endproclaim
\demo{Proof} We have already shown that $f$ is
$(H,\widetilde{H})$-harmonic. From (2.17), (5.2), (5.3), one may
find that
$$
f_0^\alpha =f_0^{\overline{\alpha }}=0,
$$
that is, $f$ is foliated.\qed
\enddemo

\definition{Definition 5.2} A map $f:(M,H(M),\theta,J)\rightarrow
(N,\widetilde{H}(N),\widetilde{\theta},\widetilde{J})$ between two
pseudohermitian manifolds is called
$(H,\widetilde{H})$-holomorphic (resp.
$(H,\widetilde{H})$-antiholomorphic) if it satisfies
$$
df_{H,\widetilde{H}}\circ J=\widetilde{J}\circ
df_{H,\widetilde{H}} \text{ (resp. }df_{H,\widetilde{H}}\circ
J=-\widetilde{J}\circ df_{H, \widetilde{H}}) \tag{5.4}
$$
Furthermore, if $f$ is foliated, then it is called a foliated
$(H,\widetilde{H})$-holomorphic map (resp. a foliated
$(H,\widetilde{H})$-antiholomorphic map).
\enddefinition

\remark{Remark 5.1} Clearly the composition of two foliated (resp.
$(H,\widetilde{H})$-holomorphic) maps is still a foliated (resp.
$(H,\widetilde{H})$-holomorphic) map. Note that the foliation of a
pseudo-Hermitian manifold $M$ is not transversally holomorphic in
general, although there is a complex structure $J$ on its
horizontal distribution $H(M)$.
\endremark

Suppose $f:(M,H(M),J,\theta )\rightarrow
(N,\widetilde{H}(N),\widetilde{J}, \widetilde{\theta })$ is a
smooth map between two pseudo-Hermitian manifolds. The
complexification of $df_{H,\widetilde{H}}$ determines various
partial horizontal differentials by the compositions with the
inclusions of $H^{1,0}(M)$ and $H^{0,1}(M)$ in $H(M)^C$
respectively and the projections of $\widetilde{H}(N)^C$ on
$H^{1,0}(N)$ and $H^{0,1}(N)$ respectively. Thus we have the
following bundle morphisms (cf. [Si1], [Do])

$$
\aligned
&\partial f_{H,\widetilde{H}}:H^{1,0}(M)\rightarrow
\widetilde{H} ^{1,0}(N),\quad \overline{\partial
}f_{H,\widetilde{H}}:H^{0,1}(M)
\rightarrow \widetilde{H}^{1,0}(N), \\
&\partial \overline{f_{H,\widetilde{H}}}:H^{1,0}(M)\rightarrow
\widetilde{H} ^{0,1}(N),\quad \overline{\partial
f_{H,\widetilde{H}}}:H^{0,1}(M) \rightarrow
\widetilde{H}^{0,1}(N),
\endaligned \tag{5.5}
$$
which can be locally expressed as follows
$$
\aligned
\partial f_{H,\widetilde{H}}=f_j^\alpha \theta ^j\otimes \widetilde{\eta }
_\alpha ,\quad \overline{\partial
}f_{H,\widetilde{H}}=f_{\overline{j}}^\alpha
\theta ^{\overline{j}}\otimes \widetilde{\eta }_\alpha,\\
\partial
\overline{f_{H,\widetilde{H}}}=f_j^{\overline{\alpha }}\theta
^j\otimes \widetilde{\eta }_{\overline{\alpha }},\quad
\overline{\partial f_{H,\widetilde{H}}}=f_{\overline{j}}^{
\overline{\alpha }}\theta ^{\overline{j}}\otimes \widetilde{\eta
}_{ \overline{\alpha }}.
\endaligned\tag{5.6}
$$
From (5.4), it is clear to see that $f:M\rightarrow N$ is
$(H,\widetilde{H})$-holomorphic (resp.
$(H,\widetilde{H})$-antiholomorphic) if and only if
$\overline{\partial }f_{H, \widetilde{H}}=0$ (resp. $\partial
f_{H,\widetilde{H}}=0$).

\proclaim{Proposition 5.2} Suppose that $f:M\rightarrow N$ is
either $(H,\widetilde{H})$-holomorphic or
$(H,\widetilde{H})$-antiholomorphic. Then $f$ is
$(H,\widetilde{H})$-harmonic if and only if $f$ is foliated.
\endproclaim
\demo{Proof} Without loss of generality, we assume that $f$ is
$(H,\widetilde{H})$-holomorphic. Then the
$(H,\widetilde{H})$-holomorphicity of $f$ means that
$f_{\overline{k}}^\alpha =f_k^{\overline{\alpha }}=0$, and thus
$f_{\overline{k}j}^\alpha =0$. Consequently, the
$(H,\widetilde{H})$-harmonicity equation (3.14) becomes
$$
f_{k\overline{k}}^\alpha +\widehat{A}_{\overline{\beta }}^\alpha
f_k^0f_{\overline{k}}^{\overline{\beta }}=0.  \tag{5.7}
$$
On the other hand, the fourth equation of (2.17) yields that
$$
f_{k\overline{k}}^\alpha +\widehat{A}_{\overline{\beta }}^\alpha
f_k^0f_{ \overline{k}}^{\overline{\beta }}=mif_0^\alpha  \tag{5.8}
$$
Therefore we conclude from (5.7) and (5.8) that $f$ is
$(H,\widetilde{H})$-harmonic if and only if it is foliated. \qed
\enddemo

\proclaim{Theorem 5.3} Let $f:(M,H(M),\theta ,J)\rightarrow
(N,\widetilde{H}(N),\widetilde{\theta },\widetilde{J})$ be either
a foliated $(H,\widetilde{H})$-holomorphic map or a foliated
$(H,\widetilde{H})$-antiholomorphic map. Then $f$ is a
$(H,\widetilde{H})$-pluriharmonic map.
\endproclaim

\demo{Proof} Without loss of generality, we may assume that $f$ is
$(H,\widetilde{H})$-holomorphic. Then we get
$f_{\overline{j}}^\alpha =f_j^{\overline{ \alpha
}}=f_{\overline{j}l}^\alpha =0$ as in Proposition 5.2. It follows
that
$$
f_{\overline{j}l}^\alpha +\widehat{A}_{\overline{\beta }}^\alpha
f_{ \overline{j}}^0f_l^{\overline{\beta }}=0.  \tag{5.9}
$$
Since $f$ is foliated, we have $f_0^\alpha =0$. Consequently,
(2.17) and (5.9) give that
$$
f_{l\overline{j}}^\alpha +\widehat{A}_{\overline{\beta }}^\alpha
f_l^0f_{ \overline{j}}^{\overline{\beta }}=0.
$$
This proves the $(H,\widetilde{H})$-pluriharmonicity of $f$. \qed
\enddemo

Recall that a map $f:M\rightarrow N$ is called a CR map if $f$ is
horizontal and $df_{H,\widetilde{H}}\circ J=\widetilde{J}\circ
df_{H, \widetilde{H}}$ ([DTo]). Thus CR maps provide us many
examples of $(H,\widetilde{H})$-holomorphic maps. Since a CR map
is not necessarily a foliated map, it is not
$(H,\widetilde{H})$-harmonic in general. A map $f:M\rightarrow N$
between two pseudo-Hermitian manifolds is called a CR-holomorphic
map (cf. [Dr], [IP], [Ur]) if
$$
df\circ J=\widetilde{J}\circ df.  \tag{5.10}
$$
From [IP], we know that a CR-holomorphic map is a harmonic map
between the two Riemannian manifolds $(M,g_\theta )$ and
$(N,\widetilde{g}_{\widetilde{ \theta }})$ in the usual sense.
Clearly (5.10) implies that $\widetilde{J} \circ df(\xi)=0$ and
$df(H(M))\subset \widetilde{H}(N)$. Hence a CR-holomorphic map is
just a foliated CR map. In addition, a special kind of CR maps,
called pseudo-Hermitian immersions in [Dr], also provide us a lot
of examples for foliated $(H,\widetilde{H})$-holomorphic maps.

\definition{Definition 5.3}
We call a diffeomorphism $f:(M^{2m+1},J,\theta )\rightarrow
(N^{2m+1}, \widetilde{J},\widetilde{\theta })$ between two
pseudo-Hermitian manifolds a $(H,\widetilde{H})$-biholomorphism if
$f$ and $f^{-1}$ are $(H,\widetilde{H})$-holomorphic and
$(\widetilde{H},H)$-holomorphic respectively. Furthermore, if $f$
is foliated, then it is called a foliated
$(H,\widetilde{H})$-biholomorphism.
\enddefinition

\example{Example 5.1} Let $(M,H(M),J,\theta )$ be a
pseudo-Hermitian manifold. For any positive function $u$ on $M$,
we set $\widetilde{\theta }=u\theta $. Then
$$
d\widetilde{\theta }=du\wedge \theta +ud\theta.
$$
Write $\widetilde{\xi}=\lambda \xi+\widetilde{T}^\alpha \eta
_\alpha +\widetilde{ T}^{\overline{\alpha }}\eta
_{\overline{\alpha }}$. By requiring that $i_{\widetilde{\xi}}
\widetilde{\theta }=1$ and $i_{\widetilde{\xi}}d\widetilde{\theta
}=0$, one gets
$$
\lambda =u^{-1},\quad \widetilde{T}^\alpha =iu^{-1}\eta
_{\overline{ \alpha }}(\log u),\quad
\widetilde{T}^{\overline{\alpha }}=-iu^{-1}\eta _\alpha (\log u).
$$
Obviously
$$
id_M:(M,H(M),J,\theta )\rightarrow (M,H(M),J,\widetilde{\theta })
$$
is a $(H,\widetilde{H})$-biholomorphism. Note that
$\widetilde{H}(M)=H(M)$ in this example. If $u$ is not constant,
then $\widetilde{\xi} \nparallel \xi$, and thus $id_M$ is not
foliated. This provides us an example of non-foliated
$(H,\widetilde{H})$-biholomorphisms. In this circumstance, we know
from Proposition 5.2 that $id_M:(M,H(M),J,\theta )\rightarrow
(M,H(M),J,\widetilde{\theta })$ is not
$(H,\widetilde{H})$-harmonic. When $u$ is constant, $
id_M:(M,H(M),J,\theta )$ $\rightarrow (M,H(M),J,\widetilde{\theta
})$ is clearly a foliated $(H,\widetilde{H})$-biholomorphism.
\endexample

\example{Example 5.2} Let $S=(\xi ,\theta ,J,g_\theta )$ be a
Sasakian structure on $M$. Then we have its Reeb foliation $F_\xi$
and the associated quotient vector bundle $\nu (F_\xi)$. Let
$\pi_\nu:TM\rightarrow \nu(F_\xi )\text{ }(X\mapsto [X] \text{ for
any } X\in TM)$ be the natural projection. The background
structure $(\xi ,\theta ,J,g_\theta )$ induces a transversal
complex structure $J_\nu$ on $\nu(F_\xi )$ by $J_\nu[X]=[JX]$.
According to [BG], [BGS], we consider the following space of
Sasakian structures with fixed Reeb vector field $\xi$ and fixed
transversal complex structure $J_\nu$:
$$S(\xi ,J_\nu)=\{\text{Sasakian structure } \widetilde{S}=
(\widetilde{\xi },\widetilde{\theta },\widetilde{J},g_{\widetilde{\theta
}})\text{ on } M \text{ }|\text{ } \widetilde{\xi }=\xi,
\widetilde{J}_\nu=J_\nu\}
$$
where $\widetilde{J}_{\nu}$ denotes the transversal complex
structure on $\nu (F_\xi)$ induced by $\widetilde{J}$. For any
$\widetilde{S}\in S(\xi ,J_\nu)$, we assert that
$id_M:(M,S)\rightarrow (M,\widetilde{S})$ is a foliated
$(H,\widetilde{H})$-biholomorphism. To prove this assertion, we
may decompose any $X\in TM$ as $X=a\xi +X_H=b\xi
+X_{\widetilde{H}}$ with $X_H\in H(M)=\ker \theta $ and
$X_{\widetilde{H}}\in \widetilde{H}(M)=\ker \widetilde{\theta }$
for some $a,b\in R$. Then $\pi
_{\widetilde{H}}(X_H)=X_{\widetilde{H}}$, and (1.7) implies that
$JX=JX_H$ and $\widetilde{J}X=\widetilde{J}X_{\widetilde{H}}$. By
definition, $J_\nu[X]=[JX_H]$ and
$\widetilde{J}_\nu[X]=[\widetilde{J}X_{\widetilde{H}}]$. Since
$\pi _\nu :H(M)\rightarrow \nu (F_\xi )$ and $\pi _\nu
:\widetilde{H}(M)\rightarrow \nu (F_\xi )$ are both vector bundle
isomorphisms, we see that $\widetilde{J}_\nu =J_\nu $ if and only
if $\pi _\nu \widetilde{J}X_{\widetilde{H}}=\pi _\nu JX_H$ for any
$X\in TM$. On the other hand, $di_{H,\widetilde{H}}\circ
JX=\widetilde{J}di_{H,\widetilde{H}}X$ if and only if $\pi
_{\widetilde{H}}JX_H=\widetilde{J}\pi _{\widetilde{H}}(X_H)$ that
is, $\pi _{\widetilde{H}}JX_H =\widetilde{J}X_{\widetilde{H}}$.
Taking projection $\pi_\nu$ on both sides of the previous
equality, we get the result.

Recall that a function $u$ on a foliated manifold $M$ is called
basic if it is constant along the leaves. Denote by
$C_B^{\infty}(M)$ the space of smooth basic functions on $M$. In
terms of [BG1,2], we know that the space $S(\xi ,J_\nu )$ is an
affine space modeled on $(C_B^\infty(M) /R)\times(C_B^\infty(M)
/R)\times H^1(M,Z)$. Indeed, if $S=(\xi ,\theta ,J,g_\theta )$ is
a given Sasakian structure in $S(\xi ,J_\nu )$, any other Sasakian
structure $\widetilde{S}=(\widetilde{\xi },\widetilde{\theta
},\widetilde{J},g_{\widetilde{\theta }})$ in it is determined by
real valued basic functions $\varphi $, $\psi \in C_B^\infty (M)$
up to a constant and $\alpha \in H^1(M,Z)$ a harmonic $1$-form
such that
$$
\widetilde{\theta }=\theta +d^c\varphi +\alpha +d\psi \tag{5.11}
$$
where $d^c=\frac{\sqrt{-1}}2(\overline{\partial }-\partial )$.
Thus one may denote the Sasakian structure $\widetilde{S}$ by
$\widetilde{S}_{\varphi ,\psi ,\alpha }$ if $S$ is fixed. In order
to use the notations in \S 2, we write $f=id_M$ and express (5.11)
as
$$
f^{*}\widetilde{\theta }=\theta +\frac{\sqrt{-1}}2\left( \varphi
_{\overline{k}}\theta ^{\overline{k}}-\varphi _k\theta ^k\right)
+\alpha _k\theta ^k+\alpha _{\overline{k}}\theta
^{\overline{k}}+\psi _k\theta ^k+\psi _{\overline{k}}\theta
^{\overline{k}}\tag{5.12}
$$
where $\varphi _k=\eta _k(\varphi )$, $\psi _k=\eta _k(\psi )$ and
$\alpha _k=\alpha (\eta _k)$. Then
$$
f_k^0=-\frac{\sqrt{-1}}2\varphi _k+\alpha _k+\psi _k,\
f_{\overline{k}}^0=\frac{\sqrt{-1}}2\varphi _{\overline{k}}+\alpha
_{\overline{k}}+\psi _{\overline{k}}.
$$
It follows that
$$
\aligned
f_{k\overline{k}}^0+f_{\overline{k}k}^0&=\frac{\sqrt{-1}}2\left(
\varphi _{\overline{k}k}-\varphi _{k\overline{k}}\right) +\alpha
_{k\overline{k}}+\alpha _{\overline{k}k}+\psi
_{k\overline{k}}+\psi _{\overline{k}k}\\
&=-\frac{m\sqrt{-1}}2\xi (\varphi )+\alpha _{k\overline{k}}+\alpha
_{\overline{k}k}+\psi _{k\overline{k}}+\psi _{\overline{k}k}
\endaligned\tag {5.13}
$$
where we use the known formula $\varphi _{k\overline{l}}-\varphi
_{\overline{l}k}=\sqrt{-1}\delta _{kl}\xi (\varphi )$ for a
function in the second equality (see for example [Le], [CDRY]).
Since $\varphi $ is basic and $\alpha $ is harmonic, we find from
(5.13) that
$$
f_{k\overline{k}}^0+f_{\overline{k}k}^0=\bigtriangleup_H(\psi).
$$
It follows that $id_M:(M,S)\rightarrow (M,\widetilde{S}_{\varphi
,0,\alpha })$ is a special foliated
$(H,\widetilde{H})$-biholomorphism. In general, if $\psi$ is not
constant, then $id_M:(M,S)\rightarrow (M,\widetilde{S}_{\varphi
,\psi ,\alpha })$ is not special in the sense of Definition 3.2.
\endexample

\proclaim{Theorem 5.4} Suppose that $M$, $N$ are Sasakian
manifolds and $f:M\rightarrow N$ is a foliated
$(H,\widetilde{H})$-harmonic map. Let $U$ be a nonempty open
subset of $M$. If $f$ is $(H,\widetilde{H})$-holomorphic (resp.
$(H,\widetilde{H})$-antiholomorphic) on $U$, then $f $ is
$(H,\widetilde{H})$-holomorphic (resp.
$(H,\widetilde{H})$-antiholomorphic) on $M$.
\endproclaim

\demo{Proof} Without loss of generality, we assume that $f$ is
$(H,\widetilde{H})$-holomorphic on $U$. Although $f$ satisfies a
PDE system of `subelliptic type', to the author's knowledge, the
unique continuation theorem is still open for such kind of PDE
systems. We will try to show this theorem by using the Aroszajin's
continuation theorem for elliptic PDE systems and the moving frame
method.

Let $\Omega $ be the largest connected open subset of $M$
containing $U$ such that $\overline{\partial }f_{H,\widetilde{H}}$
vanishes identically on $\Omega $. Suppose $\Omega $ has a
boundary point $q$. Let $W$ be a connected open neighborhood of
$q$ in $M$ such that

i) there exists a frame field $\{\xi,\eta _1,...\eta _m,\eta
_{\overline{1} },...,\eta _{\overline{m}}\}$ of $TM^C$ on some
open neighborhood of the closure of $W$ and

ii) there exists a frame field $\{\widetilde{\xi},\widetilde{\eta
}_1,..., \widetilde{\eta }_n,\widetilde{\eta
}_{\overline{1}},...,\widetilde{\eta }_{ \overline{n}}\}$ of
$TN^C$ on some open neighborhood of the closure of $f(W)$ .

The assumption that $f$ is foliated means that
$$
f_0^\alpha =f_0^{\overline{\alpha }}=0.  \tag{5.11}
$$
Since $N$ is Sasakian, the $(H,\widetilde{H})$-harmonic equation
for $f$ becomes
$$
f_{k\overline{k}}^\alpha +f_{\overline{k}k}^\alpha =0.  \tag{5.12}
$$
By definition of covariant derivatives, we have
$$
D\overline{\partial }f_{H,\widetilde{H}}=df_{\overline{j}}^\alpha
\otimes \theta ^{\overline{j}}\otimes \widetilde{\eta }_\alpha
+f_{\overline{ j}}^\alpha \nabla \theta ^{\overline{j}}\otimes
\widetilde{\eta }_\alpha +f_{\overline{j}}^\alpha \theta
^{\overline{j}}\otimes \widetilde{\nabla } \widetilde{\eta
}_\alpha\tag{5.13}
$$
and
$$
\aligned
D^2\overline{\partial }f_{H,\widetilde{H}}=&\nabla
df_{\overline{j}}^\alpha \otimes \theta ^{\overline{j}}\otimes
\widetilde{\eta }_\alpha +df_{ \overline{j}}^\alpha \otimes \nabla
\theta ^{\overline{j}}\otimes \widetilde{ \eta }_\alpha
+df_{\overline{j}}^\alpha \otimes \theta ^{\overline{j}
}\otimes \widetilde{\nabla }\widetilde{\eta }_\alpha \\
&+df_{\overline{j}}^\alpha \otimes \nabla \theta
^{\overline{j}}\otimes \widetilde{\eta }_\alpha
+f_{\overline{j}}^\alpha \nabla ^2\theta ^{\overline{j}}\otimes
\widetilde{\eta }_\alpha +f_{\overline{j}}^\alpha \nabla \theta
^{\overline{j}}\otimes \widetilde{\nabla }\widetilde{\eta }
_\alpha \\
&+df_{\overline{j}}^\alpha \otimes \theta ^{\overline{j}}\otimes
\widetilde{ \nabla }\widetilde{\eta }_\alpha
+f_{\overline{j}}^\alpha \nabla \theta ^{ \overline{j}}\otimes
\widetilde{\nabla }\widetilde{\eta }_\alpha +f_{
\overline{j}}^\alpha \theta ^{\overline{j}}\otimes
\widetilde{\nabla }^2 \widetilde{\eta }_\alpha.
\endaligned\tag{5.14}
$$
Now we compute the trace Laplacian of the section
$\overline{\partial }f_{H,\widetilde{H}}$ as follows:
$$
\aligned \bigtriangleup \overline{\partial
}f_{H,\widetilde{H}}=&tr_{g_\theta }D^2 \overline{\partial
}f_{H,\widetilde{H}}\\=&(f_{\overline{j}k\overline{k} }^\alpha
+f_{\overline{j}\overline{k}k}^\alpha +f_{\overline{j}00}^\alpha
)\theta ^{\overline{j}}\otimes \widetilde{\eta }_\alpha \\
=&(\bigtriangleup _Mf_{\overline{j}}^\alpha )\theta
^{\overline{j}}\otimes \widetilde{\eta }_\alpha +tr_{g_\theta
}\{df_{\overline{j}}^\alpha \otimes \nabla \theta
^{\overline{j}}\otimes \widetilde{\eta }_\alpha +df_{\overline{
j}}^\alpha \otimes \theta ^{\overline{j}}\otimes \widetilde{\nabla
}\widetilde{\eta }_\alpha \\
&+df_{\overline{j}}^\alpha \otimes \nabla \theta
^{\overline{j}}\otimes \widetilde{\eta }_\alpha
+f_{\overline{j}}^\alpha \nabla ^2\theta ^{ \overline{j}}\otimes
\widetilde{\eta }_\alpha +f_{\overline{j}}^\alpha \nabla \theta
^{\overline{j}}\otimes \widetilde{\nabla }\widetilde{\eta }
_\alpha \\
&+df_{\overline{j}}^\alpha \otimes \theta ^{\overline{j}}\otimes
\widetilde{ \nabla }\widetilde{\eta }_\alpha
+f_{\overline{j}}^\alpha \nabla \theta ^{ \overline{j}}\otimes
\widetilde{\nabla }\widetilde{\eta }_\alpha +f_{
\overline{j}}^\alpha \theta ^{\overline{j}}\otimes
\widetilde{\nabla }^2 \widetilde{\eta }_\alpha \}
\endaligned \tag{5.15}
$$
where $\bigtriangleup _M$ denotes the Laplace-Beltrami operator
acting on functions. Since $M$ and $N$ are Sasakian, we derive
from the second equation of (2.17) and (5.11) that
$$
f_{\overline{j}0}^\alpha =0 \tag{5.16}
$$
which yields
$$
f_{\overline{j}00}^\alpha =0.  \tag{5.17}
$$
Using (2.17), (2.38), (2.44) and (5.12), we discover
$$
\aligned f_{\overline{j}k\overline{k}}^\alpha
+f_{\overline{j}\overline{k}k}^\alpha
=&f_{k\overline{j}\overline{k}}^\alpha
+f_{\overline{k}\overline{j}k}^\alpha
\\
=&f_k^\beta \widehat{R}_{\beta \gamma \overline{\delta }}^\alpha
(f_{ \overline{k}}^\gamma f_{\overline{j}}^{\overline{\delta
}}-f_{\overline{j} }^\gamma f_{\overline{k}}^{\overline{\delta
}})-f_{\overline{t}}^\alpha R_{
\overline{k}k\overline{j}}^{\overline{t}}+f_{\overline{k}}^\beta
\widehat{R} _{\beta \gamma \overline{\delta }}^\alpha (f_k^\gamma
f_{\overline{j}}^{ \overline{\delta }}-f_{\overline{j}}^\gamma
f_k^{\overline{\delta }}).
\endaligned\tag{5.18}
$$
Consequently
$$
|f_{\overline{j}k\overline{k}}^\alpha
+f_{\overline{j}\overline{k}k}^\alpha +f_{\overline{j}00}^\alpha
|\leq C\sum_{l,\beta }|f_{\overline{l}}^\beta | \tag{5.19}
$$
on $W$ for some positive number $C$. From (5.15) and (5.19), we
find that there is a positive number $C^{\prime }$ such that
$$
|\bigtriangleup _M(f_{\overline{j}}^\alpha )|\leq C^{\prime
}(\sum_{l,\beta }|\nabla f_{\overline{l}}^\beta |+\sum_{l,\alpha
}|f_{\overline{l}}^\alpha |) \tag{5.20}
$$
where $\nabla $ denotes the gradient of the functions
$\{f_{\overline{j} }^\alpha \}$. By applying the Aronszajn's
unique continuation theorem (cf. [Ar], [PRS]) to the system of
functions $Re\{f_{\overline{j}}^\alpha \}$,
$Im\{f_{\overline{j}}^\alpha \}$ ($1\leq j\leq m$, $1\leq \alpha
\leq n$) and to the elliptic operator $\bigtriangleup _M$, we
conclude from the identical vanishing of
$Re\{f_{\overline{j}}^\alpha \}$, $Im\{f_{\overline{j}}^\alpha \}$
on $W\cap \Omega $ that $Re\{f_{\overline{j}}^\alpha \}$,
$Im\{f_{\overline{j}}^\alpha \}$ vanish identically on $W$. This
contradicts the fact that $q$ is a boundary point of $\Omega $.
Hence $\Omega =M$, which implies that $\overline{\partial
}f_{H,\widetilde{H}}=0$ on the whole domain manifold.\qed

\enddemo

\remark{Remark 5.2} Note that we verify the structural assumptions
of Aronszajn-Cordes in the proof of Theorem 5.4 by adopting the
moving frame method, whose advantage is its operability. This
method will be used again in the appendix.
\endremark

Note that $f$ is both $(H,\widetilde{H})$-holomorphic and
$(H,\widetilde{H})$-antiholomorphic if and only if
$df_{H,\widetilde{H}}=0$, that is, $f$ is horizontally constant.

\proclaim{Corollary 5.5} Suppose that $M$, $N$ are Sasakian
manifolds and $f:M\rightarrow N$ is a foliated
$(H,\widetilde{H})$-harmonic map. Let $U$ be a nonempty open
subset of $M$. If $ f $ is horizontally constant on $U$, then $f$
is horizontally constant on $M$ .
\endproclaim
\demo{Proof} Since $f$ is horizontally constant on $U$, it is both
$(H,\widetilde{H})$--holomorphic and
$(H,\widetilde{H})$-antiholomorphic on $U$. It follows from
Theorem 5.4 that $f$ is both $(H,\widetilde{H})$-holomorphic and
$(H,\widetilde{H})$-anitholomorphic on the whole $M$. Hence we
conclude that $f$ is horizontally constant on $M$. \qed
\enddemo

It would be interesting to know whether the Sasakian and foliated
conditions in Theorem 5.4 or Corollary 5.5 can be removed or not.
In terms of Proposition 5.2, one way to give an answer is to
establish a unique continuation result for the foliated property.
We would like to propose the following question:

\proclaim{Question} Suppose $f:M\rightarrow N$ is a
$(H,\widetilde{H})$-harmonic map between two pseudo-Hermitian
manifolds or even Sasakian manifolds . If $f$ is foliated on a
nonempty open subset $U$ of $M$, can we deduce that $f$ is
foliated on the whole $M$?
\endproclaim

Though a general Sasakian manifold is not a global Riemannian
submersion over a K\"ahler manifold, the following result will
help us to understand the general local picture and properties
about $(H,\widetilde{H})$-holomorphic maps between Sasakian
manifolds.

\proclaim{Proposition 5.6} Suppose $(M^{2m+1},H(M),J,\theta )$ and
$(N^{2n+1},\widetilde{H}(N),\widetilde{J},\widetilde{\theta })$
are compact Sasakian manifolds which are the total spaces of
Riemannian submersions $\pi :M\rightarrow B$ and $\widetilde{\pi
}:N\rightarrow \widetilde{B}$ over compact K\"ahler manifolds $B$
and $\widetilde{B}$ respectively. Suppose $f:M\rightarrow N$ is a
foliated map which induces a map $h:B\rightarrow \widetilde{B}$
between the base manifolds. Then $f$ is a foliated
$(H,\widetilde{H})$-holomorphic (resp.
$(H,\widetilde{H})$-antiholomorphic) map if and only if $h$ is a
holomorphic (resp. anti-holomorphic) map.
\endproclaim
\demo{Proof} Since $f$ is foliated, we have
$h\circ \pi =\widetilde{\pi }\circ f$, and thus $dh\circ d\pi
=d\widetilde{\pi }\circ df$. Denote by $J_1$ and $J_2$ the complex
structures of $B$ and $\widetilde{B}$ respectively. Since $d\pi
\circ J=J_1\circ d\pi $ and $d\widetilde{\pi }\circ
\widetilde{J}=J_2\circ d\widetilde{\pi }$, we immediately have
that $f$ is $(H,\widetilde{H})$-holomorphic (resp.
$(H,\widetilde{H})$-antiholomorphic) if and only if $h$ is
holomorphic (resp. antiholomorphic).\qed\enddemo

\remark{Remark 5.2} Suppose now that $M$ and $N$ are two general
Sasakian manifolds, which are not necessarily total spaces of
Riemannian submersions, and $f:M\rightarrow N$ is a foliated map.
Let $p\ $be any point in $M$ and $q=f(p)$. We have foliated
neighborhoods $U_1$ and $\widetilde{U}_2$ of $p$ and $q$
respectively, together with Riemannian submersions $\pi
_1:U_1\rightarrow W_1$ and $\widetilde{\pi
}_2:\widetilde{U}_2\rightarrow \widetilde{W}_2$ over two K\"ahler
manifolds $W_1$ and $\widetilde{W}_2$. Assuming that
$f(U_1)\subset \widetilde{U}_2$, then $f$ induces locally a map
$h_p:W_1\rightarrow \widetilde{W}_2$. According to Proposition
5.6, we find that $f:M\rightarrow N$ is
$(H,\widetilde{H})$-holomorphic (resp.
$(H,\widetilde{H})$-antiholomorphic) if and only if the locally
induced map $h_p$ for each $p\in M$ is holomorphic (resp.
anti-holomorphic).
\endremark

\definition{Definition 5.3}A foliated map $f:(M,H(M),\theta
,J)\rightarrow (N,\widetilde{H}(N),\widetilde{J},\widetilde{\theta
})$ between two pseudo-Hermtian manifolds is called a horizontally
one-to-one map if it induces a one-to-one map $h:M/F_\xi
\rightarrow N/\widetilde{F}_{\widetilde{\xi }}$ between the spaces
of leaves.
\enddefinition

\proclaim{Proposition 5.7}Let $M$ and $N$ be Sasakian manifolds
and let $f:M\rightarrow N$ be a foliated
$(H,\widetilde{H})$-holomorphic map. If $f$ is both one-to-one and
horizontally one-to-one, then $f^{-1}$ is a
$(\widetilde{H},H)$-holomorphic map, that is, $f$ is a foliated
$(H,\widetilde{H})$-biholomorphism.
\endproclaim
\demo{Proof} For any $p\in M$, we let $\pi _1:U_1\rightarrow W_1$,
$\pi _2:\widetilde{U}_2\rightarrow \widetilde{W}_2$ and
$h_p:W_1\rightarrow \widetilde{W}_2$ be as in Remark 5.2. Then
$h_p:W_1\rightarrow \widetilde{W}_2$ is holomorphic. Since $f$ is
horizontally one-to-one, we see that $h_p:W_1\rightarrow h_p(W_1)$
is one-to-one. By a known result for holomorphic maps (cf.
Proposition 1.1.13 in [Hu]), we know that $h_p$ is a local
biholomorphic, that is, $h_p(W_1)$ is an open set of
$\widetilde{W}_2$ and $h_p:W_1\rightarrow h_p(W_1)$ is
biholomorphic.

Set $\widetilde{V_2}=f(U_1)\subset \widetilde{U}_2$. Clearly
$f^{-1}:N\rightarrow M$ is also a foliated map with
$f^{-1}(\widetilde{V_2})=U_1$, and  $f^{-1}$ induces the
holomorphic map $h_p^{-1}:\widetilde{h}_p(W_1)\rightarrow W_1$.
Using Proposition 5.6, we find that  $f^{-1}\ $is
$(\widetilde{H},H)$-holomorphic at $q=f(p)$. Since $p$ is
arbitrary, we conclude that $f$ is a foliated
$(H,\widetilde{H})$-biholomorphism.\qed
\enddemo

\heading{\bf 6. Lichnerowicz type results}
\endheading
\vskip 0.3 true cm

By definition, one gets from (5.6) that (cf. [Do])
$$
\aligned |\partial f_{H,\widetilde{H}}|^2=&\langle\partial
f_{H,\widetilde{H}}(\eta _j),\overline{\partial f_{H,\widetilde{H}}(\eta _j)}\rangle\\
=&\frac14\sum_{j=1}^m\{\langle
df_{H,\widetilde{H}}(e_j),df_{H,\widetilde{H}
}(e_j)\rangle+\langle df_{H,\widetilde{H}}(Je_j),df_{H,\widetilde{H}}(Je_j)\rangle\\
&+2\langle
df_{H,\widetilde{H}}(Je_j),\widetilde{J}df_{H,\widetilde{H}}(e_j)\rangle\}
\endaligned\tag{6.1}
$$
and
$$
\aligned |\overline{\partial
}f_{H,\widetilde{H}}|^2=&\langle\overline{\partial }f_{H,
\widetilde{H}}(\eta _{\overline{j}}),\overline{\overline{\partial
}f_{H,
\widetilde{H}}(\eta _{\overline{j}})}\rangle \\
=&\frac 14\sum_{j=1}^m\{\langle
df_{H,\widetilde{H}}(e_j),df_{H,\widetilde{H}
}(e_j)\rangle+\langle df_{H,\widetilde{H}}(Je_j),df_{H,\widetilde{H}}(Je_j)\rangle \\
&-2\langle
df_{H,\widetilde{H}}(Je_j),\widetilde{J}df_{H,\widetilde{H}}(e_j)\rangle\}.
\endaligned\tag{6.2}
$$
Thus
$$
\frac 12|df_{H,\widetilde{H}}|^2=|\partial
f_{H,\widetilde{H}}|^2+| \overline{\partial
}f_{H,\widetilde{H}}|^2.\tag{6.3}
$$
Define
$$
E_{H,\widetilde{H}}^{\prime }(f)=\int_M|\partial
f_{H,\widetilde{H} }|^2dv_{\theta } ,\quad
E_{H,\widetilde{H}}^{\prime \prime }(f)=\int_M|\overline{\partial
}f_{H,\widetilde{H}}|^2 dv_{\theta } .  \tag{6.4}
$$
Then $E_{H,\widetilde{H}}(f)=E_{H,\widetilde{H}}^{\prime
}(f)+E_{H, \widetilde{H}}^{\prime \prime }(f)$. Set
$$
k_{H,\widetilde{H}}(f)=|\partial
f_{H,\widetilde{H}}|^2-|\overline{
\partial }f_{H,\widetilde{H}}|^2,\quad K_{H,\widetilde{H}}(f)=E_{H,
\widetilde{H}}^{\prime }(f)-E_{H,\widetilde{H}}^{\prime \prime
}(f).\tag{6.5}
$$
\proclaim{Lemma 6.2} Set $\omega ^M=d\theta $ and $\omega
^N=d\widetilde{\theta }$. Then
$$
k_{H,\widetilde{H}}(f)=\langle\omega ^M,f^{*}\omega ^N\rangle.
$$
\endproclaim

\demo{Proof} Choose an orthonormal frame
$\{e_1,...,e_m,Je_1,...,Je_m\}$ of $H(M)$. Using (1.9), we deduce
that
$$
\aligned \langle\omega ^M,f^{*}\omega
^N\rangle=&\sum_{i<j}\{(f^{*}\omega ^N)(e_i,e_j)\omega
^M(e_i,e_j)+(f^{*}\omega ^N)(Je_i,Je_j)\omega ^M(Je_i,Je_j)\} \\
&+\sum_{i,j}(f^{*}\omega ^N)(e_i,Je_j)\omega ^M(e_i,Je_j) \\
=&\sum_i\langle\widetilde{J}df(e_i),df(Je_i)\rangle \\
=&\sum_i\langle\widetilde{J}df_{H,\widetilde{H}}(e_i),df_{H,\widetilde{H}}(Je_i)\rangle.
\endaligned\tag{6.6}
$$
Consequently (6.1), (6.2) and (6.6) imply the lemma.\qed
\enddemo

We need the following lemma:

\proclaim{Lemma 6.3 (\text{Homotopy Lemma, }cf. [Lic], [EL])} Let
$f_t:M\rightarrow N$ be a smooth family of maps between the smooth
manifolds $M$ and $N$, parameterized by the real number $t$, and
let $\omega $ be a closed two-form on $N$. Then
$$
\frac \partial {\partial t}(f_t^{*}\omega )=d\big(
f_t^{*}i(\frac{\partial f_t}{\partial t})\omega \big)
$$
where $i(X)\omega $ denotes the interior product of the vector $X$
with the two-form $\omega $.
\endproclaim

\proclaim{Lemma 6.4} Let $f_t:(M^{2m+1},H(M),J,\theta )\rightarrow
(N^{2n+1},\widetilde{H}(N), \widetilde{J},\widetilde{\theta })$ be
a family of smooth maps between two pseudo-Hermitian manifolds.
Then
$$
\frac d{dt}K_{H,\widetilde{H}}(f_t)=m\int_Md\widetilde{\theta }
(v_t,df_t(\xi))dv_\theta
$$
where $v_t=\partial f_t/\partial t$.
\endproclaim

\demo{Proof} Clearly Lemmas 6.2, 6.3 imply that
$$
\aligned \frac d{dt}K_{H,\widetilde{H}}(f_t)&=\int_M\langle\frac
\partial {\partial
t}f_t^{*}\omega ^N,\omega ^M\rangle dv_\theta \\
&=\int_M\langle d\sigma _t,\omega ^M\rangle dv_\theta \\
&=\int_M\langle \sigma _t,\delta \omega ^M\rangle dv_\theta
\endaligned \tag{6.7}
$$
where $\sigma _t=f_t^{*}i(\frac{\partial f_t}{\partial t})\omega
^N$. Choose an orthonormal frame field $\{e_A\}_{A=0}^{2m}=$
$\{\xi,e_1,...,e_{2m}\}$ on $ M $. From Lemma 1.2, we get
$$
\nabla _{e_A}^\theta X=\nabla _{e_A}X-(\frac 12d\theta
(e_A,X)+A(e_A,X))\xi \tag{6.8}
$$
and
$$
\nabla _{e_A}^\theta \xi=\tau (e_A)+\frac 12Je_A  \tag{6.9}
$$
for any $X\in \Gamma (H(M))$ and $1\leq A\leq 2m$. Using (1.4),
(1.13) and (6.8), we compute the codifferential $(\delta \omega
^M)(X)$ for $X\in \Gamma (H(M))$ as follows:
$$
\aligned (\delta \omega ^M)(X)&=-\sum_{A=0}^{2m}(\nabla
_{e_A}^\theta d\theta )(e_A,X)
\\
&=-\sum_{A=0}^{2m}\{e_A[d\theta (e_A,X)]-d\theta (\nabla
_{e_A}^\theta
e_A,X)-d\theta (e_A,\nabla _{e_A}^\theta X)\} \\
&=-\sum_{A=1}^{2m}\{e_A[d\theta (e_A,X)]-d\theta (\nabla
_{e_A}e_A,X)\}+\sum_{A=1}^{2m}d\theta (e_A,\nabla _{e_A}X) \\
&=-\sum_{A=1}^{2m}(\nabla _{e_A}d\theta )(e_A,X) \\
&=0,
\endaligned\tag{6.10}
$$
due to the fact that $\nabla d\theta=0$. Next
$$
\aligned
(\delta \omega ^M)(\xi)&=\sum_{A=0}^{2m}d\theta (e_A,\nabla _{e_A}^\theta \xi) \\
&=\sum_{A=1}^{2m}d\theta (e_A,\tau (e_A)+\frac 12Je_A) \\
&=\sum_{A=1}^{2m}[g_\theta (Je_A,\tau (e_A))+\frac 12g_\theta (Je_A,Je_A)] \\
&=m.
\endaligned \tag{6.11}
$$
It follows from (6.7), (6.10) and (6.11) that
$$
\aligned
\frac d{dt}K_{H,\widetilde{H}}(f_t)&=m\int_M\sigma _t(\xi)dv_\theta \\
&=m\int_Md\widetilde{\theta }(v_t,df_t(\xi))dv_\theta .
\endaligned
$$
\qed
\enddemo

\definition{Definition 6.2} Let $f_0$ and $f_1$ be two maps between two
pseudo-Hermitian manifolds $M$ and $N$. We say that $f_0$ and
$f_1$ are vertically homotopic if there exists a map $F:M\times
[0,1]\rightarrow N$ such that $F(\cdot ,0)=f_0,\ F(\cdot ,1)=f_1$
and for each point $x\in M$, the tangent vector at each point
along the curve $F(x,\cdot )$ is vertical.
\enddefinition
\proclaim{Theorem 6.5} Let $(M,H(M),J,\theta)$ and
$(N,\widetilde{H}(N),\widetilde{J},\widetilde{\theta })$ be two
pseudoHermitian manifolds. Suppose that $M$ is compact. Then
$K_{H,\widetilde{H}}(f)$ is a smooth vertical homotopy invariant,
that is, if $f_t$ is any smooth $1$-parameter vertical variation
$f_t$ of $f$, then $t\mapsto K_{H,\widetilde{H}}(f_t)$ is a
constant map.
\endproclaim
\demo{Proof} Let $f_0$ and $f_1$ be two maps from $M$ to $N$
through a family of maps $f_t:M\rightarrow N$, $t\in [0,1]$ with
the property that $\partial f_t/\partial t$ is vertical. By Lemma
6.4 and (1.4), we get
$$
\frac d{dt}K_{H,\widetilde{H}}(f_t)=0.
$$
Consequently $t\mapsto K_{H,\widetilde{H}}(f_t)$ is a constant
map. \qed \enddemo

\remark{Remark 6.2} When $N$ is Sasakian, we even have stronger
results. Let $f$ be a map from a pseudo-Hermitian manifold to a
Sasakian manifold and $\{f_t\}_{|t|< \varepsilon }$ a vertical
variation of $f$. Set $v_t=\frac{\partial f_t}{
\partial t}$ and $\Phi (\cdot ,t)=f_t(\cdot )$. Using (1.10), Lemma 2.1 and
a direct computation, we may derive from (6.2) that
$$\aligned
\frac \partial {\partial t}|\overline{\partial
}f_{tH,\widetilde{H} }|^2 =&\frac
12\sum_{j=1}^m\langle\widetilde{\nabla }_{e_j}d\Phi
_{H,\widetilde{H} }(\frac \partial {\partial t}),d\Phi
_{H,\widetilde{H}}(e_j)\rangle+\langle\widetilde{\nabla
}_{Je_j}d\Phi _{H,\widetilde{H}}(\frac \partial {\partial
t}),d\Phi
_{H,\widetilde{H}}(Je_j)\rangle \\
&-\langle\widetilde{\nabla }_{Je_j}d\Phi _{H,\widetilde{H}}(\frac
\partial
{\partial t}),\widetilde{J}d\Phi
_{H,\widetilde{H}}(e_j)\rangle-\langle\Phi _{H,
\widetilde{H}}(Je_j),\widetilde{J}\widetilde{\nabla }_{e_j}d\Phi
_{H,\widetilde{H}}(\frac \partial {\partial t})\rangle] \\
=&0.
\endaligned
$$
Similarly we have $\frac \partial {\partial t}|\partial
f_{tH,\widetilde{H} }|^2=0$. This shows that the horizontal
partial energy densities are preserved under the vertical
deformation. Consequently $E_{H,\widetilde{H} }(f_t)$,
$E_{H,\widetilde{H}}^{\prime }(f_t)$ and $E_{H,\widetilde{H}
}^{\prime \prime }(f_t)$ are invariant under the vertical
variation of $f$.
\endremark
\proclaim{Theorem 6.6} Let $(M,H(M),J,\theta )$ and
$(N,\widetilde{H}(N),\widetilde{J},\widetilde{ \theta })$ be two
pseudoHermitian manifolds and let $f:M\rightarrow N$ be a foliated
map. Suppose that $M$ is compact. Then $K_{H,\widetilde{H}}(f)$ is
a smooth foliated homotopy invariant, that is, $t\mapsto
K_{H,\widetilde{H} }(f_t)$ is a constant map for any smooth
$1$-parameter of foliated maps $f_t$ with $f_0=f$.
\endproclaim
\demo{Proof} Suppose $f_t$ is a smooth $1$-parameter of foliated
maps with $f_0=f$ . Since $f_t$ is foliated, $df_t(T)$ is
vertical. Hence Lemma 6.4 yields
$$
\frac d{dt}K_{H,\widetilde{H}}(f_t)=0
$$
that is, $t\mapsto K_{H,\widetilde{H}}(f_t)$ is constant. \qed
\enddemo

\remark{Remark 6.3} Although pseudo-Hermitian foliations are not
K\"ahler foliations in general, we would mention that the authors
in [BD] proved a similar result for foliated maps between K\"ahler
foliations.
\endremark

\proclaim{Theorem 6.7} Suppose $f:(M,H(M),J,\theta )\rightarrow
(N,\widetilde{H}(N),\widetilde{J}, \widetilde{\theta })$ is either
a $(H,\widetilde{H})$-holomorphic map or a
$(H,\widetilde{H})$-antiholomorphic map. Then

(i) $f$ is an absolute minimum of $E_{H,\widetilde{H}}$ among its
vertical homotopy class;

(ii) If $f$ is foliated, then it is also an absolute minimum of
$E_{H, \widetilde{H}}$ among its foliated homotopy class.
\endproclaim

\demo{Proof} Without loss of generality, we assume that $f$ is a
$(H,\widetilde{H})$-holomorphic map. Let $\widetilde{f}$ be any
smooth map in the vertical homotopy class of $f$. By Theorem 6.5,
we have
$$
E_{H,\widetilde{H}}^{\prime }(\widetilde{f})-E_{H,\widetilde{H}
}^{\prime \prime }(\widetilde{f})=E_{H,\widetilde{H}}^{\prime
}(f)-E_{H, \widetilde{H}}^{\prime \prime
}(f)=E_{H,\widetilde{H}}^{\prime }(f)\tag{6.12}
$$
Then $E_{H,\widetilde{H}}^{\prime }(f)\leq
E_{H,\widetilde{H}}^{\prime }( \widetilde{f})$ and thus
$E_{H,\widetilde{H}}(f)\leq E_{H,\widetilde{H}}( \widetilde{f})$.
This proves that $f$ is an absolute minimum of $E_{H,
\widetilde{H}}$ among its vertical homotopy class. Similarly one
may prove that $f$ is an absolute minimum of $E_{H,\widetilde{H}}$
among its foliated homotopy class, provided that $f$ is foliated.
\qed
\enddemo

\remark{Remark 6.4} If $f:M\rightarrow N$ is a foliated map, then
any vertical variation $f_t$ of $f$ is clearly a foliated
variation.
\endremark

\proclaim{Corollary 6.8} Let $f:(M^{2m+1},H(M),J,\theta
)\rightarrow (N,\widetilde{H}(N),\widetilde{J} ,\widetilde{\theta
})$ be either a foliated $(H,\widetilde{H})$-holomorphic map or a
foliated $(H,\widetilde{H})$-antiholomorphic map between two
pseudo-Hermitian manifolds. Then

(i) $f$ is a pseudo-harmonic map in the sense of [Pe2], that is,
$f$ is a critical point of $E_{H,\widetilde{H}}(f_t)$ for any
variation $\{f_t\}$ with $(\partial f_t/\partial t)|_{t=0}\in
\Gamma (f^{-1}TN)$;

(ii) If $f_t$ is a foliated variation of $f$, then
$\frac{d^2}{dt^2}E_{H, \widetilde{H}}(f_t)|_{t=0}\geq 0.$
\endproclaim
\demo{Proof} Without loss of generality, we assume that $f$ is a
foliated $(H,\widetilde{H})$-holomorphic map. From Proposition
5.2, one knows that $f$ is pseudo-harmonic, that is, $f$ is a
critical point of $E_{H,\widetilde{H}}(f_t)$ for any horizontal
variation $\{f_t\}$. From Theorem 6.7, it follows that $f$ is also
a critical point of $E_{H,\widetilde{H}}(f_t)$ for any vertical
variation $\{f_t\}$. Hence $f$ is a critical point of $E_{H,
\widetilde{H}}(f_t)$ for any variation $\{f_t\}$. This proves (i).
It is clear that (ii) follows directly from Theorem 6.7. \qed
\enddemo

\heading{\bf 7. Existence of $(H,\widetilde{H})$-harmonic maps
under $\widetilde{K}^H\leq 0$}
\endheading
\vskip 0.3 true cm

We will introduce a subelliptic heat flow for maps between
pseudo-Hermitian manifolds in order to find special
$(H,\widetilde{H})$-harmonic maps between these manifolds. We
always assume that both $M$ and $N$ are compact, and $N$ is
Sasakian in this section.

For a map $f:(M,H(M),J,\theta )\rightarrow
(N,\widetilde{H}(N),\widetilde{J},\widetilde{\theta })$ between
two pseudo-Hermitian manifolds, besides the horizontal
differential $df_{H,\widetilde{H}}:H(M)\rightarrow
\widetilde{H}(N)$, we have the following partial differentials
$df_{L,\widetilde{H}}:L\rightarrow \widetilde{H}(N)$,
$df_{L,\widetilde{L}}:L\rightarrow \widetilde{L}$ and
$df_{H,\widetilde{L}}:H(M)\rightarrow \widetilde{L}$ defined
respectively by:
$$
\aligned
df_{L,\widetilde{H}}=&\pi _{\widetilde{H}}\circ df\circ i_L,\\
df_{L,\widetilde{L}}=&\pi _{\widetilde{L}}\circ df\circ i_L,\\
df_{H,\widetilde{L}}=&\pi_{\widetilde{L}}\circ df\circ i_H,
\endaligned
$$
where $i_L:L\rightarrow TM$, $i_H:H(M)\rightarrow TM$ are the
inclusion morphisms, and $\pi _{\widetilde{H}}:TN\rightarrow
\widetilde{H}(N)$, $\pi _{\widetilde{L}}:TN\rightarrow
\widetilde{L}$ are the natural projection morphisms. The
corresponding energy densities are given respectively by
$$
\aligned e_{L,\widetilde{H}}=&\frac
12|df_{L,\widetilde{H}}|^2=f_0^\alpha f_0^{\overline{\alpha}},\\
e_{L,\widetilde{L}}=&\frac 12|df_{L,\widetilde{L}}|^2=\frac
12(f_0^0)^2,\\
e_{H,\widetilde{L}}=&\frac
12|df_{H,\widetilde{L}}|^2=f_j^0f_{\overline{j}}^0.\endaligned\tag{7.1}
$$
Let us introduce the following partial second fundamental forms:
$$
\aligned &\beta _H=\beta \big(i_H(\cdot),i_H(\cdot)\big),\quad
\beta_{H,\widetilde{H}}=\pi_{\widetilde{H}}(\beta_H),\quad\beta_{H,\widetilde{L}}=\pi_{\widetilde{L}}(\beta _H),\\
&\beta _{L\times H}=\beta \big(i_L(\cdot ),i_H(\cdot )\big),\quad \beta_{H\times L}=\beta \big(i_H(\cdot ),i_L(\cdot)\big),\\
&\beta_{L\times H,\widetilde{H}}=\pi
_{\widetilde{H}}(\beta_{L\times H}),\quad \beta _{H\times
L,\widetilde{H}}=\pi
_{\widetilde{H}}(\beta _{H\times L}),\\
&\beta _{L\times H,\widetilde{L}}=\pi _{\widetilde{L}}(\beta
_{L\times H}),\quad \beta _{H\times L,\widetilde{L}}=\pi
_{\widetilde{L}}(\beta _{H\times L}).\endaligned\tag{7.2}
$$
and set
$$
\tau_H=tr_{G_\theta }\beta_H,\quad \tau
_{H,\widetilde{L}}=tr_{G_\theta }\beta_{H,\widetilde{L}}.\tag{7.3}
$$
Recalling that $\tau _{H,\widetilde{H}}(f)=tr_{G_\theta }\beta
_{H,\widetilde{H}}$ (see the notations in \S 3), we have the
following decomposition
$$
\tau _H(f)=\tau _{H,\widetilde{H}}(f)+\tau
_{H,\widetilde{L}}(f).\tag{7.4}
$$
Hence $\tau _H(f)=0$ if and only if $\tau
_{H,\widetilde{H}}(f)=\tau _{H,\widetilde{L}}(f)=0$, that is, $f$
is a special $(H,\widetilde{H})$-harmonic map.

Now we consider the following evolution problem on $M\times
[0,T)$:
$$
\cases \frac{\partial f_t}{\partial t}&=\tau_H(f_t)\\
f|_{t=0}&=h \endcases\tag{7.5}
$$
where $h:M\rightarrow N$ is a smooth map. Since the horizontal
part of $\tau _H(f)$ is the gradient of the functional
$E_{H,\widetilde{H}}$, the flow (7.5) has a partial variational
structure. In the appendix, we show that, in terms of the Nash
embedding of $N$ into some Euclidean space $R^K$, the PDE system
in (7.5) can be equivalently expressed as the following type of
subelliptic parabolic system (see Theorem B4):
$$
(\bigtriangleup _{H}-\frac{\partial }{\partial
t})u^{a}=P_{bc}^{a}(u)\langle \nabla _{H}u^{b},\nabla
_{H}u^{c}\rangle ,\quad 1\leq a,b,c\leq K \tag{7.6}
$$
for a map $u:M\times \lbrack 0,T)\rightarrow R^{K}$, where $
P_{bc}^{a}:B(N)\rightarrow R$ are functions on a tubular
neighborhood of $N\subset R^{K}$.

We shall apply the regularity theory in [RS] to investigate
solutions of (7.6). Let $U$ be a relatively compact open subset of
$M$, on which there is an orthonormal frame field
$\{e_{A}\}_{A=1,...,2m}$ for $H(M)$. Set $X_{0}=\frac{\partial }{
\partial t}$, $X_{A}=e_{A}$ ($A=1,...,2m$). Clearly $
\{X_{0},X_{1},...,X_{2m}\}$ together with their commutators of $2$
order span the tangent space of $U\times (0,\infty )$ at any
point. In terms of (1.14), we know from [H\"{o}], [FS] that the
operator $\bigtriangleup _{H}-\frac{\partial }{\partial t}$ is
hypoelliptic on $M\times (0,\infty )$.

Let us recall briefly the function spaces adapted with the
differential operator $\bigtriangleup _{H}-\frac{\partial
}{\partial t}$. Let $U_{T}=U\times (0,T)$ for some $T>0$. For a
monomial $X_{A_{1}}\cdots X_{A_{l}}$ with $0\leq A_{s}\leq 2m$,
$s=1,..,l$, its weight is defined as an integer $r_{1}+2r_{2}$,
where $r_{1}$ is the number of $X_{j}$'s that enter with $j$
between $1$ and $2m$, and $r_{2}$ is the number of $X_{0}$'s. We
also write $w(A_{1},...,A_{l})=r_{1}+2r_{2}$. For any integer
$k\geq 0$ and any $p$, $1<p<\infty $, we define $S_{k}^{p}(U_{T})$
to consist of all $u\in L^{p}(U_{T})$ such that
$(X_{i_{1}}X_{i_{2}}\cdots X_{i_{l}})u\in L^{p}(U_{T})$ for all
monomials of weight $\leq k$. For the norm, we take
$$
\Vert u \Vert_{S_{k}^{p}(U_{T})}=\sum_{w(A_{1},...,A_{l})\leq
k}\Vert X_{A_{1}}\cdots X_{A_{l}}f \Vert_{L^{p}(U_{T})},
$$
that is, the sum is taken over all ordered monomials
$X_{A_{1}}\cdots X_{A_{l}}$ of weight $\leq k$. Using a $C^{\infty
}$ partition of unity subordinate to a finite open cover
$\{U_{j}\}$ of $M$, one may define the space $S_{k}^{p}(M\times
(0,T))$.

For any two points $x,y\in M$, the Carnot-Carath\'{e}odory
distance is defined by
$$
\aligned d_{C}(x,y)=\inf \{& L(\gamma )\ |\ \gamma
:[0,T]\rightarrow
M\text{ is a horizontal } C^1 \text { curve with }\\
& \gamma (0)=x,\gamma (T)=y\}\endaligned
$$
where $L(\gamma )$ denotes the length of $\gamma$ defined by the
Webster metric $ g_{\theta }$. The parabolic
Carnot-Carath\'{e}dory distance on $M\times (0,\infty )$ is
defined by (cf. [BB])
$$
d_{P}\big( (x,t),(y,s\big) )=\sqrt{d_{C}(x,y)+|t-s|}.
$$
We now define the parabolic H\"{o}lder spaces adapted to the
operator $\bigtriangleup _{H}-\frac{\partial }{\partial t}$. Let
$\Omega \subset U_{T}\ $be any open subset. For any integer $k\geq
0$ and any $\alpha >0$, let
$$
\aligned & C_{P}^{k,\alpha }(\Omega )=\left\{ u:\Omega \rightarrow
R:\Vert u\Vert_{C_{P}^{k,\alpha }}<\infty \right\} , \\
&\Vert u\Vert _{C_{P}^{k,\alpha }(\Omega
)}=\sum_{w(A_{1},...,A_{l})\leq k}\Vert X_{A_{1}}\cdots
X_{A_{l}}u\Vert_{C_{P}^{\alpha }(\Omega )},\\
&\Vert u\Vert _{C_{P}^{\alpha }(\Omega )}=|u|_{C_{P}^{\alpha
}(\Omega)}+\Vert u\Vert _{L^{\infty }(\Omega )}, \\
&|u|_{C_{P}^{\alpha }(\Omega )}=\sup \left\{
\frac{|u(x,t)-u(y,s)|}{ d_{P}\big((x,t),(y,s)\big)^{\alpha
}}:(x,t),(y,s)\in \Omega ,(x,t)\neq (y,s)\right\},
\endaligned
$$
where $0\leq A_{s}\leq 2m$, $1\leq s\leq l$. Similarly one may use
a $ C^{\infty }$ partition of unity to define the function space $
C_{P}^{k,\alpha }(M\times \lbrack T_{1},T_{2}])$ for any $
[T_{1},T_{2}]\subset (0,\infty )$. Let $d(x,y)$ be the Riemannian
distance of $x$ and $y$ in $(M,g_{\theta })$. For the relatively
compact open set $U\subset M$, there exist positive constants
$c_{1}$, $c_{2}$ depending on $U$ such that (cf. [NSW])
$$
c_{1}d(x,y)\leq d_{C}(x,y)\leq c_{2}d(x,y)^{1/2}
$$
for any $x,y\in U$. We have a natural Riemannian distance
$$\widehat{d}((x,t),(y,s))=\sqrt{d(x,y)+|t-s|^{2}}$$ on $M\times
(0,\infty)$. Using the distance $\widehat{d}$, one may define the
usual H\"{o}lder space $C^{k,\alpha}(\Omega)$ for any open subset
$\Omega \subset U_{T}$. Clearly there exist two positive constants
$C_{1}$, $C_{2}$ such that
$$
C_{1}\widehat{d}((x,t),(y,s))\leq d_{P}((x,t),(y,s))\leq
C_{2}\widehat{d}((x,t),(y,s))^{\frac{1}{2}}
$$
for any $x,y\in U$ and $0\leq |t-s|<<1$. This implies the
following relations between the parabolic H\"{o}rmander H\"{o}lder
spaces and the usual H\"{o}lder spaces
$$
C^{\alpha }(\Omega )\subset C_{P}^{\alpha }(\Omega )\subset
C^{\frac{\alpha }{2}}(\Omega ),\quad C^{k,\alpha }(\Omega )\subset
C_{P}^{k,\alpha }(\Omega ),\quad C_{P}^{2k,\alpha }(\Omega
)\subset C^{k,\frac{\alpha }{2}}(\Omega ).\tag{7.7}
$$

\proclaim{Proposition 7.1} (cf. [RS, Theorem 18], [BB, Theorem
1.1]) Let $U_{T}=U\times (0,T)$ ($T>0$) and let $\Omega \Subset
U_{T}$ be a relatively compact open subset of $U_{T}$. Suppose $u$
is locally in $L^{p}(U_{T})$, and $(\bigtriangleup
_{H}-\frac{\partial }{\partial t})u=v$.

a) If $v\in S_{k}^{p}(U_{T})$, then $\chi u\in S_{k+2}^{p}(U_{T})$
for any $\chi \in C_{0}^{\infty }(U_{T})$. In particular, there
exists a constant $c>0 $ such that
$$
\Vert u\Vert _{S_{k+2}^{p}(\Omega )}\leq c\left( \Vert u\Vert
_{L^{p}(U_{T})}+\Vert v\Vert _{S_{k}^{p}(U_{T})}\right) .
$$

b) If $v\in C_{P}^{k,\alpha }(U_{T})$, then there exists a
constant $c$ such that
$$
\Vert u\Vert _{C_{P}^{k+2,\alpha }(\Omega )}\leq c\left\{ \Vert
v\Vert _{C_{P}^{k,\alpha }(U_{T})}+\Vert u\Vert _{L^{\infty
}(U_{T})}\right\} .
$$
\endproclaim
\remark{Remark 7.1}
\newline (i) It is known that if $kp$ is large enough, then the Sobolev type space $
S_{k}^{p}$ is contained in some H\"{o}lder space (cf. [RS], [FS],
[DT], [FGN]). For example, let $k=2$ and $p>2n+4$. If $u\in
S_{2}^{p}(U_{T})$ , then for any $\chi \in C_{0}^{\infty
}(U_{T})$, we have $\chi u\in \Omega _{P}^{1,\alpha }(U_{T})$ with
$\alpha =1-\frac{2n+4}{p}$. In particular,  for any relatively
compact open subset $\Omega$ of $U_{T}$, there exists a positive
constant $c$ such that $\Vert u\Vert _{C_{P}^{1,\alpha }(\Omega
)}\leq c\Vert u\Vert _{S_{2}^{p}(U_{T})\text{. }}$
\newline (ii) Combining (7.7) and Proposition 7.1(b), we have
$$
\Vert u\Vert _{C^{l+1,\frac{\alpha }{2}}(\Omega )}\leq C\{\Vert
v\Vert _{C_{P}^{2l,\alpha }(U_{T})}+\Vert u\Vert _{L^{\infty
}(U_{T})}\}.
$$
\endremark

Since the linearization of (7.6) is a linear subelliptic parabolic
system, the short time existence and uniqueness of solution to
(7.6) follow from a standard argument. By Proposition 7.1 and a
bootstrapping argument, one can always assume that the short-time
solution $u$ of (7.6) (or (7.5)) is smooth on $M\times \lbrack
0,T)$ for some $T>0$.

\proclaim{Lemma 7.2}Let $M$ be a compact pseudo-Hermitian manifold
and let $N$ be a Sasakian manifold. For any $0<T\leq \infty $, if
$f\in C^\infty (M\times [0,T);N)$ solves (7.5), then
$$
E_{H,\widetilde{H}}(f_t)+\int_0^t\int_M|\tau
_{H,\widetilde{H}}(f_s)|^2dv_\theta ds=E_{H,\widetilde{H}}(h)
$$
for any $t\in [0,T)$. In particular, the energy
$E_{H,\widetilde{H}}$ decays along the flow.
\endproclaim
\demo{Proof} By Proposition 3.1, we get
$$
\frac{dE_{H,\widetilde{H}}(f_s)}{ds}=-\int_M\langle
\frac{\partial f}{\partial s},\tau _{H,\widetilde{H}}(f_s)\rangle
dv_\theta
$$
Therefore (7.4) and (7.5) imply that
$$
\frac{dE_{H,\widetilde{H}}(f_s)}{ds}=-\int_M|\tau
_{H,\widetilde{H}}(f_s)|^2dv_\theta
$$
Integrating the above equality over $[0,t]$ then proves this
lemma.\qed
\enddemo

Let $f:M\times [0,T)\rightarrow N$ be a $C^\infty $ solution of
(7.5). In terms of Lemma 2.1, (1.10) and the assumption that $N$
is Sasakian, we have
$$
\aligned \langle \widetilde{\nabla }\tau
_{H,\widetilde{H}},df_{H,\widetilde{H}}\rangle &=\big\langle
\widetilde{\nabla }(\tau _{H,\widetilde{H}}+\tau
_{H,\widetilde{L}}),df_{H,\widetilde{H}}\big\rangle\\
&=\sum_{A=1}^{2m}\big\langle \widetilde{\nabla }_{e_A}df(\frac
\partial {\partial t}),df_{H,\widetilde{H}}(e_A)\big\rangle
\\&=\sum_{A=1}^{2m}\big\langle \widetilde{\nabla }_{\frac \partial
{\partial t}}df(e_A)+\widetilde{T}_{\widetilde{\nabla
}}\big(df(e_A),df(\frac \partial {\partial
t})\big),df_{H,\widetilde{H}}(e_A)\big\rangle\\
&=\frac \partial {\partial t}e_{H,\widetilde{H}}.\endaligned
$$
Using (2.17), we get from (4.16) that
$$
\aligned &(\bigtriangleup_H-\frac
\partial {\partial t})e_{H,\widetilde{H}}\\
&=|\beta_{H,\widetilde{H}}|^2-2i(f_j^{\overline{\alpha
}}f_{0\overline{j}}^\alpha +f_j^\alpha
f_{0\overline{j}}^{\overline{\alpha }}-f_{\overline{j}}^\alpha
f_{0j}^{\overline{\alpha }}-f_{\overline{j}}^{\overline{\alpha
}}f_{0j}^\alpha)+(2m-4)i(f_{\overline{j}}^\alpha
f_{\overline{k}}^{\overline{\alpha }}A_{jk}-f_j^\alpha
f_k^{\overline{\alpha }}A_{\overline{j}\overline{k}})\\
&+2\big\langle
df_{H,\widetilde{H}}(Ric(\eta_j)),df_{H,\widetilde{H}}(\eta
_{\overline{j}})\big\rangle-2\widetilde{R}\big(df_{H,\widetilde{H}}(\eta_j),df_{H,\widetilde{H}}(\eta
_k),df_{H,\widetilde{H}}(\eta_{\overline{j}}),df_{H,\widetilde{H}}(\eta_{\overline{k}})\big)\\
&-2\widetilde{R}\big(df_{H,\widetilde{H}}(\eta_j),df_{H,\widetilde{H}}(\eta
_{\overline{k}}),df_{H,\widetilde{H}}(\eta_{\overline{j}}),df_{H,\widetilde{H}}(\eta _k)\big)\\
&=|\beta _{H,\widetilde{H}}|^2+\beta _{L\times
H,\widetilde{H}}*df_{H,\widetilde{H}}+A*(df_{H,\widetilde{H}})*(df_{H,\widetilde{H}})\\
&+2\big\langle
df_{H,\widetilde{H}}(Ric(\eta_j),df_{H,\widetilde{H}}(\eta_{\overline{j}})\big\rangle-2\widetilde{R}\big(df_{H,\widetilde{H}}(\eta_j),df_{H,\widetilde{H}}(\eta
_k),df_{H,\widetilde{H}}(\eta_{\overline{j}}),df_{H,\widetilde{H}}(\eta
_{\overline{k}})\big)\\
&-2\widetilde{R}\big(df_{H,\widetilde{H}}(\eta_j),df_{H,\widetilde{H}}(\eta_{\overline{k}}),df_{H,\widetilde{H}}(\eta
_{\overline{j}}),df_{H,\widetilde{H}}(\eta_k)\big),\endaligned\tag{7.8}
$$
where the notation $\Phi \ast \Psi $ denotes some contraction of
two tensors $\Phi $ and $\Psi$.

Now we want to compute $(\bigtriangleup_H-\frac \partial {\partial
t})e_{L,\widetilde{H}}$. In terms of (2.17), (2.35), (2.41), we
deduce that
$$
\aligned \bigtriangleup_H e_{L,\widetilde{H}}&=(f_0^\alpha
f_0^{\overline{\alpha }})_{k\overline{k}}+(f_0^\alpha
f_0^{\overline{\alpha }})_{\overline{k}k}\\
&=(f_{0k}^\alpha f_0^{\overline{\alpha }}+f_0^\alpha
f_{0k}^{\overline{\alpha
}})_{\overline{k}}+(f_{0\overline{k}}^\alpha f_0^{\overline{\alpha
}}+f_0^\alpha f_{0\overline{k}}^{\overline{\alpha
}})_k\\
&=2(f_{0k}^\alpha f_{0\overline{k}}^{\overline{\alpha
}}+f_{0\overline{k}}^\alpha f_{0k}^{\overline{\alpha
}})+f_0^{\overline{\alpha }}f_{0k\overline{k}}^\alpha
+f_0^{\overline{\alpha }}f_{0\overline{k}k}^\alpha +f_0^\alpha
f_{0\overline{k}k}^{\overline{\alpha }}+f_0^\alpha
f_{0k\overline{k}}^{\overline{\alpha }}\\
&=2(f_{0k}^\alpha f_{0\overline{k}}^{\overline{\alpha
}}+f_{0\overline{k}}^\alpha f_{0k}^{\overline{\alpha
}})+f_0^{\overline{\alpha }}f_{k\overline{k}0}^\alpha
+f_0^{\overline{\alpha }}f_{\overline{k}k0}^\alpha +f_0^\alpha
f_{\overline{k}k0}^{\overline{\alpha }}+f_0^\alpha
f_{k\overline{k}0}^{\overline{\alpha }}\\
&+f_0^{\overline{\alpha }}f_k^\beta \widehat{R}_{\beta \gamma
\overline{\delta }}^\alpha (f_{\overline{k}}^\gamma
f_0^{\overline{\delta }}-f_0^\gamma
f_{\overline{k}}^{\overline{\delta }})+f_0^{\overline{\alpha
}}f_{\overline{k}}^\beta \widehat{R}_{\beta \gamma
\overline{\delta }}^\alpha (f_k^\gamma f_0^{\overline{\delta
}}-f_0^\gamma f_k^{\overline{\delta }})\\
&+f_0^\alpha f_{\overline{k}}^{\overline{\beta
}}\widehat{R}_{\overline{\beta }\overline{\gamma }\delta
}^{\overline{\alpha }}(f_k^{\overline{\gamma }}f_0^\delta
-f_0^{\overline{\gamma }}f_k^\delta )+f_0^\alpha
f_k^{\overline{\beta }}\widehat{R}_{\overline{\beta
}\overline{\gamma }\delta }^{\overline{\alpha
}}(f_{\overline{k}}^{\overline{\gamma }}f_0^\delta
-f_0^{\overline{\gamma }}f_{\overline{k}}^\delta )\\
&+A_{\overline{k}}^jf_0^{\overline{\alpha }}f_{kj}^\alpha
+W_{k\overline{k}}^tf_0^{\overline{\alpha }}f_t^\alpha
+A_k^{\overline{j}}f_0^{\overline{\alpha
}}f_{\overline{j}\overline{k}}^\alpha
+A_{k,\overline{k}}^{\overline{j}}f_0^{\overline{\alpha
}}f_{\overline{j}}^\alpha\\
&+A_k^{\overline{j}}f_0^{\overline{\alpha
}}f_{\overline{k}\overline{j}}^\alpha
+W_{\overline{k}k}^{\overline{t}}f_0^{\overline{\alpha
}}f_{\overline{t}}^\alpha +A_{\overline{k}}^jf_0^{\overline{\alpha
}}f_{jk}^\alpha +A_{\overline{k},k}^jf_0^{\overline{\alpha
}}f_j^\alpha\\
&+A_k^{\overline{j}}f_0^\alpha
f_{\overline{k}\overline{j}}^{\overline{\alpha
}}+W_{\overline{k}k}^{\overline{t}}f_0^\alpha
f_{\overline{t}}^{\overline{\alpha }}+A_{\overline{k}}^jf_0^\alpha
f_{jk}^{\overline{\alpha }}+A_{\overline{k},k}^jf_0^\alpha
f_j^{\overline{\alpha }}\\
&+A_{\overline{k}}^jf_0^\alpha f_{kj}^{\overline{\alpha
}}+W_{k\overline{k}}^tf_0^\alpha f_t^{\overline{\alpha
}}+A_k^{\overline{j}}f_0^\alpha
f_{\overline{j}\overline{k}}^{\overline{\alpha
}}+A_{k,\overline{k}}^{\overline{j}}f_0^\alpha
f_{\overline{j}}^{\overline{\alpha }}.
\endaligned
$$
Consequently
$$
\aligned &\bigtriangleup_H e_{L,\widetilde{H}}\\
&=|\beta _{L\times
H,\widetilde{H}}|^2+\big\langle\widetilde{\nabla }_{\xi}\tau
_H,df_{L,\widetilde{H}}(\xi)\big\rangle-2\widetilde{R}\big(df_{L,\widetilde{H}}(\xi
),df_{H,\widetilde{H}}(\eta _k),df_{L,\widetilde{H}}(\xi
),df_{H,\widetilde{H}}(\eta_{\overline{k}})\big)\\
&+A*(df_{L,\widetilde{H}})*\beta_{H,\widetilde{H}}+\nabla
A*(df_{L,\widetilde{H}})*(df_{H,\widetilde{H}}).\endaligned\tag{7.9}
$$
Using Lemma 2.1 and (7.4), we get
$$
\aligned \widetilde{\nabla }_{\xi}\tau _H&=\widetilde{\nabla
}_{\xi}df(\frac\partial {\partial t})\\
&=\widetilde{\nabla }_{\frac\partial {\partial t}}df(\xi
)+T_{\widetilde{\nabla }}\big(df(\xi ),df(\frac\partial {\partial
t})\big).\endaligned \tag{7.10}
$$
From (7.9), (7.10), (1.10) and the assumption that $N$ is
Sasakian, we obtain
$$
\aligned (\bigtriangleup_H-\frac \partial {\partial
t})e_{L,\widetilde{H}}&=|\beta _{L\times
H,\widetilde{H}}|^2-2\widetilde{R}\big(df_{L,\widetilde{H}}(\xi
),df_{H,\widetilde{H}}(\eta _k),df_{L,\widetilde{H}}(\xi
),df_{H,\widetilde{H}}(\eta
_{\overline{k}})\big)\\
&+A*(df_{L,\widetilde{H}})*\beta _{H,\widetilde{H}}+\nabla
A*(df_{L,\widetilde{H}})*(df_{H,\widetilde{H}}).\endaligned\tag{7.11}
$$

\proclaim{Lemma 7.3} If $N$ has non-positive horizontal sectional
curvature, then
$$\widetilde{R}\big(df_{L,\widetilde{H}}(\xi
),df_{H,\widetilde{H}}(\eta _k),df_{L,\widetilde{H}}(\xi
),df_{H,\widetilde{H}}(\eta _{\overline{k}})\big)\leq 0.
$$
\endproclaim
\demo{Proof}Write $df_{H,\widetilde{H}}(\eta _k)=X_k+iY_k$
($k=1,...,m$). Since $\widetilde{K}^H\leq 0$, we find
$$
\aligned
&\widetilde{R}\big(df_{L,\widetilde{H}}(\xi),X_k+iY_k,df_{L,\widetilde{H}}(\xi
),X_k-iY_k\big)\\
&=\widetilde{R}\big(df_{L,\widetilde{H}}(\xi
),X_k,df_{L,\widetilde{H}}(\xi),X_k)+\widetilde{R}(df_{L,\widetilde{H}}(\xi
),Y_k,df_{L,\widetilde{H}}(\xi ),Y_k\big)\\
&\leq 0.\endaligned
$$\qed
\enddemo

\proclaim{Lemma 7.4} Let $N$ be a Sasakian manifold with
$\widetilde{K}^H\leq 0$ and let $f\in C^\infty (M\times [0,T),N)$
be a solution of (7.5). Set
$e_{\widetilde{H}}=e_{H,\widetilde{H}}+e_{L,\widetilde{H}}$. Then
$$
(\bigtriangleup_H-\frac \partial {\partial
t})e_{\widetilde{H}}\geq -Ce_{\widetilde{H}}.
$$
Here $C$ is a positive constant depending only on the
pseudo-Hermitian Ricci curvature and the torsion of
$(M,H(M),J,\theta )$.
\endproclaim
\demo{Proof} Utilizing (7.8), (7.11), Lemma 7.3 and Cauchy-Schwarz
inequality, we deduce that
$$
\aligned &(\bigtriangleup_H-\frac \partial {\partial
t})e_{\widetilde{H}}\\
&\geq |\beta _{H,\widetilde{H}}|^2+|\beta _{L\times
H,\widetilde{H}}|^2+\beta _{L\times
H,\widetilde{H}}*df_{H,\widetilde{H}}+A*(df_{H,\widetilde{H}})*(df_{H,\widetilde{H}})\\
&+2\big\langle df_{H,\widetilde{H}}(Ric(\eta
_j),df_{H,\widetilde{H}}(\eta
_{\overline{j}})\big\rangle+A*(df_{L,\widetilde{H}})*\beta
_{H,\widetilde{H}}+\nabla A*(df_{L,\widetilde{H}})*(df_{H,\widetilde{H}})\\
&\geq (1-\frac 12\varepsilon _1)|\beta
_{H,\widetilde{H}}|^2+(1-\frac 12\varepsilon _1)|\beta _{L\times
H,\widetilde{H}}|^2-\frac{C_1}{\varepsilon
_1}(e_{H,\widetilde{H}}+e_{L,\widetilde{H}})\\
&\geq (1-\frac 12\varepsilon _1)|\beta
_{H,\widetilde{H}}|^2+(1-\frac 12\varepsilon _1)|\beta _{L\times
H,\widetilde{H}}|^2-\frac{C_1}{\varepsilon_1}e_{\widetilde{H}}\endaligned
$$
for any $\varepsilon _1>0$, where $C_1$ is a positive constant
depending only on $Ric$, $A$ and $\nabla A$. Taking $\varepsilon
_1=2$ in the above inequality, we prove this lemma. \qed
\enddemo

In order to estimate $e_{\widetilde{H}}$, let us recall Moser's
Harnack inequality ([Mo]). For any $z_0=(x_0,t_0)\in M\times
(0,T)$, let $0<\delta <inj(M)$ (the injectivity radius), $0<\sigma
<t_0$ and let $R(z_0,\delta ,\sigma )$ be the following cylinder
$$
R(z_0,\delta ,\sigma )=\big\{(x,t)\in M\times [0,\infty
):d(x,x_0)<\delta ,\  t_0-\sigma <t<t_0\big\}
$$
where $d$ denotes the distance function of the Webster metric
$g_\theta$.

\proclaim{Lemma 7.5} Let $u$ be a non-negative smooth solution of
$$
(\bigtriangleup_H-\frac \partial {\partial t})u\geq 0
$$
on $M$. Then
$$
u(z_0)\leq C(m,\delta ,\sigma )\int_{R(z_0,\delta ,\sigma
)}u(x,t)dv_\theta dt,
$$
where $C$ is a positive constant depending only on $m,\delta $ and
$\sigma $.
\endproclaim

Since $M$ is compact, it follows from Lemma 7.5 that
$$
u(z_0)\leq C(m,\sigma )\int_{t_0-\sigma
}^{t_0}\int_Mu(x,t)dv_\theta dt\tag{7.12}
$$
for any $z_0=(x_0,t_0)\in M\times [\sigma ,T)$.

Henceforth in this section, we assume that $M$ is also Sasakian.

\proclaim{Lemma 7.6}Let $M$ be a compact Sasakian manifold and let
$N$ be a Sasakian manifold with $\widetilde{K}^H\leq 0$. Suppose
$f\in C^\infty (M\times [0,T),N)$ is a solution of (7.5). Then the
energy $E_{L,\widetilde{H}}(f(t))$ is decreasing in $t$. In
particular, if the initial map $h$ is foliated, then $f(t)$ is
foliated for each $t\in [0,T)$.
\endproclaim
\demo{Proof} From (7.11), the divergence theorem and Lemma 7.3, we
have
$$
\aligned
\frac d{dt}E_{L,\widetilde{H}}(f(t))&=\int_M\frac
\partial {\partial t}\big(e_{L,\widetilde{H}}(f)\big)dv_\theta\\
&\leq \int_M\big\{-|\beta _{L\times
H,\widetilde{H}}|^2+2\widetilde{R}\big(df_{L,\widetilde{H}}(\xi),df_{H,\widetilde{H}}(\eta
_k),df_{L,\widetilde{H}}(\xi),df_{H,\widetilde{H}}(\eta
_{\overline{k}})\big)\big\}dv_\theta\\
&\leq 0.\endaligned
$$
\qed
\enddemo

\proclaim{Lemma 7.7}Let $M$ be a compact Sasakian manifold and let
$N$ be a Sasakian manifold with $\widetilde{K}^H\leq 0$. Suppose
$f\in C^\infty (M\times [0,T),N)$ ($0<T\leq \infty $) is a
solution of (7.5). Then $e_{\widetilde{H}}(f)$ is uniformly
bounded.
\endproclaim
\demo{Proof} From Lemma 7.4, we know that $e_{\widetilde{H}}$
satisfies
$$
(\bigtriangleup_H-\frac
\partial {\partial t})e_{\widetilde{H}}\geq
-Ce_{\widetilde{H}}
$$
for some constant $C$. Let
$$F(x,t):=e^{-Ct}e_{\widetilde{H}}(f),\quad (x,t)\in
M\times [0,T).
$$
It follows that
$$
(\bigtriangleup_H-\frac \partial {\partial t})F(x,t)\geq 0.
$$
Let $0<\sigma <T$. Then for any $z_0=(x_0,t_0)\in M\times
[\sigma,T)$, (7.12) implies that
$$
\aligned e_{\widetilde{H}}(f)(z_0)&\leq
C_1e^{Ct_0}\int_{t_0-\sigma}^{t_0}\int_Me^{-Ct}e_{\widetilde{H}}(f)dv_\theta dt\\
&\leq
C_1e^{C\sigma}\int_{t_0-\sigma}^{t_0}\int_Me_{\widetilde{H}}(f)dv_\theta dt\\
&\leq C_2\int_{t_0-\sigma}^{t_0}E_{\widetilde{H}}(f(t))dt\\
&\leq C_2E_{\widetilde{H}}(h)\endaligned
$$
since $E_{H,\widetilde{H}}(f(t))$ and $E_{L,\widetilde{H}}(f(t))$
are decreasing in $t$ in view of Lemmas 7.2, 7.6. \qed
\enddemo

Next we want to derive Bochner formulas for
$e_{H,\widetilde{L}}(f_t)$ and $e_{L,\widetilde{L}}(f_t)$.
According to (7.1) and the definition of $\bigtriangleup_H$, one
has
$$
\aligned
\bigtriangleup_He_{H,\widetilde{L}}&=(f_j^0f_{\overline{j}}^0)_{k\overline{k}}+(f_j^0f_{\overline{j}}^0)_{\overline{k}k}\\
&=2(f_{jk}^0f_{\overline{j}\overline{k}}^0+f_{j\overline{k}}^0f_{\overline{j}k}^0)
+f_{\overline{j}}^0f_{jk\overline{k}}^0+f_j^0f_{\overline{j}k\overline{k}}^0
+f_{\overline{j}}^0f_{j\overline{k}k}^0+f_j^0f_{\overline{j}\overline{k}k}^0
\endaligned \tag{7.13}
$$
and
$$
\aligned \bigtriangleup_H e_{L,\widetilde{L}}&=\frac
12\big((f_0^0)_{k\overline{k}}^2+(f_0^0)_{\overline{k}k}^2\big)\\
&=2f_{0k}^0f_{0\overline{k}}^0+f_0^0(f_{0k\overline{k}}^0+f_{0\overline{k}k}^0)
\endaligned
\tag{7.14}
$$
Using (2.14) and (2.24), we deduce from (7.16) and (7.17) that
$$
\aligned \bigtriangleup_H e_{H,\widetilde{L}}
&=2(f_{jk}^0f_{\overline{j}\overline{k}}^0+f_{j\overline{k}}^0f_{\overline{j}k}^0)
+f_{\overline{j}}^0f_{k\overline{k}j}^0+f_j^0f_{k\overline{k}\overline{j}}^0
+f_{\overline{j}}^0f_{\overline{k}kj}^0+f_j^0f_{\overline{k}k\overline{j}}^0\\
&+f_{\overline{j}}^0f_{j0}^0+f_j^0f_{\overline{j}0}^0+if_{\overline{j}}^0f_{0j}^0-if_j^0f_{0\overline{j}}^0
+f_{\overline{j}}^0f_t^0R_{kj\overline{k}}^t+f_j^0f_{\overline{t}}^0R_{\overline{k}\overline{j}k}^{\overline{t}}\\
&+if_{\overline{j}}^0\big(f_j^{\overline{\alpha }}f_k^\alpha
-f_j^\alpha f_k^{\overline{\alpha
}}\big)_{\overline{k}}-if_j^0\big(f_k^{\overline{\alpha
}}f_{\overline{j}}^\alpha -f_k^\alpha
f_{\overline{j}}^{\overline{\alpha }}\big)_{\overline{k}}\\
&+if_{\overline{j}}^0\big(f_{\overline{k}}^\alpha
f_j^{\overline{\alpha }}-f_{\overline{k}}^{\overline{\alpha
}}f_j^\alpha \big)_k-if_j^0\big(f_{\overline{j}}^\alpha
f_{\overline{k}}^{\overline{\alpha
}}-f_{\overline{j}}^{\overline{\alpha }}f_{\overline{k}}^\alpha\big)_k\\
&=|\beta _{H\times H,\widetilde{L}}|^2+\langle \widetilde{\nabla
}\tau _{H,\widetilde{L}},df_{H,\widetilde{L}}\rangle +\beta
_{H\times L,\widetilde{L}}*df_{H,\widetilde{L}}+\beta _{L\times
H,\widetilde{L}}*df_{H,\widetilde{L}}\\
&+\big\langle df_{H,\widetilde{L}}(Ric_H(\eta
_j),df_{H,\widetilde{L}}(\eta _{\overline{j}})\big\rangle
+\big\langle df_{H,\widetilde{L}}(Ric_H(\eta
_{\overline{j}}),df_{H,\widetilde{L}}(\eta _j)\big\rangle\\
&+\beta
_{H,\widetilde{H}}*(df_{H,\widetilde{H}})*(df_{H,\widetilde{L}})
\endaligned\tag{7.15}
$$
and
$$
\aligned \bigtriangleup _H
e_{L,\widetilde{L}}&=2f_{0k}^0f_{0\overline{k}}^0+f_0^0f_{k\overline{k}0}^0+f_0^0f_{\overline{k}k0}^0+if_0^0\big((f_0^{\overline{\alpha
}}f_{\overline{k}}^\alpha -f_0^\alpha
f_{\overline{k}}^{\overline{\alpha }})_k-(f_0^\alpha
f_k^{\overline{\alpha }}-f_0^{\overline{\alpha }}f_k^\alpha
)_{\overline{k}}\big)\\
&=|\beta _{L\times H,\widetilde{L}}|^2+\langle \widetilde{\nabla
}_\xi \tau _{H,\widetilde{L}},df_{L,\widetilde{L}}(\xi )\rangle
+\beta _{L\times
H,\widetilde{H}}*df_{H,\widetilde{H}}*df_{L,\widetilde{L}}\\
&+\beta_{H,\widetilde{H}}*df_{L,\widetilde{L}}*df_{L,\widetilde{H}}.
\endaligned\tag{7.16}
$$
In terms of (1.7), (1.9), (1.10), Lemma 2.1 and (7.5), we have
$$
\langle \widetilde{\nabla }\tau
_{H,\widetilde{L}},df_{H,\widetilde{L}}\rangle =\frac \partial
{\partial t}e_{H,\widetilde{L}}+\sum_{A=1}^{2m}\big\langle
\widetilde{J}df_{H,\widetilde{H}}(e_A),\tau
_{H,\widetilde{H}}(f)\big\rangle \widetilde{\theta
}\big(df_{H,\widetilde{L}}(e_A)\big)
$$
and
$$
\langle \widetilde{\nabla }_\xi \tau
_{H,\widetilde{L}},df_{L,\widetilde{L}}(\xi )\rangle =\frac
\partial {\partial t}e_{L,\widetilde{L}}+\big\langle
\widetilde{J}df_{L,\widetilde{H}}(\xi ),\tau
_{H,\widetilde{H}}(f)\big\rangle \widetilde{\theta
}\big(df_{L,\widetilde{L}}(\xi )\big).
$$
In addition, (2.14) implies
$$
\beta _{H\times L,\widetilde{L}}*df_{H,\widetilde{L}}=\beta
_{L\times
H,\widetilde{L}}*df_{H,\widetilde{L}}+df_{H,\widetilde{H}}*df_{L,\widetilde{H}}*df_{H,\widetilde{L}}.
$$
Consequently
$$
\aligned (\bigtriangleup_H-\frac \partial {\partial
t})e_{H,\widetilde{L}}&=|\beta _{H\times
H,\widetilde{L}}|^2+df_{H,\widetilde{H}}*\tau
_{H,\widetilde{H}}*df_{H,\widetilde{L}}+\beta _{L\times
H,\widetilde{L}}*df_{H,\widetilde{L}}\\
&+df_{H,\widetilde{H}}*df_{L,\widetilde{H}}*df_{H,\widetilde{L}}+\beta
_{H,\widetilde{H}}*(df_{H,\widetilde{H}})*(df_{H,\widetilde{L}})\\
&+\big\langle df_{H,\widetilde{L}}(Ric_H(\eta
_j),df_{H,\widetilde{L}}(\eta
_{\overline{j}})\big\rangle+\big\langle
df_{H,\widetilde{L}}(Ric_H(\eta
_{\overline{j}}),df_{H,\widetilde{L}}(\eta _j)\big\rangle
\endaligned\tag{7.17}
$$
and
$$
\aligned (\bigtriangleup_H-\frac \partial {\partial
t})e_{L,\widetilde{L}}&=|\beta _{L\times
H,\widetilde{L}}|^2+df_{L,\widetilde{H}}*\tau
_{H,\widetilde{H}}*df_{L,\widetilde{L}}+\beta _{L\times
H,\widetilde{H}}*df_{H,\widetilde{H}}*df_{L,\widetilde{L}}\\
&+\beta
_{H,\widetilde{H}}*df_{L,\widetilde{L}}*df_{L,\widetilde{H}}.
\endaligned\tag{7.18}
$$

Clearly the usual energy density of the map $f$ is given by
$$e(f)=e_{H,\widetilde{H}}(f)+e_{L,\widetilde{H}}(f)+e_{H,\widetilde{L}}(f)+e_{L,\widetilde{L}}(f)
$$
and thus
$E(f)=E_{H,\widetilde{H}}(f)+E_{L,\widetilde{H}}(f)+E_{H,\widetilde{L}}(f)+E_{L,\widetilde{L}}(f)$.

\proclaim{Lemma 7.8} Let $M$ be a compact Sasakian manifold and
$N$ be a Sasakian manifold with $\widetilde{K}^H\leq 0$. Suppose
$f\in C^\infty (M\times [0,T),N)$ is a solution of (7.5). Then
$$
e(f_t)\leq C(\sigma )\sup_{[t-\sigma ,t]}E(f_t)
$$ for $t\in[\sigma ,T)$.
\endproclaim
\demo{Proof} From the proof of Lemma 7.4, we have
$$
\aligned (\bigtriangleup_H-\frac \partial {\partial t})\left(
e_{H,\widetilde{H}}(f_t)+e_{L,\widetilde{H}}(f_t)\right)& \geq
(1-\frac 12\varepsilon _1)|\beta _{H,\widetilde{H}}|^2+(1-\frac
12\varepsilon _1)|\beta _{L\times
H,\widetilde{H}}|^2\\
&-\frac{C_1}{\varepsilon _1}\big(
e_{H,\widetilde{H}}(f_t)+e_{L,\widetilde{H}}(f_t)\big)\endaligned\tag{7.19}
$$
Utilizing Lemma 7.7 and Cauchy-Schwarz inequality, we derive from
(7.17) and (7.18) that
$$
\aligned &(\bigtriangleup_H-\frac \partial {\partial t})\big(
e_{H,\widetilde{L}}(f_t)+e_{L,\widetilde{L}}(f_t)\big)\\
&=|\beta _{H\times H,\widetilde{L}}|^2+df_{H,\widetilde{H}}*\tau
_{H,\widetilde{H}}*df_{H,\widetilde{L}}+\beta _{L\times
H,\widetilde{L}}*df_{H,\widetilde{L}}+df_{H,\widetilde{H}}*df_{L,\widetilde{H}}*df_{H,\widetilde{L}}\\
&+|\beta _{L\times H,\widetilde{L}}|^2+\big\langle
df_{H,\widetilde{L}}(Ric_H(\eta _j),df_{H,\widetilde{L}}(\eta
_{\overline{j}})\big\rangle+\big\langle
df_{H,\widetilde{L}}(Ric_H(\eta_{\overline{j}}),df_{H,\widetilde{L}}(\eta _j)\big\rangle\\
&+\beta_{H,\widetilde{H}}*(df_{H,\widetilde{H}})*(df_{H,\widetilde{L}})+df_{L,\widetilde{H}}*\tau
_{H,\widetilde{H}}*df_{L,\widetilde{L}} +\beta _{L\times
H,\widetilde{H}}*df_{H,\widetilde{H}}*df_{L,\widetilde{L}}\\
&+\beta
_{H,\widetilde{H}}*df_{L,\widetilde{L}}*df_{L,\widetilde{H}}\\
&\geq -\frac 12\varepsilon _2|\beta _{H,\widetilde{H}}|^2+(1-\frac
12\varepsilon _2)|\beta _{L\times H,\widetilde{L}}|^2-\frac
12\varepsilon _2\big(
e_{H,\widetilde{H}}(f_t)+e_{L,\widetilde{H}}(f_t)\big)\\
&-\frac 12\varepsilon _2|\beta _{L\times
H,\widetilde{H}}|^2-\frac{C_2}{\varepsilon _2}\big(
e_{H,\widetilde{L}}(f_t)+e_{L,\widetilde{L}}(f_t)\big).
\endaligned\tag{7.20}
$$
Taking $\varepsilon _1=\varepsilon _2=1$, it follows from (7.19)
and (7.20) that
$$
(\bigtriangleup_H-\frac \partial {\partial t})e(f_t)\geq
-\widetilde{C}e(f_t)\tag{7.21}
$$
for some positive $\widetilde{C}$ depending only on $M,N$ and $h$.
Therefore this lemma follows immediately from (7.21) and Lemma
7.5.\qed
\enddemo

Now we want to estimate the partial energies
$E_{H,\widetilde{L}}(f_t)$ and $E_{L,\widetilde{L}}(f_t)$.

\proclaim{Lemma 7.9} Let $M$, $N$ and $f\in C^\infty (M\times
[0,T),N)$ be as in Lemma 7.8. Then
$$
\frac d{ds}E_{H,\widetilde{L}}(f_t)=-\int_M|\tau
_{H,L}(f_t)|^2dv_\theta +\int_M\big\langle J\tau
_{H,\widetilde{H}}(f_t),df_{H,\widetilde{H}}(e_A)\big\rangle
\widetilde{\theta }\big(df_{H,\widetilde{L}}(e_A)\big)dv_\theta .
$$
Furthermore, we have
$$E_{H,\widetilde{L}}(f_t)\leq C_2t+C_3
$$
for any $t\in [0,T)$ and some constants
$C_2,C_3$ depending on $M$, $N$ and $h$.
\endproclaim
\demo{Proof} By definition, $E_{H,\widetilde{L}}(f_t)$ is given by
$$
E_{H,\widetilde{L}}(f_t)=\frac 12\int_M\sum_{A=1}^{2m}\big\langle
df_{H,\widetilde{L}}(e_A),df_{H,\widetilde{L}}(e_A)\big\rangle
dv_\theta
$$
where $\{e_A\}_{A=1}^{2m}$ is an orthonormal frame field of
$H(M)$. Then, by using Lemma 2.1, (1.9), (1.10) and the divergence
theorem, we derive that
$$
\aligned &\frac d{ds}E_{H,\widetilde{L}}(f_t)\\
&=\frac 12\int_M\sum_{A=1}^{2m}\frac \partial {\partial t}\langle
df_{H,\widetilde{L}}(e_A),df_{H,\widetilde{L}}(e_A)\rangle dv_\theta\\
&=\int_M\sum_{A=1}^{2m}\langle\widetilde{\nabla }_{\frac \partial
{\partial t}}df(e_A),df_{H,\widetilde{L}}(e_A)\rangle dv_\theta \\
&=\int_M\sum_{A=1}^{2m}\big\{\langle\widetilde{\nabla
}_{e_A}df(\frac
\partial {\partial
t}),df_{H,\widetilde{L}}(e_A)\rangle+\big\langle\widetilde{T}_{\widetilde{\nabla
}}\big(df(\frac \partial {\partial
t}),df(e_A)\big),df_{H,\widetilde{L}}(e_A)\big\rangle\big\}dv_\theta \\
&=\int_M\sum_{A=1}^{2m}\big\{\langle\widetilde{\nabla }_{e_A}\tau
_H(f_t),df_{H,\widetilde{L}}(e_A)\rangle+\big\langle
d\widetilde{\theta}(\tau_{H,\widetilde{H}}(f_t),df_{H,\widetilde{H}}(e_A))\widetilde{\xi
},df_{H,\widetilde{L}}(e_A)\big\rangle\big\}dv_\theta\\
&=-\int_M|\tau _{H,L}(f_t)|^2dv_\theta +\int_M\big\langle J\tau
_{H,\widetilde{H}}(f_t),df_{H,\widetilde{H}}(e_A)\big\rangle
\widetilde{\theta }(df_{H,\widetilde{L}}(e_A))dv_\theta.
\endaligned\tag{7.22}
$$
In terms of Lemma 7.7, (7.25) and H\"older's inequality, we have
$$
\aligned \frac d{dt}E_{H,\widetilde{L}}(f_t)&\leq C\int_M|\tau
_{H,\widetilde{H}}(f_t)||df_{H,\widetilde{L}}|dv_\theta\\
&\leq \sqrt{2}C\big( \int_M|\tau
_{H,\widetilde{H}}(f_t)|^2dv_\theta \big) ^{1/2}\big(
\int_Me_{H,\widetilde{L}}(f_t)dv_\theta \big) ^{1/2}
\endaligned
$$
which implies
$$
\int_0^t\frac{dE_{H,\widetilde{L}}(f_s)}{2\sqrt{E_{H,\widetilde{L}}(f_s)}}
\leq \frac C{\sqrt{2}}\int_0^t\big( \int_M|\tau
_{H,\widetilde{H}}(f_s)|^2dv_\theta \big) ^{1/2}ds.
$$
Consequently
$$
\sqrt{E_{H,\widetilde{L}}(f_t)}-\sqrt{E_{H,\widetilde{L}}(h)}\leq
\frac C{\sqrt{2}}\big( \int_0^t\int_M|\tau
_{H,\widetilde{H}}(f_s)|^2dv_\theta \big)
^{1/2}\sqrt{t}.\tag{7.23}
$$
Thus Lemma 7.2 and (7.23) imply
$$
\sqrt{E_{H,\widetilde{L}}(f_t)}\leq \frac
C{\sqrt{2}}\sqrt{E_{H,\widetilde{H}}(h)}\sqrt{t}+\sqrt{E_{H,\widetilde{L}}(h)}.
$$
\qed
\enddemo

\proclaim{Lemma 7.10}Let $M$, $N$ and $f\in C^\infty (M\times
[0,T),N)$ be as in Lemma 7.7. Suppose $h$ is foliated. Then
$e_{L,\widetilde{L}}(f_t)$ is uniformly bounded, and
$E_{L,\widetilde{L}}(f_t)$ is decreasing in $t$.
\endproclaim
\demo{Proof} Since $h$ is foliated, we see from Lemma 7.6 that
$f_t$ is foliated for each $t\in [0,T)$. As a result, (7.18)
becomes
$$
(\bigtriangleup_H-\frac
\partial {\partial t})e_{L,\widetilde{L}}(f_t) =|\beta _{L\times
H,\widetilde{L}}|^2\geq 0.\tag{7.24}
$$
By Maximum principle, we have $e_{L,\widetilde{L}}(f_t)\leq
\sup_Me_{L,\widetilde{L}}(h)$. In terms of the divergence theorem,
(7.24) implies that
$$
\frac d{dt}\int_Me_{L,\widetilde{L}}(f_t)dv_{\theta } =\int_M\big(
\bigtriangleup_H e_{L,\widetilde{L}}(f_t)-|\beta _{L\times
H,\widetilde{L}}|^2\big) dv_\theta \leq 0.
$$
Hence $E_{L,\widetilde{L}}(f_t)$ is decreasing in $t$. \qed
\enddemo

From Lemmas 7.2, 7.6, 7.8, 7.9 and 7.10, we immediately get the
following global existence of (7.5).

\proclaim{Proposition 7.11} Let $M$ and $N$ be compact Sasakian
manifolds. Suppose $N$ has non-positive horizontal curvature and
the initial map $h:M\rightarrow N$ is foliated. Then the solution
$f$ of (7.5) exits for all $t\geq 0$.
\endproclaim

Now we are able to establish some existence results for
$(H,\widetilde{H})$-harmonic maps when the target manifold $N$ is
a compact regular Sasakian manifold, that is, $N$ can be realized
as a Riemannian submersion $\pi :(N,g_{\widetilde{\theta
}})\rightarrow (B,g_B)$ over a compact K\"ahler manifold. Let
$i(B)$ be the injectivity radius of $B$. We denote by $B_r(y)$ the
geodesic ball centered at $y$ with radius $r$ in $B$. Hence, if
$r\leq i(B)$, then any two points in $B_r(y)$ can be joined by a
unique geodesic in $B_r(y)$.

\proclaim{Theorem 7.12} Let $M$ be a compact Sasakian manifold and
let $N$ be a compact Sasakian manifold with non-positive
horizontal sectional curvature. Suppose $\pi :N\rightarrow B$ is a
Riemannian submersion over a K\"ahler manifold $B$ and
$h:M\rightarrow N$ is a given foliated map. Then there exists a
smooth foliated $(H,\widetilde{H})$-harmonic map in the same
homotopy class as $h$.
\endproclaim

\demo{Proof} From Proposition 7.11, we know that there is a global
solution $f:M\times [0,\infty )\rightarrow N$ of (7.5) with the
initial map $h$. Set $\varphi _t=\pi \circ f_t$ and $\psi =\pi
\circ h$. Since $f_t$ is foliated for each $t\in [0,\infty )$ in
view of Lemma 7.6, it follows from Proposition 3.8 that
$\varphi:M\times [0,\infty)\rightarrow N$ satisfies the following
harmonic heat flow
$$
\cases\frac{\partial \varphi }{\partial t}=\tau ^{g_B}(\varphi_t)\\
\varphi |_{t=0}=\psi.\endcases\tag{7.25}
$$
Observe from (1.25) that $B$ has non-positive sectional curvature
duo to the non-positive horizontal curvature condition on $N$. By
Eells-Sampson theorem, the solution $\varphi $ of (7.25) converges
in $C^\infty (M,B)$ to a harmonic map $\varphi _\infty $ as
$t\rightarrow \infty $. Therefore there is a sufficiently large
$T>0$ such that if $t\geq T$, then $\varphi _t(x)\in
B_{i(B)}(\varphi _\infty (x))$. In particular, $\varphi _T(x)\in
B_{i(B)}(\varphi _\infty (x))$) for each $x\in M$. Set $v(x)=\exp
_{\varphi _T(x)}^{-1}(\varphi _\infty (x))$. Clearly $v\in \Gamma
(\varphi ^{-1}(TB))$. We may define a horizontal vector field
$\widetilde{v}$ along $f_T$ such that $\widetilde{v}(x)\in
H_{f_T(x)}(N)$ and $d\pi ()=v(x)$ for $x\in M$. Set
$$
\widehat{f}(x)=\exp _{f_T(x)}^{\bot }(\widetilde{v}(x)),\quad x\in
M,
$$
where $\exp ^{\bot }$ denotes the normal exponential map along the
Reeb leaf $\pi ^{-1}(\varphi _T(x))$. Note that $\varphi _T(x)$
and $\varphi _\infty (x)$ can be joined by a unique geodesic
$\gamma _x$ lying in $B_{i(B)}(\varphi _\infty (x))$. Actually
$\widehat{f}(x)$ is the image of $f_T(x)$ under the honolomy
diffeomorphism $h_{\gamma _x}:\pi ^{-1}($ $\varphi
_T(x))\rightarrow \pi ^{-1}(\varphi _\infty (x))$ associated to
the geodesic $\gamma _x$ from $\varphi _T(x)$ to $\varphi _\infty
(x)$. Thus $\widehat{f}$ is a foliated map. Obviously the map
$\widehat{f}:M\rightarrow N$ lies in the same homotopy class as
$h$ and satisfies
$$
\pi \circ
\widehat{f}=\varphi _\infty .
$$
Thus Proposition 3.8 implies that $\widehat{f}$ is a
$(H,\widetilde{H})$-harmonic map.\qed
\enddemo

\remark{Remark 7.2}Let $V$ be any vertical vector field on $N$ and
let $\zeta _s$ denote the one parameter transformation group
generated by $V$. Then $\zeta _s\circ \widehat{f}$ is also a
$(H,\widetilde{H})$-harmonic map in the same homotopy class as
$h$. Hence we do not have the uniqueness for
$(H,\widetilde{H})$-harmonic maps in a fixed homotopy class in
general.
\endremark

\proclaim{Lemma 7.13} Let $M,N$ and $B$ be as in Theorem 7.12. Let
$f:M\times [0,\infty )\rightarrow N$ be a solution of (7.5) with
initial map $h$. Suppose $h:M\rightarrow N$ is a foliated
$(H,\widetilde{H})$-harmonic map. Then $f_t:M\rightarrow N$ is a
foliated $(H,\widetilde{H})$-harmonic map for each $t\in [0,\infty
)$.
\endproclaim
\demo{Proof} Set $\varphi _t=\pi \circ f_t$ and $\psi =\pi \circ
h$, where $\pi :N\rightarrow B$ is the Riemannian submersion. From
the proof of Theorem 7.12, we know that $\varphi _t$ satisfies the
harmonic map heat flow (7.25). Since $h$ is assumed to be a
foliated $(H,\widetilde{H})$-harmonic map, Proposition 3.8 implies
that $\psi :M\rightarrow B$ is a harmonic map, which may be also
regarded as a solution of (7.25) independent of the time $t$. By
the uniqueness for solutions of (7.25), we find that $\varphi _t$
is harmonic for each $t$. It follows from Proposition 3.8 again
that $f_t:M\rightarrow N$ is $(H,\widetilde{H})$-harmonic.\qed
\enddemo

\proclaim{Theorem 7.14} Let $M$ and $N$ be compact Sasakian
manifolds and let $h:M\rightarrow N$ be a foliated map. Suppose
$N$ is regular with non-positive horizontal sectional curvature.
Then there exists a foliated special $(H,\widetilde{H})$-harmonic
map in the same foliated homotopy class as $h$.
\endproclaim

\demo{Proof} Without loss of generality, we may assume that $h$ is
a foliated $(H,\widetilde{H})$-harmonic map in view of Theorem
7.12. Suppose $f:M\times [0,\infty )$ is a solution of (7.25) with
the initial map $h$. By Lemma 7.13 , each $f_t$ is a foliated
$(H,\widetilde{H})$-harmonic map. Then Lemma 7.9 gives
$$
\frac d{dt}E_{H,\widetilde{L}}(f_t)=-\int_M|\tau
_{H,\widetilde{L}}(f_t)|^2dv_\theta.
$$
This yields that $E_{H,\widetilde{L}}(f_t)$ is decreasing in $t$
and
$$
\int_0^\infty \int_M|\tau _{H,\widetilde{L}}(f_t)|^2dv_\theta
dt<\infty.\tag{7.26}
$$
Consequently Lemmas 7.2, 7.6, 7.8, 7.10 imply that $e(f_t)$ is
uniformly bounded. Applying Proposition 7.1 to the solution of
(7.6), we find that all higher derivatives of $f$ are uniformly
bounded. On the other hand, (7.26) implies that there exists a
sequence $\{t_k\}$ such that
$$
|\tau _{H,\widetilde{L}}(f_{t_k})|_{L^2(M)}\rightarrow 0\quad
\text{as }t_k\rightarrow \infty .\tag{7.27}
$$
In terms of the Arzela-Ascoli theorem, by passing a subsequence
$\{t_{k_l}\}$ of $\{t_k\}$, we conclude that $f(\cdot ,t_{k_l})$
converges in $C^\infty (M,N)$ to a limit $f_\infty $ (as
$t_{k_l}\rightarrow \infty $), which satisfies both $\tau
_{H,\widetilde{H}}(f_\infty)=0$ and $\tau
_{H,\widetilde{L}}(f_\infty)=0$. Clearly $f_\infty $ lies in the
same foliated homotopy class as $h$. \qed
\enddemo

\heading{\bf 8. Foliated rigidity and Siu-type strong rigidity
results}
\endheading
\vskip 0.3 true cm

First, we introduce the following

\definition{Definition 8.1}
We say that a map $f:(M^{2m+1},H(M),J,\theta )\rightarrow
(N^{2n+1}, \widetilde{H}(N),\widetilde{J},\widetilde{\theta })$
has split horizontal second fundamental form if $\pi
_{\widetilde{H}}(\beta (T,X))=0$ for any $ X\in H(M)$, that is,
$f_{0k}^\alpha =f_{0\overline{k}}^\alpha =0$ for $k=1,...,m$.
\enddefinition

\proclaim{Lemma 8.1} Let $f:M^{2m+1}\rightarrow N^{2n+1}$ be a
$(H,\widetilde{H})$-harmonic map with split horizontal second
fundamental form. Suppose that $M$ is compact and $N$ is Sasakian.
Then $f$ is foliated.
\endproclaim

\demo{Proof} By (2.17), we have
$$
f_{k\overline{k}}^\alpha =f_{\overline{k}k}^\alpha +mif_0^\alpha
\tag{8.1}
$$
It follows from Corollary 3.3 and (8.1) that
$$
f_{k\overline{k}}^\alpha =\frac{mi}2f_0^\alpha ,\quad
f_{\overline{k}k}^\alpha =- \frac{mi}2f_0^\alpha.  \tag{8.2}
$$
Let us define a global $1$-form on $M$ by $\widetilde{\rho
}=-i(f_0^\alpha f_k^{\overline{\alpha }}\theta ^k-f_0^{
\overline{\alpha }}f_{\overline{k}}^\alpha \theta
^{\overline{k}})$. Using Lemma 1.3, (8.2) and the assumption that
$f_{0k}^\alpha =f_{0\overline{k}}^\alpha =0$, we find that
$$
\aligned
0&=\int_M\delta\widetilde{\rho}\\
&=i\int_M\big\{(f_0^\alpha f_k^{\overline{\alpha
}})_{\overline{k}}-(f_0^{\overline{\alpha }}f_{\overline{k}}^\alpha )_k\big\}dv_{\theta }  \\
&=i\int_M\big\{f_{0\overline{k}}^\alpha f_k^{\overline{\alpha
}}+f_0^\alpha f_{k \overline{k}}^{\overline{\alpha
}}-f_{0k}^{\overline{\alpha }}f_{\overline{k}
}^\alpha -f_0^{\overline{\alpha }}f_{\overline{k}k}^\alpha \big\}dv_{\theta } \\
&=i\int_M\big\{\frac{mi}2f_0^\alpha f_0^{\overline{\alpha
}}+\frac{mi}2f_0^{
\overline{\alpha }}f_0^\alpha \big\}dv_{\theta } \\
&=-m\int_M\sum_\alpha |f_0^\alpha |^2dv_{\theta },
\endaligned
$$
which implies $f_0^\alpha =0$, that is, $f$ is foliated.\qed
\enddemo

\proclaim{Theorem 8.2} Let $(M^{2m+1},H(M),J,\theta )$ and
$(N^{2n+1},\widetilde{H}(N),\widetilde{J} ,\widetilde{\theta })$
be compact Sasakian manifolds and let $f:M\rightarrow N$ be a
$(H,\widetilde{H})$-harmonic map. If $(N^{2n+1},\widetilde{\nabla
})$ has non-positive horizontal sectional curvature, then $f$ is
foliated.
\endproclaim

\demo{Proof} First we define a global $1$-form $\widetilde{\rho
_1}$ by
$$
\widetilde{\rho _1}=-\{(f_{0k}^\alpha f_0^{\overline{\alpha
}}+f_{0k}^{\overline{ \alpha }}f_0^\alpha )\theta
^k+(f_{0\overline{k}}^{\overline{\alpha } }f_0^\alpha
+f_{0\overline{k}}^\alpha f_0^{\overline{\alpha }})\theta ^{
\overline{k}}\}.  \tag{8.3}
$$
Using (1.20), (1.22), (2.17), (2.35) and (2.41), we deduce that
$$
\aligned \delta\widetilde{\rho _1}=&(f_{0k}^\alpha
f_0^{\overline{\alpha }}+f_{0k}^{\overline{ \alpha }}f_0^\alpha
)_{\overline{k}}+(f_{0\overline{k}}^{\overline{\alpha }
}f_0^\alpha +f_{0\overline{k}}^\alpha f_0^{\overline{\alpha }})_k \\
=&2|f_{0k}^\alpha |^2+2|f_{0\overline{k}}^\alpha
|^2+f_0^{\overline{\alpha } }(f_{0k\overline{k}}^\alpha
+f_{0\overline{k}k}^\alpha )+f_0^\alpha (f_{0k
\overline{k}}^{\overline{\alpha
}}+f_{0\overline{k}k}^{\overline{\alpha }})
\\
=&2|f_{0k}^\alpha |^2+2|f_{0\overline{k}}^\alpha
|^2+f_0^{\overline{\alpha } }(f_{k0\overline{k}}^\alpha
+f_{\overline{k}0k}^\alpha )+f_0^\alpha (f_{k0
\overline{k}}^{\overline{\alpha
}}+f_{\overline{k}0k}^{\overline{\alpha }})
\\
=&2|f_{0k}^\alpha |^2+2|f_{0\overline{k}}^\alpha
|^2+f_0^{\overline{\alpha } }(f_{k\overline{k}0}^\alpha
+f_{\overline{k}k0}^\alpha )+f_0^\alpha (f_{k
\overline{k}0}^{\overline{\alpha
}}+f_{\overline{k}k0}^{\overline{\alpha }})
\\
&+f_0^{\overline{\alpha }}f_k^\beta \widehat{R}_{\beta \gamma
\overline{ \delta }}^\alpha (f_{\overline{k}}^\gamma
f_0^{\overline{\delta } }-f_0^\gamma
f_{\overline{k}}^{\overline{\delta }})+f_0^\alpha f_{\overline{k
}}^{\overline{\beta }}\widehat{R}_{\overline{\beta
}\overline{\gamma }\delta }^{\overline{\alpha
}}(f_k^{\overline{\gamma }}f_0^\delta -f_0^{\overline{
\gamma }}f_k^\delta ) \\
&+f_0^{\overline{\alpha }}f_{\overline{k}}^\beta
\widehat{R}_{\beta \gamma \overline{\delta }}^\alpha (f_k^\gamma
f_0^{\overline{\delta }}-f_0^\gamma f_k^{\overline{\delta
}})+f_0^\alpha f_k^{\overline{\beta }}\widehat{R}_{
\overline{\beta }\overline{\gamma }\delta }^{\overline{\alpha
}}(f_{ \overline{k}}^{\overline{\gamma }}f_0^\delta
-f_0^{\overline{\gamma }}f_{
\overline{k}}^\delta ) \\
=&2|f_{0k}^\alpha |^2+2|f_{0\overline{k}}^\alpha
|^2+f_0^{\overline{\alpha } }(f_{k\overline{k}0}^\alpha
+f_{\overline{k}k0}^\alpha )+f_0^\alpha (f_{k
\overline{k}0}^{\overline{\alpha
}}+f_{\overline{k}k0}^{\overline{\alpha }})
\\
&-\widetilde{Q}(df_{L,\widetilde{H}}(\xi)\wedge
df_{H,\widetilde{H}}(\eta _k),
\overline{df_{L,\widetilde{H}}(\xi)\wedge df_{H,\widetilde{H}}(\eta _k)}) \\
&-\widetilde{Q}(df_{L,\widetilde{H}}(\xi)\wedge
df_{H,\widetilde{H}}(\eta _{
\overline{k}}),\overline{df_{L,\widetilde{H}}(\xi)\wedge
df_{H,\widetilde{H} }(\eta _{\overline{k}})})
\endaligned \tag{8.4}
$$
where $df_{L,\widetilde{H}}(\xi)=\pi _{\widetilde{H}}(df(\xi))$.
It follows from Corollary 3.3 and (8.4) that
$$
\aligned \delta\widetilde{\rho _1}=&2|f_{0k}^\alpha
|^2+2|f_{0\overline{k}}^\alpha |^2-
\widetilde{Q}(df_{L,\widetilde{H}}(\xi)\wedge
df_{H,\widetilde{H}}(\eta _k),
\overline{df_{L,\widetilde{H}}(\xi)\wedge df_{H,\widetilde{H}}(\eta _k)}) \\
&-\widetilde{Q}(df_{L,\widetilde{H}}(\xi)\wedge
df_{H,\widetilde{H}}(\eta _{
\overline{k}}),\overline{df_{L,\widetilde{H}}(\xi)\wedge
df_{H,\widetilde{H} }(\eta _{\overline{k}})})\endaligned \tag{8.5}
$$
Note that
$$
\aligned \widetilde{Q}&(df_{L,\widetilde{H}}(\xi)\wedge
df_{H,\widetilde{H}}(\eta _k),
\overline{df_{L,\widetilde{H}}(\xi)\wedge
df_{H,\widetilde{H}}(\eta
_k)})\\
+ &\widetilde{Q}(df_{L,\widetilde{H}}(\xi)\wedge
df_{H,\widetilde{H}}(\eta _{
\overline{k}}),\overline{df_{L,\widetilde{H}}(\xi)\wedge
df_{H,\widetilde{H}
}(\eta _{\overline{k}})}) \\
=&\frac 12\widetilde{Q}(df_{L,\widetilde{H}}(\xi)\wedge
df_{H,\widetilde{H}
}(e_k-iJe_k),df_{L,\widetilde{H}}(\xi)\wedge df_{H,\widetilde{H}}(e_k+iJe_k)) \\
+&\frac 12\widetilde{Q}(df_{L,\widetilde{H}}(\xi)\wedge
df_{H,\widetilde{H}
}(e_k+iJe_k),df_{L,\widetilde{H}}(\xi)\wedge df_{H,\widetilde{H}}(e_k-iJe_k)) \\
=&\widetilde{Q}(df_{L,\widetilde{H}}(\xi)\wedge
df_{H,\widetilde{H}}(e_k),df_{L,\widetilde{H}}(\xi)\wedge
df_{H,\widetilde{H}}(e_k))\\+&\widetilde{Q}(df_{L,\widetilde{H}}(\xi)\wedge
df_{H,\widetilde{H}}(Je_k),df_{L,\widetilde{H}}(\xi)\wedge
df_{H,\widetilde{H}}(Je_k))\\
\leq &0
\endaligned \tag{8.6}
$$
in view of the curvature condition on $N$. Consequently
$$
\delta\widetilde{\rho _1}\geq 2|f_{0k}^\alpha
|^2+2|f_{0\overline{k}}^\alpha |^2
$$
which gives that $f_{0k}^\alpha =f_{0\overline{k} }^\alpha =0$ by
the divergence theorem. In terms of Lemma 8.1, we see that $f$ is
foliated. \qed

\enddemo

In terms of Remark 3.2 and Theorem 8.2, we find that
$(H,\widetilde{H})$-harmonic maps seem to be more sensitive to the
foliated structures than to the horizontal distributions, at least
when the target manifolds are Sasakian manifolds with non-positive
horizontal sectional curvature. However, we will see that under
some further conditions on the maps and the target manifolds,
these critical maps may be related to the CR structures in a
horizontally projective way.

\proclaim{Theorem 8.3} Let $(M^{2m+1},H(M),J,\theta )$ and
$(N^{2n+1},\widetilde{H}(N),\widetilde{J} ,\widetilde{\theta })$
be compact Sasakian manifolds and $f:M\rightarrow N$ be a
$(H,\widetilde{H})$-harmonic map. If $(N^{2n+1},\widetilde{\nabla
})$ has strongly seminegative horizontal curvature, then $f$ is
$(H,\widetilde{H})$-pluriharmonic and
$$
\widetilde{Q}(df_{H,\widetilde{H}}(\eta _j)\wedge
df_{H,\widetilde{H}}(\eta _k),\overline{df_{H,\widetilde{H}}(\eta
_j)\wedge df_{H,\widetilde{H}}(\eta _k)})=0\tag{8.7}
$$
for any unitary frame $\{\eta _j\}$ of $H^{1,0}(M)$.
\endproclaim

\demo{Proof} Let us define a global $1$-form by
$$
\widetilde{\rho _2}=-(f_{\overline{j}k}^\alpha
f_j^{\overline{\alpha }}\theta ^k+f_{j
\overline{k}}^{\overline{\alpha }}f_{\overline{j}}^\alpha \theta
^{\overline{ k}}).  \tag{8.8}
$$
According to (1.20), (1.22), (2.17) and (2.38), we compute that
$$
\aligned \delta\widetilde{\rho_2}=&(f_{\overline{j}k}^\alpha
f_j^{\overline{\alpha
}})_{\overline{k}}+(f_{j\overline{k}}^{\overline{\alpha
}}f_{\overline{j}}^\alpha
)_k \\
=&f_{\overline{j}k}^\alpha f_{j\overline{k}}^{\overline{\alpha
}}+f_{j \overline{k}}^{\overline{\alpha }}f_{\overline{j}k}^\alpha
+f_{\overline{j}k \overline{k}}^\alpha f_j^{\overline{\alpha
}}+f_{j\overline{k}k}^{\overline{
\alpha }}f_{\overline{j}}^\alpha \\
=&2|f_{\overline{j}k}^\alpha |^2+(f_{k\overline{j}}^\alpha
-if_0^\alpha \delta _k^j)_{\overline{k}}f_j^{\overline{\alpha
}}+(f_{\overline{k}j}^{ \overline{\alpha }}+if_0^{\overline{\alpha
}}\delta _{\overline{k}}^{
\overline{j}})_kf_{\overline{j}}^\alpha \\
=&2|f_{\overline{j}k}^\alpha
|^2+(f_{k\overline{j}\overline{k}}^\alpha -if_{0
\overline{k}}^\alpha \delta _k^j)f_j^{\overline{\alpha
}}+(f_{\overline{k} jk}^{\overline{\alpha
}}+if_{0k}^{\overline{\alpha }}\delta _{\overline{k}}^{
\overline{j}})f_{\overline{j}}^\alpha \\
=&2|f_{\overline{j}k}^\alpha |^2+i(f_{0j}^{\overline{\alpha
}}f_{\overline{j} }^\alpha -f_{0\overline{j}}^\alpha
f_j^{\overline{\alpha }})+f_{k\overline{j} \overline{k}}^\alpha
f_j^{\overline{\alpha }}+f_{\overline{k}jk}^{\overline{
\alpha }}f_{\overline{j}}^\alpha \\
=&2|f_{\overline{j}k}^\alpha |^2+i(f_{0j}^{\overline{\alpha
}}f_{\overline{j} }^\alpha -f_{0\overline{j}}^\alpha
f_j^{\overline{\alpha }})+f_{k\overline{k} \overline{j}}^\alpha
f_j^{\overline{\alpha }}+f_{\overline{k}kj}^{\overline{
\alpha }}f_{\overline{j}}^\alpha \\
&+f_j^{\overline{\alpha }}f_k^\beta \widehat{R}_{\beta \gamma
\overline{\delta }}^\alpha (f_{\overline{k}}^\gamma
f_{\overline{j}}^{\overline{\delta }}-f_{\overline{j}}^\gamma
f_{\overline{k}}^{\overline{\delta }})+f_{ \overline{j}}^\alpha
f_{\overline{k}}^{\overline{\beta }}\widehat{R}_{ \overline{\beta
}\overline{\gamma }\delta }^{\overline{\alpha }}(f_k^{
\overline{\gamma }}f_j^\delta -f_j^{\overline{\gamma }}f_k^\delta ) \\
=&2|f_{\overline{j}k}^\alpha |^2+i(f_{0j}^{\overline{\alpha
}}f_{\overline{j} }^\alpha -f_{0\overline{j}}^\alpha
f_j^{\overline{\alpha }})+f_{k\overline{k} \overline{j}}^\alpha
f_j^{\overline{\alpha }}+f_{\overline{k}kj}^{\overline{
\alpha }}f_{\overline{j}}^\alpha \\
&-\widetilde{Q}(df_{H,\widetilde{H}}(\eta _j)\wedge
df_{H,\widetilde{H}}(\eta _k),\overline{df_{H,\widetilde{H}}(\eta
_j)\wedge df_{H,\widetilde{H}}(\eta _k)}).
\endaligned \tag{8.9}
$$
Note that a Sasakian manifold with strongly semi-negative
horizontal curvature has automatically non-positive horizontal
sectional curvature. Then Theorem 8.2 yields that $f_0^\alpha
=f_0^{\overline{\alpha } }=0$, and thus
$$
f_{0j}^\alpha =f_{0\overline{j}}^\alpha =0. \tag{8.10}
$$
From the fourth equation of (2.17), we derive that
$$
\aligned 2f_{k\overline{k}}^\alpha &=f_{k\overline{k}}^\alpha
+f_{\overline{k}
k}^\alpha +mif_0^\alpha \\
&=mif_0^\alpha \\
&=0.
\endaligned\tag{8.11}
$$
It follows from (8.9), (8.10) and 8.11) that
$$
\delta\widetilde{\rho_2}= 2|f_{\overline{j}k}^\alpha
|^2-\widetilde{Q}(df_{H, \widetilde{H}}(\eta _j)\wedge
df_{H,\widetilde{H}}(\eta _k),\overline{df_{H, \widetilde{H}}(\eta
_j)\wedge df_{H,\widetilde{H}}(\eta _k)}). \tag{8.12}
$$
In terms of the curvature condition of $(N,\widetilde{\nabla })$
and the divergence theorem, we get from (8.12) that
$$
f_{\overline{j}k}^\alpha =0
$$
and
$$
\widetilde{Q}(df_{H,\widetilde{H}}(\eta _j)\wedge
df_{H,\widetilde{H}}(\eta _k),\overline{df_{H,\widetilde{H}}(\eta
_j)\wedge df_{H,\widetilde{H}}(\eta _k)})=0
$$
for any $1\leq j,k\leq m$. Since $f_0^\alpha =0$, (2.17) implies
that $ f_{k\overline{j}}^\alpha=f_{\overline{j}k}^\alpha=0$.
According to (5.2) and (5.3), $f$ is
$(H,\widetilde{H})$-pluriharmonic. Hence we complete the proof.
\qed
\enddemo

\proclaim{Theorem 8.4} Suppose $(M^{2m+1},H(M),J,\theta )(m\ge 2)$
and $(N^{2n+1},\widetilde{H}(N),\widetilde{J},\widetilde{\theta
})$ are compact Sasakian manifolds and $N$ has strongly negative
horizontal curvature. Suppose $f:M\rightarrow N$ is a
$(H,\widetilde{H})$-harmonic map with
$\max_Mrank_R\{df_{H,\widetilde{H}}\}\geq 3$. Then $f$ is either a
foliated $(H,\widetilde{H})$-holomorphic map or a foliated
$(H,\widetilde{H})$-antiholomorphic map.
\endproclaim

\demo{Proof} From Theorem 8.2, we know that $f$ is foliated. The
rank condition for $f$ means that there exists a point $p\in M$
such that $rank_R\{df_{H,\widetilde{H}}(p)\}\geq 3$. Consequently
$ rank_R\{df_{H,\widetilde{H}}\}\geq 3$ in some connected open
neighborhood $U$ of $p$. Write $df_{H,\widetilde{H}}(\eta
_{\overline{j}})=f_{\overline{j} }^\alpha \widetilde{\eta }_\alpha
+f_{\overline{j}}^{\overline{\alpha }} \widetilde{\eta
}_{\overline{\alpha }}$. Then
$$
\left( df_{H,\widetilde{H}}(\eta _{\overline{j}})\wedge
df_{H,\widetilde{H} }(\eta _{\overline{k}})\right)
^{(1,1)}=(f_{\overline{j}}^\alpha f_{
\overline{k}}^{\overline{\beta }}-f_{\overline{k}}^\alpha
f_{\overline{j}}^{ \overline{\beta }})\widetilde{\eta }_\alpha
\wedge \widetilde{\eta }_{ \overline{\beta }}.  \tag{8.13}
$$
Since $N$ has strongly negative horizontal curvature, it follows
from (8.7) and (8.13) that
$$
f_{\overline{j}}^\alpha f_{\overline{k}}^{\overline{\beta
}}-f_{\overline{k} }^\alpha f_{\overline{j}}^{\overline{\beta }}=0
\tag{8.14}
$$
for $1\leq j, k\leq m$ and $1\leq \alpha , \beta \leq n$. We want
to show that for every point $q\in U$, either $\partial
f_{H,\widetilde{H}}(q)=0$ or $ \overline{\partial
}f_{H,\widetilde{H}}(q)=0$.

Without loss of generality, we assume that $\partial
f_{H,\widetilde{H}}(q)\neq 0$. This means that $ f_k^\gamma (q)
\neq 0$ for some $1\leq k\leq m$ and some $1\leq \gamma \leq n$.
Therefore (8.14) yields that
$$
\{f_{\overline{1}}^\alpha(q) ,\cdots ,f_{\overline{m}}^\alpha(q)
\}=c^\alpha \{f_{ \overline{1}}^{\overline{\gamma}}(q),\cdots
,f_{\overline{m}}^{\overline{ \gamma}}(q)\}  \tag{8.15}
$$
for each $1\leq \alpha \leq m$, where $c^\alpha
=f_{\overline{k}}^\alpha(q) /f_{
\overline{k}}^{\overline{\gamma}}(q)$. If $\overline{\partial
}f_{H,\widetilde{ H}}(q)\neq 0$ too, then $f_{\overline{j}}^\delta
(q)\neq 0$ for some $1\leq j\leq m\ $and some $1\leq \delta \leq
n$. Using (8.14) again, we get
$$
\{f_{\overline{1}}^{\overline{\beta }}(q),\cdots
,f_{\overline{m}}^{\overline{ \beta }}(q)\}=d^{\overline{\beta
}}\{f_{\overline{1}}^\delta (q),\cdots ,f_{ \overline{m}}^\delta
(q)\} \tag{8.16}
$$
for each $1\leq \beta \leq m$, where $d^{\overline{\beta
}}=f_{\overline{j} }^{\overline{\beta
}}(q)/f_{\overline{j}}^\delta (q)$. In terms of (8.15) and (8.16),
we find
$$
\aligned df_{H,\widetilde{H}}(\eta
_{\overline{i}})&=f_{\overline{i}}^\alpha \eta _\alpha
+f_{\overline{i
}}^{\overline{\alpha }}\eta _{\overline{\alpha }} \\
&=f_{\overline{i}}^{\overline{\gamma }}(c^\alpha \eta _\alpha
)+f_{\overline{i
}}^\delta (d^{\overline{\alpha }}\eta _{\overline{\alpha }}) \\
&=f_{\overline{i}}^{\overline{\gamma }}(c^\alpha \eta _\alpha
)+c^\delta f_{ \overline{i}}^{\overline{\gamma
}}(d^{\overline{\alpha }}\eta _{\overline{
\alpha }}) \\
&=f_{\overline{i}}^{\overline{\gamma }}(c^\alpha \eta _\alpha
+c^\delta d^{ \overline{\alpha }}\eta _{\overline{\alpha }})
\endaligned\tag{8.17}
$$
at $q\in U$. From (8.17), we have
$$
\aligned span_C\{df_{H,\widetilde{H}}(\eta _{\overline{i}})\}_q
&=span_C\{c^\alpha \eta _\alpha +c^\delta d^{\overline{\alpha
}}\eta _{\overline{
\alpha }}\} \\
&\leq 1
\endaligned
$$
which implies that $rank_R\{df_{H,\widetilde{H}}(q)\}\leq 2$,
contradicting $q\in U$. Thus, for each $q\in U$, either $\partial
f_{H,\widetilde{H}}(q)=0$ or $\overline{
\partial }f_{H,\widetilde{H}}(q)=0$.

Since $rank_R\{df_{H,\widetilde{H}}>0\}$ in $U$, the two sets
$U\cap \{\partial f_{H, \widetilde{H}}=0\}$ and $U\cap
\{\overline{\partial }f_{H,\widetilde{H}}=0\}$ are disjoint closed
subsets of $U$, and their union is the connected set $U$ . It
follows that either $\partial f_{H,\widetilde{H}}\equiv 0$ on $U$
or $ \overline{\partial }f_{H,\widetilde{H}}\equiv 0$ on $U$. From
Theorem 5.4 and Theorem 8.2, we conclude that either $\partial
f_{H,\widetilde{H}}\equiv 0$ on $M$ or $\overline{\partial
}f_{H,\widetilde{H}}\equiv 0$ on $M$.\qed
\enddemo

\proclaim{Theorem 8.5} Let $k\geq 2$. Suppose
$(M^{2m+1},H(M),J,\theta )$ and
$(N^{2n+1},\widetilde{H}(N),\widetilde{J},\widetilde{\theta })$
are compact Sasakian manifolds and the horizontal curvature of $N$
is negative of order $k$. Suppose $f:M\rightarrow N$ is a
$(H,\widetilde{H})$-harmonic map and
$\max_M\{df_{H,\widetilde{H}}\}\geq 2k$. Then $f$ is either a
foliated $(H,\widetilde{H})$-holomorphic map or a foliated
$(H,\widetilde{H})$-antiholomorphic map.
\endproclaim
\demo{Proof} By a similar argument as that in Theorem 8.4, we may
deduce the conclusion of this theorem from (8.7), (8.13) and
Definition 8.2.\qed
\enddemo

For a manifold $M$, we use $H_{l}(M,R)$ to denote its usual
singular homology.

\proclaim{Corollary 8.6} Suppose $f:(M^{2m+1},H(M),J,\theta
)\rightarrow
(N^{2n+1},\widetilde{H}(N),\widetilde{J},\widetilde{\theta} )$ is
a foliated $(H,\widetilde{H})$-harmonic map between compact
Sasakian manifolds. Suppose the horizontal curvature of $N$ is
negative of order $k$ and $k\ge 2$. If the induced map
$f_{*}:H_l(M,R)\rightarrow H_l(N,R)$ is nonzero for some $l\geq
2k+1$, then $f$ is either $(H,\widetilde{H})$-holomorphic or
$(H,\widetilde{H})$-anti-holomorphic.
\endproclaim
\demo{Proof} The assumption that $f_{*}:H_l(M,R)\rightarrow
H_l(N,R)$ is nonzero for some $l\geq 2k+1$ implies that
$rank_R\{df\}\geq 2k+1$ at some point $p$ of $M$. Since $f$ is
foliated,
$$
rank_R\{\pi _{\widetilde{L}}\circ
df\}+rank_R\{df_{H,\widetilde{H}}\}\ge rank_R\{df\}
$$
at each point of $M$. Hence $rank_R\{df_{H,\widetilde{H}}\}_p\geq
2k$. By Theorem 8.5, we find that $f$ is either a foliated
$(H,\widetilde{H})$-holomorphic or a foliated
$(H,\widetilde{H})$-anti-holomorphic map. \qed
\enddemo

Before investigating the strong rigidity of Sasakian manifolds
with some kind of negative horizontal curvature, let us recall
some basic notions and results for foliations, especially for
Riemannian foliations. Given a foliation $F$ on a manifold $M$, a
differential form $\omega \in \Omega ^{*}(M)$ is called basic if
for every vector field $X$ tangent to the leaves, $i_X\omega =0$
and $i_Xd\omega =0$, where $i_X$ denotes the interior product with
respect to $X$. In particular, a basic form of degree $0$ is just
a basic function. Clearly the exterior derivative of a basic form
is again basic, so all basic forms $(\Omega _B^{*}(M),d_B)$
constitutes a subcomplex of the de Rham complex $(\Omega
^{*}(M),d)$, where $d_B=d|_{\Omega^{*} _B(M)}$. The cohomology
$H_B^{*}(M/F)=ker d_B/Im d_B$ of this subcomplex is called the
basic cohomology of $(M,F)$. In general, the basic cohomology
groups are not always finite dimensional. It is known that
Riemannian foliations on compact manifolds form a large class of
foliations for which the basic cohomology groups are
finite-dimensional (cf. [KHS], [KT], [PR]).

Let $(M,F)$ and $(N,\widetilde{F})$ be two manifolds endowed with
complete Riemannian foliations $F$ and $\widetilde{F}$
respectively. Similar to Definition 3.3, we have the notion of
continuous foliated map between $(M,F)$ and $(N,\widetilde{F})$,
that is, a continuous map from $M$ to $N$ mapping leaves of $F$
into leaves of $\widetilde{F}$. A homotopy between $M$ and $N$
consisting of continuous (resp. smooth) foliated maps is called a
continuous (resp. smooth) foliated homotopy, and the corresponding
homotopy equivalence can be defined in a natural way. Actually we
may work in smooth category due to the following result.

\proclaim{Lemma 8.7} ([LM]) Any continuous foliated map between
complete Riemannian foliations is foliatedly homotopic to a
$C^\infty $ foliated map.
\endproclaim

In view of Lemma 8.7, two complete Riemannian foliated manifolds
$(M,F)$ and $(N,\widetilde{F})$ have the same foliated homotopy
type in the $C^\infty $ sense if and only if they have the same
foliated homotopy type in the usual continuous sense. Another
important property of basic cohomologies of Riemannian foliations
is the homotopy invariance.

\proclaim{Lemma 8.8} (cf. [Vi]) Let $(M,F)$ and
$(N,\widetilde{F})$ be two compact Riemannian foliated manifolds.
Suppose $f,g:(M,F)\rightarrow (N,\widetilde{F})$ are foliated
homotopic. Then $f^{*}=g^{*}:H_B^{*}(N/\widetilde{F})\rightarrow
H_B^{*}(M/F)$. In particular, if $f:(M,F)\rightarrow
(N,\widetilde{F})$ is a foliated homotopy equivalence, then
$f^{*}:H_B^{*}(N/\widetilde{F})\rightarrow H_B^{*}(M/F)$ is an
isomorphism.
\endproclaim

As mentioned in \S 1, the Reeb foliation of a Sasakian manifold is
a Riemannian foliation of dimension 1. Suppose
$(M^{2m+1},H(M),J,\theta )$ is a compact Sasakian manifold with
the Reeb foliation $F_\xi$. Then each basic cohomology group
$H_{B}^{k}(M/F_{\xi})$ ($k=0,1,...,2m$) is finite dimensional.
Using the second property in (1.4), it is easy to verify that
$d\theta$ is a basic form, and thus so is $(d\theta )^k$ for
$2\leq k\leq m$.

\proclaim{Lemma 8.9} Suppose $(M^{2m+1},H(M),J,\theta )$ is a
compact Sasakian manifold. Then $0\neq [(d\theta )^k]_B\in
H_B^{2k}(M/F)$ ($1\leq k\leq m$).
\endproclaim
\demo{Proof} We prove this lemma by contradiction. Suppose
$(d\theta )^k=d\alpha $ for some $\alpha \in \Omega _B^{2k-1}(M)$.
Then
$$
(d\theta)^m=d\alpha \wedge (d\theta )^{m-k} =d\big(\alpha \wedge
(d\theta)^{m-k}\big).\tag{8.18}
$$
Using (8.18) and the Stokes formula, the volume of $M$ is given,
up to a positive constant, by
$$
\aligned
\int_M\theta \wedge (d\theta )^m&=\int_M\theta\wedge d\big(\alpha\wedge(d\theta)^{m-k}\big)\\
&=\int_Md\theta \wedge \alpha \wedge (d\theta )^{m-k}\\
&=0, \endaligned
$$
where the last equality follows from the fact that
$d\theta\wedge\alpha\wedge (d\theta)^{m-k}=0$, since it is a basic
form of degree $2m+1$. Hence we get a contradiction. \qed
\enddemo
\remark{Remark 8.1} The above result is actually true for a
general compact pseudo-Hermitian manifold.
\endremark

Recall that a smooth complex-valued function $f:M\rightarrow C$ on
a CR manifold $M$ is called a CR function if $Z(f)=0$ for any
$Z\in H^{0,1}(M)$.

\definition {Definition 8.2}Let $M^{2m+1}$ be a Sasakian manifold.
We say that a subset $V$ of $M$ is a foliated analytic subvariety
if, for any point $p\in V$, there exists a foliated coordinate
chart $(U,\Phi ;\varphi )$ of $p$ such that $V\cap U$ is the
common zero locus of a finite collection of basic CR functions
$f_1,...,f_k$ on $U$. In particular, $V$ is called a foliated
analytic hypersurface if $V$ is locally the zero locus of a single
nonzero basic CR function $f$. \enddefinition

More explicitly, let $\varphi :U\rightarrow W\subset C^m$ be the
submersion associated with the foliated coordinate chart $(U,\Phi
;\varphi )$, that is, $\varphi =\pi \circ \Phi $ (see (1.23) in \S
1). Since the CR functions $f_1,..,f_k$ are constant along the
leaves, there are holomorphic function
$\widetilde{f}_1,...,\widetilde{f}_k$ on $W$ such that
$f_i=\widetilde{f}_i\circ \pi $ ($i=1,...,k$). Set
$\widetilde{V}_{\varphi}=\varphi (V\cap U)$. Thus
$\widetilde{V}_{\varphi}$ is a complex analytic subvariety in $W$
defined by the common zero locus of
$\widetilde{f}_1,..,\widetilde{f}_k$.

A point $p\in V$ is called a smooth point of $V$ if $V$ is a
submanifold of $M$ near $p$. The locus of smooth points of $V$ is
denoted by $V^{*}$. A point $p\in V-V^{*}$ is called a singular
point of $V$; the singular locus $V-V^{*}$ of $V$ is denoted by
$V^s$. Similarly we have the notions of smooth points and singular
points for $\widetilde{V}_\varphi$. Let
$\widetilde{V}^{*}_\varphi$ (resp. $\widetilde{V}^s_\varphi$)
denote the locus of smooth points (resp. singular points) of
$\widetilde{V}_\varphi$. Clearly
$\widetilde{V}^{*}_\varphi=\varphi(V^{*}\cap U)$ and
$\widetilde{V}_\varphi ^s=\varphi (V^s\cap U)$ of $\widetilde{V}$.
By the proposition on page 32 in [GH], we know that
$\widetilde{V}^{*}_\varphi$ has finite volume in bounded regions.
Consequently $V^{*}$ has finite volume in bounded regions too.
Therefore we may define the integral of a differential form
$\omega $ on $M$ over $V$ to be the integral of $\omega $ over the
smooth locus $V^{*}$ of $V$.

We need the following Stokes' formula for foliated analytic
subvarieties, which is a generalization of the usual Stokes'
formula for analytic subvarieties.

\proclaim{Proposition 8.10} Let $M$ be a Sasakian manifold and let
$V\subset M$ be a foliated analytic subvarieties of real dimension
$2k+1$. Suppose $\alpha $ is a differential form of degree $2k$
with compact support in $M$. Then
$$
\int_Vd\alpha =0.
$$
\endproclaim
\demo{Proof} The question is local, it will be sufficient to show
that for every point $p\in V$, there exists a neighborhood $U$ of
$p$ such that for every $\alpha \in A_c^{2k}(U)$ (the space of
differential forms of degree $2k$ with compact support in $U$)
$$\int_Vd\alpha =0.$$
Suppose $(U,\Phi ;\varphi )$ is a bounded foliated coordinate
chart around $p$. Let $\varphi :U\rightarrow W$ be the submersion
associated with $(U,\Phi ;\varphi )$ and let
$\widetilde{p}=\varphi (p)\in \widetilde{V}_\varphi \subset W$. By
the local structure of an analytic subvariety in $C^m$ (cf. [GH]),
we may find a coordinate system $z=(z_1,...,z_m)$ and a
polycylinder $\bigtriangleup $ around $\widetilde{p}$ such that
the projection $\widetilde{\pi }:(z_1,...,z_m)\rightarrow
(z_1,...,z_k)$ expresses $\widetilde{V}_\varphi \cap
\bigtriangleup $ as a branched cover of $\bigtriangleup ^{\prime
}=\widetilde{\pi }(\bigtriangleup )$, branched over an analytic
hypersurface $\Sigma \subset \bigtriangleup ^{\prime }$. Let
$T^\varepsilon $ be the $\varepsilon $-neighborhood of $\Sigma $
in $\bigtriangleup ^{\prime }$ and
$$
\widetilde{V}_\varphi ^\varepsilon =(\widetilde{V}_\varphi \cap
\bigtriangleup )-\widetilde{\pi }^{-1}(T^\varepsilon )\subset
\widetilde{V}_\varphi ^{*}.
$$
Set $V^\varepsilon =\varphi ^{-1}(\widetilde{V}_\varphi
^\varepsilon )$. Clearly $\varphi ^{-1}(\bigtriangleup )\subset U$
is a foliated neighborhood of $p$. For $\alpha \in
A_c^{2k}(\varphi ^{-1}(\bigtriangleup ))$, we have
$$
\aligned \int_Vd\alpha &=\int_{V\cap \varphi ^{-1}(\bigtriangleup
)}d\alpha \\
&=\lim_{\varepsilon \rightarrow 0}\int_{V^\varepsilon }d\alpha \\
&=\lim_{\varepsilon \rightarrow 0}\int_{\partial V^\varepsilon}\alpha\\
&=\lim_{\varepsilon \rightarrow 0}\int_{\varphi ^{-1}(\partial
\widetilde{\pi }^{-1}(T^\varepsilon ))}\alpha.
\endaligned
$$
Thus to prove the result, we simply have to prove that
$vol(\varphi ^{-1}(\partial \widetilde{\pi }^{-1}(T^\varepsilon
)))\rightarrow 0$ as $\varepsilon \rightarrow 0$ or equivalently
$vol(\partial \widetilde{\pi }^{-1}(T^\varepsilon ))\rightarrow 0$
as $\varepsilon \rightarrow 0$. However, the latter one has
already be shown on page 33 in [GH] .\qed
\enddemo

We will use special $(H,\widetilde{H})$-harmonic maps to establish
strong rigidity results for Sasakian manifolds.

\proclaim{Lemma 8.11} Suppose $f:(M^{2m+1},H(M),J,\theta
)\rightarrow
(N^{2n+1},\widetilde{H}(N),\widetilde{J},\widetilde{\theta })$ is
a foliated special $(H,\widetilde{H})$-harmonic map between
Sasakian manifolds. If $M$ is compact, then
$df(\xi)=\lambda\widetilde{\xi}$ with $\lambda$ constant.
\endproclaim
\demo{Proof} Since f is foliated, $df(\xi)=\lambda\widetilde{\xi}$
for some smooth function $\lambda$ on M, that is, $f_0^0=\lambda$
and $f_0^\alpha =f_0^{\overline{\alpha }}=0$. Then (2.14) yields
that
$$
f_{0j}^0=f_{j0}^0,\quad
f_{0\overline{j}}^0=f_{\overline{j}0}^0.\tag{8.19}
$$
In terms of (2.24) and (8.19), we get
$$
f_{0j\overline{j}}^0=f_{j0\overline{j}}^0=f_{j\overline{j}0}^0,\quad
f_{0\overline{j}j}^0=f_{\overline{j}0j}^0=f_{\overline{j}j0}^0.\tag{8.20}
$$
Since f is special, we derive from (8.19) and (8.20) that
$$
\aligned \bigtriangleup
_b\lambda=&f_{0j\overline{j}}^0+f_{0\overline{j}j}^0=f_{j\overline{j}0}^0+f_{\overline{j}j0}^0\\
=&0.\endaligned
$$
Due to the compactness of $M$, it follows that $\lambda$ is
constant.\qed
\enddemo

\proclaim{Theorem 8.12}Let $f:(M,H(M),J,\theta )\rightarrow
(N,\widetilde{H}(N),\widetilde{J},\widetilde{\theta })$ be a
foliated special $(H,\widetilde{H})$-harmonic map of Sasakian
manifolds, both of $CR$ dimension $m\geq 2$. Suppose the
horizontal curvature of $N$ is either strongly negative or
adequately negative. Suppose $f$ is of degree $1$ and the map
$f^{*}:H_B^{2m-2}(N)\rightarrow H_B^{2m-2}(M)$ induced by $f$ is
surjective. Then $f:M\rightarrow N$ is either
$(H,\widetilde{H})$-biholomorphic or
$(H,\widetilde{H})$-anti-biholomorphic.
\endproclaim
\demo{Proof} Let $V$ be the set of points of $M$ where $f$ is not
locally diffeomorphic. Suppose $V$ is nonempty. We will derive a
contradiction in the following steps.

{\it Step 1}. Since $\deg (f)=1$, $V\neq M$ and $f$ maps
$M-f^{-1}(f(V))$ bijectively onto $N-f(V)$. By either Theorem 8.4
or Theorem 8.5,  $f$ is either $(H,\widetilde{H})$-holomorphic or
$(H,\widetilde{H})$-antiholomorphic. Without loss of generality,
we assume now that $f$ is a foliated special
$(H,\widetilde{H})$-holomorphic map. In terms of Lemma 8.11 and
$V\neq M$, we see that $df(\xi )=\lambda \widetilde{\xi }$ for
some nonzero constant $\lambda $. For any point $p\in M$, let
$q=f(p)\in N$. Let $ (U_{1},\Phi _{1};\varphi _{1})$ and
$(U_{2},\Phi _{2};\varphi _{2})$ be foliated coordinate charts
around $p$ and $q$ respectively, and let $\varphi
_{i}:U_{i}\rightarrow W_{i}\subset C^{m}$ be the submersion
associated with $ (U_{i},\Phi _{i};\varphi _{i})$ ($i=1,2$).
Suppose $f(U_{1})\subset U_{2}$. Then $f$ induces a holomorphic
map $\widetilde{f}:W_{1}\rightarrow W_{2}$ such that $\varphi
_{2}\circ f=\widetilde{f}\circ \varphi _{1}$. Clearly $ p\in V$ if
and only if $\widetilde{f}$ is not diffeomorphic at $\varphi
_{1}(p)$. Hence $V$ is defined locally by the zero locus of $\det
(\partial w^{\alpha }\circ \widetilde{f}/\partial z^{i})\circ
\varphi _{1}$, where $ (z^{i})$ and $(w^{\alpha })$ are
holomorphic coordinate systems of $W_{1}$ and $W_{2}$
respectively. This shows that $V$ is a foliated analytic
hypersurface in $M$. Since $f$ is a foliated
$(H,\widetilde{H})$-holomorphic map, $ f(V)$ is a foliated
analytic subvariety in $N$. It is obvious that both $
M-f^{-1}(f(V))$ and $N-f(V)$ are foliated open submanifolds of $M$
and $N$ respectively.

{\it Step 2}. We claim that $f$ is a horizontally one-to-one map
from $M-f^{-1}(f(V))$ to $N-f(V)$. Let $\widehat{
f}:(M-f^{-1}(f(V)))/F_{\xi }\rightarrow
(N-f(V))/\widetilde{F}_{\widetilde{ \xi }}$ denote the induced map
of $f|_{M-f^{-1}(f(V))}$. Suppose there are two leaves $L_{1}$,
$L_{2}\subset M-f^{-1}(f(V))$ and a leaf $\widetilde{L} \subset
N-f(V)$ such that $f(L_{1}),f(L_{2})\subset \widetilde{L}$. Let
$\gamma _{1}(t), \gamma_{2}(t)$ ($t\in (-\infty ,\infty )$) be the
integral curves of $\xi $ whose images are $ L_{1}$ and $L_2$
respectively. In terms of the fact that $df(\xi )=\lambda
\widetilde{\xi }$ with constant $\lambda \neq 0$, we see that both
$f(\gamma _{1})$ and $f(\gamma_{2})$ are integral curves of
$\lambda \widetilde{\xi }$ with possibly different initial points.
As a result, their image must be $\widetilde{L}$, that is,
$f(L_{1})=f(L_{2})=\widetilde{L}$. The injectivity of
$f|_{M-f^{-1}(f(V))}$ implies that $L_{1}=L_{2}$, that is,
$\widehat{f}$ is injective. By the surjectivity of $f$ from
$M-f^{-1}(f(V))$ to $N-f(V)$, we conclude that $ \widehat{f}$ is
surjective, and thus $\widehat{f}$ is bijective. By Proposition
5.7, $f:$ $M-f^{-1}(f(V))\rightarrow N-f(V)$ is a foliated $(H,
\widetilde{H})$-biholomorphism.

{\it Step 3}. Now we assert that $f(V)$ must be a foliated
analytic subvariety with real codimension at least $4$, that is,
the transversal complex codimension is at least two. Otherwise,
suppose $f(V)$ is also a foliated analytic hypersurface, we are
going to prove that the critical points set $V$ of $f$ is actually
removable. For any $p\in V$, let $ (U_{1},\Phi _{1};\varphi _{1})$
and $(U_{2},\Phi _{2};\varphi _{2})$ be foliated coordinate charts
around, respectively, $p$ and $q=f(p)$ as in Step 1, such that
$f(U_{1})\subset U_{2}$ and $f(V)\cap U_{2}$ is defined by the
zero locus of a single basic CR function. Set $V_{\varphi
_{1}}=\varphi _{1}(V\cap U_{1})$ and $[f(V)]_{\varphi
_{2}}=\varphi _{2}\left( f(V)\cap U_{2}\right) $. From Steps 1 and
2, we know that the induced holomorphic map
$\widetilde{f}:W_1\rightarrow W_2$ is injective on
$\widetilde{f}:W_1-\widetilde{f}^{-1}([f(V)]_{\varphi_2})$.
Clearly $V_{\varphi _{1}}$ and $[f(V)]_{\varphi _{2}}$ are
analytic hypersurfaces in $W_{1}$ and $W_{2}$ respectively, and
$\widetilde{f }(V_{\varphi _{1}})\subset \lbrack f(V)]_{\varphi
_{2}}$. Then there exists $ v\in V_{\varphi _{1}}$ such that $v$
is an isolated point of $\widetilde{f} ^{-1}(\widetilde{f}(v))$
and, by using a local coordinate chart $(z^{i})$ of $W_{1}$ at $v$
and by applying the Riemann removable singularity theorem to $
z^{i}\circ \widetilde{f}^{-1}$ on $W_{2}-[f(V)]_{\varphi _{2}}$
for some open neighborhood of $\widetilde{f}(v)$ in $W_{2}$, we
find that $\widetilde{ f}$ is locally diffeomorphic at $v$. This
implies that $f$ is locally diffeomorphic at any point in $\varphi
_{1}^{-1}(v)$, contradicting $\varphi _{1}^{-1}(v)\subset V$.

{\it Step 4}. From Lemma 8.9, we know that $[(d\theta
)^{m-1}]_{B}\neq 0$. By the assumption that $f^{\ast
}:H_{B}^{2m-2}(N)\rightarrow H_{B}^{2m-2}(M)$ is surjective, there
exists an element $[\beta ]_{B}\in H_{B}^{2m-2}(N)$ such that
$f^{\ast }[\beta ]_{B}=[(d\theta )^{m-1}]_{B}$, that is,
$$
f^{\ast }\beta =(d\theta )^{m-1}+d\alpha\tag{8.21}
$$
for some $\alpha \in \Omega _{B}^{2m-3}(M)$. Since $df(\xi
)=\lambda \widetilde{\xi }$, we may write $f^{\ast
}\widetilde{\theta } =\lambda \theta +f_{j}^{0}\theta
^{j}+f_{\overline{j}}^{0}\theta ^{\overline{ j}}$. Note that
$f_{j}^{0}\theta ^{j}+f_{\overline{j}}^{0}\theta ^{\overline{j}}$
is a global 1-form on $M$. Since $V$ is foliated and $i_{\xi
}\{(f_{j}^{0}\theta ^{j}+f_{\overline{j}}^{0}\theta
^{\overline{j}})\wedge f^{\ast }\beta \}=0$, we have
$$
\int_{V}(f_{j}^{0}\theta ^{j}+f_{\overline{j}}^{0}\theta
^{\overline{j}})\wedge f^{\ast }\beta =0.  \tag{8.22}
$$
From (8.21), (8.22) and Proposition 8.12, we get
$$
\aligned \int_{V}f^{\ast }\widetilde{\theta }\wedge f^{\ast }\beta
&=\lambda \int_{V}\theta \wedge f^{\ast }\beta +\int_{V}\left(
f_{j}^{0}\theta ^{j}+f_{
\overline{j}}\theta ^{\overline{j}}\right) \wedge f^{\ast }\beta \\
&=\lambda \int_{V}\theta \wedge (d\theta )^{m-1}+\int_{V}\theta
\wedge d\alpha
\\
&=\lambda \int_{V}\theta \wedge (d\theta )^{m-1}+\int_{V}d\theta
\wedge \alpha. \endaligned  \tag{8.23}
$$
Clearly we have $i_{\xi }\left( d\theta \wedge \alpha \right) =0$,
which implies that
$$
\int_{V}d\theta \wedge \alpha =0.\tag{8.24}
$$
It follows from (8.23) and (8.24) that
$$
\int_{V}f^{\ast }\widetilde{\theta }\wedge f^{\ast }\beta =\lambda
\int_{V}\theta \wedge (d\theta )^{m-1}>0.  \tag{8.25}
$$
On the other hand, since $f(V)$ is of dimension less than $2m-3$,
one has
$$
\int_{V}f^{\ast }\widetilde{\theta }\wedge f^{\ast }\beta
=\int_{f(V)}\widetilde{\theta }\wedge \beta =0
$$
which contradicts to (8.25).

It follows from the above discussion that $V$ must be empty. In
terms of step 2, we may conclude that $f:M\rightarrow N$ is a
foliated $(H,\widetilde{H})$-biholomorphism. \qed

\enddemo
\remark{Remark 8.2} The argument for removing the critical points
of $\widetilde{f}$ in Step 3 is inspired by the related argument
in Theorem 8 of [Si1].
\endremark

\proclaim{Corollary 8.13} Let
$(N^{2n+1},\widetilde{H}(N),\widetilde{J},\widetilde{\theta })$
($n\geq 2$) be a compact regular Sasakian manifold with either
strongly negative or adequately negative horizontal curvature.
Suppose $(M,H(M),J,\theta )$ is a compact Sasakian manifold with
the same foliated homotopy type as $N$. Then $(M,H(M),J)$ is
$(H,\widetilde{H})$-biholomorphic to either $(N,
\widetilde{H}(N),\widetilde{J})$ or
$(N,\widetilde{H}(N),-\widetilde{J})$.
\endproclaim
\demo{Proof}Since $M$ is foliatedly homotopic to $N$, we have a
foliated smooth homotopy equivalent map $h:M\rightarrow N$.
Consequently $\deg (h)=1$ and $h^{\ast
}:H_{B}^{2m-2}(N)\rightarrow H_{B}^{2m-2}(M)$ is an isomorphism.
By Theorem 7.14, there exists a foliated special
$(H,\widetilde{H})$-harmonic map $f:M\rightarrow N$ which is
foliatedly homotopic to $h$. In terms of Lemma 8.8, $deg(f)=1$ and
$f^{\ast }:H_{B}^{2m-2}(N)\rightarrow H_{B}^{2m-2}(M)$ is an
isomorphism too. This corollary then follows immediately from
Theorem 8.12. \qed
\enddemo
\remark{Remark 8.3} Note that the
$(H,\widetilde{H})$-biholomorphism $f$ between $M$ and $N$ in
Corollary 8.13 is actually a vertically homothetic map, that is,
$df(\xi)=\lambda\widetilde{\xi}$ with $\lambda$ constant. In
Example 1.1, we give some Sasakian manifolds with either strongly
negative or adequately negative horizontal curvature. These
Sasakian manifolds, which may be regarded as model spaces, appear
also in the classification of contact sub-symmetric spaces by
[BFG]. As applications, Corollary 8.13 exhibits the foliated
strong rigidity of these model spaces.
\endremark

\heading{\bf Appendix}
\endheading
\vskip 0.3 true cm

{\bf A. Pseudo-Hermitian harmonic maps}
\vskip 0.2 cm

In this subsection, we introduce another natural generalized
harmonic map between pseudo-Hermitian manifolds.

\definition{Definition A1} A map $f:(M^{2m+1},H(M),J,\theta )\rightarrow
(N^{2n+1},\widetilde{H}(N),\widetilde{J},\widetilde{\theta })$ is
called a pseudo-Hermitian harmonic map if it satisfies
$$
\tau (f)=tr_{g_{\theta }}\beta =0,  \tag{A1}
$$
that is
$$
\aligned &f_{00}^{\alpha }+f_{k\overline{k}}^{\alpha
}+f_{\overline{k}k}^{\alpha }=0\\
&f_{00}^{0}+f_{k\overline{k}}^{0}+f_{\overline{k}k}^{0}=0.\endaligned\tag{A2}
$$
\enddefinition
\remark{Remark A1} Similar ideas for introducing generalized
harmonic maps as Definition A1 were also mentioned in [DT] and
[Kok].
\endremark

Clearly if $f$ is a foliated pseudo-Hermitian harmonic map, then
$f$ automatically satisfies the equation
$f_{k\overline{k}}^{\alpha }+f_{\overline{k}k}^{\alpha }=0$, that
is, $f$ is $(H,\widetilde{H})$-harmonic. Concerning the Question
proposed in section 5, we establish a continuation result about
the foliated property for pseudo-Hermitian harmonic maps.

\proclaim{Theorem A1} Let $M$ and $N$ be Sasakian manifolds and
let $f:M\rightarrow N$ be a pseudo-Hermitian harmonic map. Assume
that $U$ is a nonempty open subset of $M$. If $f$ is foliated on
$U$, then $f$ is foliated on $M$.
\endproclaim
\demo{Proof} The argument is similar to that for Theorem 5.4.
Choose any point $q\in \partial U$. Let $W$ be a connected open
neighborhood of $q$ such that there exist a frame field $\{\xi
,\eta _{1},...,\eta _{m},\eta _{\overline{1} },...,\eta
_{\overline{m}}\}$ on $W$ and a frame field $\{\widetilde{\xi },
\widetilde{\eta }_{1},...,\widetilde{\eta }_{n},\widetilde{\eta
}_{\overline{ 1}},...,\widetilde{\eta }_{\overline{n}}\}$ on some
open neighborhood of $f(W)$ respectively. Write
$$
df_{L,\widetilde{H}^{1,0}}=f_{0}^{\alpha }\theta \otimes
\widetilde{\eta }_{\alpha }.
$$
Clearly $df_{L,\widetilde{H}^{1,0}}\in
Hom(L,f^{-1}\widetilde{H}^{1,0})$. By definition, it is easy to
derive the following
$$
\aligned \bigtriangleup df_{T,\widetilde{H}^{1,0}}&=tr_{g_{\theta
}}D^{2}df_{T,\widetilde{H}^{1,0}} \\
&=(f_{000}^{\alpha }+f_{0k\overline{k}}^{\alpha
}+f_{0\overline{k}k}^{\alpha
})\theta \otimes \widetilde{\eta }_{\alpha } \\
&=(\bigtriangleup _{M}f_{0}^{\alpha })\theta \otimes
\widetilde{\eta } _{\alpha }+tr_{g_{\theta }}\{df_{0}^{\alpha
}\otimes \theta \otimes \widetilde{\nabla }\widetilde{\eta
}_{\alpha }+f_{0}^{\alpha }\theta \otimes \widetilde{\nabla
}^{2}\widetilde{\eta }_{\alpha }\}
\endaligned \tag{A3}
$$
where the property that $\nabla \theta =0$ is used. From (2.17)
and the Sasakian conditions of both $M$ and $N$, we have
$$
f_{0k}^{\alpha }=f_{k0}^{\alpha },\quad f_{0\overline{k}}^{\alpha
}=f_{\overline{k}0}^{\alpha }.  \tag{A4}
$$
In terms of (A4), (2.35) and (2.41), we derive that
$$
\aligned f_{0k\overline{k}}^{\alpha }+f_{0\overline{k}k}^{\alpha
}&=f_{k0\overline{k}
}^{\alpha }+f_{\overline{k}0k}^{\alpha } \\
&=f_{k\overline{k}0}^{\alpha }+f_{k}^{\beta }\widehat{R}_{\beta
\gamma \overline{\delta }}^{\alpha }(f_{\overline{k}}^{\gamma
}f_{0}^{\overline{ \delta }}-f_{0}^{\gamma
}f_{\overline{k}}^{\overline{\delta }})+f_{\overline{
k}k0}^{\alpha }+f_{\overline{k}}^{\beta }\widehat{R}_{\beta \gamma
\overline{\delta }}^{\alpha }(f_{k}^{\gamma
}f_{0}^{\overline{\delta }}-f_{0}^{\gamma
}f_{k}^{\overline{\delta }}) \\
&=f_{k\overline{k}0}^{\alpha }+f_{\overline{k}k0}^{\alpha
}+f_{k}^{\beta } \widehat{R}_{\beta \gamma \overline{\delta
}}^{\alpha }(f_{\overline{k} }^{\gamma }f_{0}^{\overline{\delta
}}-f_{0}^{\gamma }f_{\overline{k}}^{\overline{\delta
}})+f_{\overline{k}}^{\beta }\widehat{R}_{\beta \gamma
\overline{\delta }}^{\alpha }(f_{k}^{\gamma
}f_{0}^{\overline{\delta } }-f_{0}^{\gamma
}f_{k}^{\overline{\delta }}). \endaligned\tag{A5}
$$
It is easy to see from (A2) and (A5) that there exists a constant
$C$ such that
$$
\left. |f_{000}^{\alpha }+f_{0k\overline{k}}^{\alpha
}+f_{0\overline{k} k}^{\alpha }|\leq C\sum_{\alpha }|f_{0}^{\alpha
}|\right. \tag{A6}
$$
on a fixed open subset $W$ of $M$. Consequently
$$
|\bigtriangleup _{M}f_{0}^{\alpha }|\leq C^{\prime }\sum_{\alpha
}\big(|f_{0}^{\alpha }|+|\nabla f_{0}^{\alpha }|\big), \tag{A7}
$$
that is, $\{f_{0}^{\alpha }\}$ satisfies the structural
assumptions of Aronszajn-Cordes. Since $f_{0}^{\alpha }=0$ on
$W\cap U$, then $ f_{0}^{\alpha }=0$ on $W$, and thus we may
conclude that $f$ is foliated on $M$.\qed
\enddemo

\remark{Remark A2} From the proof of Theorem A1, we see that the
above continuation result about the foliated property still holds
if $f$ only satisfies $tr_{g_{\theta }}(\pi _{\widetilde{H}}\beta
)=0$.
\endremark

For a smooth map $f:(M^{2m+1},H(M),J,\theta )\rightarrow
(N^{2n+1}, \widetilde{H}(N),\widetilde{J},\widetilde{\theta })$,
we introduce the following two $1$-forms:
$$
\widetilde{\rho }_{3}=-(f_{00}^{\alpha }f_{0}^{\overline{\alpha
}}+f_{00}^{\overline{\alpha }}f_{0}^{\alpha })\theta
+\widetilde{\rho }_{1} \tag{A8}
$$
where $\widetilde{\rho }_{1}$ is given by (8.3), and
$$
\widetilde{\rho }_{4}=-(f_{00}^{0}f_{0}^{0}\theta
+f_{0k}^{0}f_{0}^{0}\theta
^{k}++f_{0\overline{k}}^{0}f_{0}^{0}\theta ^{\overline{k}}).
\tag{A9}
$$

\proclaim{Lemma A2} Let $f:M\rightarrow N$ be a map between two
Sasakian manifolds. Then
$$
\aligned \delta \widetilde{\rho }_{3}&=2|f_{00}^{\alpha
}|^{2}+2|f_{0k}^{\alpha }|^{2}+2|f_{0\overline{k}}^{\alpha
}|^{2}+(f_{000}^{\alpha }+f_{k\overline{k}
0}^{\alpha }+f_{\overline{k}k0}^{\alpha })f_{0}^{\overline{\alpha }} \\
&+(f_{000}^{\overline{\alpha
}}+f_{k\overline{k}0}^{\overline{\alpha }}+f_{
\overline{k}k0}^{\overline{\alpha }})f_{0}^{\alpha
}-\widehat{R}(df_{L, \widetilde{H}}(\xi )\wedge
df_{H,\widetilde{H}}(\eta _{k}),\overline{df_{L,
\widetilde{H}}(\xi )\wedge df_{H,\widetilde{H}}(\eta _{k})}) \\
&-\widehat{R}(df_{L,\widetilde{H}}(T)\wedge
df_{H,\widetilde{H}}(\eta _{
\overline{k}}),\overline{df_{L,\widetilde{H}}(T)\wedge
df_{H,\widetilde{H} }(\eta _{\overline{k}})}).
\endaligned\tag{A10}
$$
If $f$ is foliated, then
$$
\delta
\widetilde{\rho}_{4}=|f_{00}^{0}|^{2}+|f_{0k}^{0}|^{2}+|f_{0
\overline{k}}^{0}|^{2}+f_{0}^{0}(f_{000}^{0}+f_{k\overline{k}0}^{0}+f_{
\overline{k}k0}^{0}). \tag{A11}
$$
\endproclaim
\demo{Proof} From (8.4) and (A8), we immediately get (A10). Now
suppose $f$ is foliated. Then (2.14) yield
$$
f_{0k}^{0}=f_{k0}^{0},\quad
f_{0\overline{k}}^{0}=f_{\overline{k}0}^{0}. \tag{A12}
$$
Using (A12), (2.24), we deduce from (A9) that
$$
\aligned \delta \widetilde{\rho
}_{4}&=|f_{00}^{0}|^{2}+|f_{0k}^{0}|^{2}+|f_{0
\overline{k}}^{0}|^{2}+f_{0}^{0}(f_{000}^{0}+f_{0k\overline{k}}^{0}+f_{0
\overline{k}k}^{0}) \\
&=|f_{00}^{0}|^{2}+|f_{0k}^{0}|^{2}+|f_{0\overline{k}
}^{0}|^{2}+f_{0}^{0}(f_{000}^{0}+f_{k0\overline{k}}^{0}+f_{\overline{k}
0k}^{0}) \\
&=|f_{00}^{0}|^{2}+|f_{0k}^{0}|^{2}+|f_{0\overline{k}
}^{0}|^{2}+f_{0}^{0}(f_{000}^{0}+f_{k\overline{k}0}^{0}+f_{\overline{k}
k0}^{0}).\endaligned
$$
\qed
\enddemo
\remark{Remark A3} We only use the Sasakian condition on $M$ to
derive (A11), so it is still valid if $N$ is any pseudo-Hermitian
manifold.
\endremark

\proclaim{Theorem A3} Let $M$ and $N$ be two compact Sasakian
manifolds and let $f:M\rightarrow N$ be a pseudo-Hermitian
harmonic map. If $ (N^{2n+1},\widetilde{\nabla })$ has
non-positive horizontal sectional curvature, then $f$ is a
foliated special $(H,\widetilde{H})$-harmonic map with $df(\xi
)=\lambda \widetilde{\xi }$ for some constant $\lambda $.
\endproclaim
\demo{Proof} Using Lemma A2 and the divergence theorem, we get
$$
0\geq \int_{M}\{2|f_{00}^{\alpha }|^{2}+2|f_{0k}^{\alpha
}|^{2}+2|f_{0\overline{k}}^{\alpha }|^{2}\}dv_{\theta },
$$
and thus
$$
f_{00}^{\alpha }=f_{0k}^{\alpha }=f_{0\overline{k}}^{\alpha }=0.
\tag{A13}
$$
Consequently, $f$ is a $(H,\widetilde{H})$-harmonic map with split
horizontal second fundamental form. Therefore Lemma 8.1 implies
that $f$ is foliated.

Since $f$ is both foliated and pseudo-Hermitian harmonic, we know
that (A11) holds and becomes
$$
\delta \widetilde{\rho
}_{4}=|f_{00}^{0}|^{2}+|f_{0k}^{0}|^{2}+|f_{0
\overline{k}}^{0}|^{2}.  \tag{A14}
$$
Applying the divergence theorem to (A14), we obtain
$$
f_{00}^{0}=f_{0k}^{0}=f_{0\overline{k}}^{0}=0.  \tag{A15}
$$
From (A2) and (A15), we find that $f$ is special in the sense of
Definition 3.2, and $df(\xi )=\lambda \widetilde{\xi}$ with
$\lambda $ constant. \qed
\enddemo

{\bf B. An explicit formulation for special
$(H,\widetilde{H})$-harmonic maps}
\vskip 0.2cm

Now we want to give the explicit formulations for both the special
$(H,\widetilde{H})$-harmonic map equation and its parabolic
version, which are convenient for proving the existence theory of
special $(H,\widetilde{H})$-harmonic maps between two
pseudo-Hermitian manifolds $(M^{2m+1},H(M),J,\theta )$ and
$(N^{2n+1},\widetilde{H}(N), \widetilde{J},\widetilde{\theta })$.
As in the theory of harmonic maps (cf. [Li]), one can always
assume, in view of the Nash embedding theorem, that
$I:(N,g_{\widetilde{\theta }})\hookrightarrow (R^K,g_{E})$ is an
isometric embedding in some Euclidean space, where $g_{E}$ denotes
the standard Euclidean metric. Let $\nabla $ and
$\widetilde{\nabla }$ denote the Tanaka-Webster connections of $M$
and $N$ respectively, and let $\nabla ^{\widetilde{\theta }}$ and
$D$ denote the Levi-Civita connections of $(N,g_{\widetilde{\theta
}})$ and $(R^{K},g_{E})$ respectively.

For a map $f:(M,\nabla )\rightarrow (N,\widetilde{\nabla })$
between the two manifolds, the second fundamental form of $f$ with
respect to $(\nabla ,\widetilde{\nabla })$ is defined by
$$
\beta (f;\nabla ,\widetilde{\nabla })(X,Y)=\widetilde{\nabla }
_{Y}df(X)-df(\nabla _{Y}X)  \tag{B1}
$$
for any vector fields $X$, $Y$ on $M$. Applying the composition
formula for second fundamental forms (see Proposition 2.20 on page
16 of [EL]) to the maps $f:(M,\nabla )\rightarrow
(N,\widetilde{\nabla })$ and $I:(N,\widetilde{\nabla })\rightarrow
(R^{K},D)$, we have
$$
\beta (I\circ f;\nabla ,D)(\cdot ,\cdot )=dI\big(\beta (f;\nabla
,\widetilde{\nabla })(\cdot ,\cdot )\big)+\beta
(I;\widetilde{\nabla },D)\big(df(\cdot ),df(\cdot )\big). \tag{B2}
$$
Define a $2$-tensor field $S$ on $N$ by
$$
S(Z_{1},Z_{2})=\nabla _{Z_{2}}^{\widetilde{\theta
}}Z_{1}-\widetilde{\nabla }_{Z_{2}}Z_{1}  \tag{B3}
$$
where $Z_{1}$, $Z_{2}$ are any vector fields on $N$. Therefore
$$
\beta (I;\widetilde{\nabla },D)(\cdot ,\cdot )=\beta (I;\nabla
^{\widetilde{\theta }},D)(\cdot ,\cdot )+dI\big(S(\cdot ,\cdot
)\big). \tag{B4}
$$
Note that $\beta (I;\nabla ^{\widetilde{\theta }},D)$ is the usual
second fundamental form of the submanifold
$I:(N,g_{\widetilde{\theta }})\hookrightarrow (R^{K},g_{E})$. For
simplicity, we shall identify $N$ with $I(N)$, and write $I \circ
f$ as $u$, which is a map from $M$ to $R^{K}$. Set
$$
\tau _{H}(u;\nabla ,D)=tr_{g_{\theta }}\big(\beta
(u;\nabla,D)|_{H}\big). \tag{B5}
$$

\proclaim{Lemma B1} Suppose $f:M\rightarrow N$ is a map between
pseudo-Hermitian manifolds and $N$ is Sasakian. Suppose $I
:N\hookrightarrow R^{K}$ is an isometric embedding. Set $u=I\circ
f $. Then $f$ is a special $(H,\widetilde{H})$-harmonic map if and
only if
$$
\tau _{H}(u;\nabla ,D)-tr_{g_{\theta }}\beta (I;\nabla ^{
\widetilde{\theta }},D)(df_{H},df_{H})-tr_{g_{\theta }}dI\big
(S(df_{H},df_{H})\big)=0
$$
where $df_{H}$ denotes the restriction of $df$ to $H(M)$.
\endproclaim
\demo{Proof} From (B2), (B4) and (B5), we get
$$
\tau _{H}(u;\nabla ,D)=dI\big(\tau _{H}(f)\big)+tr_{g_{\theta
}}\beta (I;\nabla ^{\widetilde{\theta
}},D)(df_{H},df_{H})+tr_{g_{\theta }}dI\big(S(df_{H},df_{H})\big)
\tag{B6}
$$
where $\tau_{H}(f)=tr_{g_{\theta }}\big(\beta(f;\nabla
,\widetilde{\nabla })|_{H}\big)$ . We know from Corollary 3.3 that
if $N$ is Sasakian, then $f$ is a special
$(H,\widetilde{H})$-harmonic map if and only if $\tau _{H}(f)=0$.
Consequently this proposition follows immediately from (B6). \qed
\enddemo

Suppose now that $N$ is a compact Sasakian manifold. By
compactness of $N$, there exists a tubular neighborhood $B(N)$ of
$N$ in $R^{K}$ which can be realized as a submersion $\Pi
:B(N)\rightarrow N$ over $N$. Actually the projection map $\Pi $
is simply given by mapping any point in $B(N)$ to its closest
point in $N$. Clearly its differential $d\Pi
_{y}:T_{y}R^{K}\rightarrow T_{y}R^{K}$ when evaluate at a point
$y\in N$ is given by the identity map when restricted to the
tangent space $TN$ of $N$ and maps all the normal vectors to $N$
to the zero vector. Since $\Pi\circ I=I:N\hookrightarrow R^{K}$,
we have
$$
\beta (I;\nabla ^{\widetilde{\theta }},D)(\cdot ,\cdot )=d\Pi
(\beta (I;\nabla ^{\widetilde{\theta }},D)(\cdot ,\cdot ))+\beta
(\Pi ;D,D)(dI,dI)
$$
and thus
$$
\beta (I;\nabla ^{\widetilde{\theta }},D)=\beta (\Pi ;D,D)(d
I,dI).  \tag{B7}
$$
The tensor field $S$ may be extended to a tensor field
$\widehat{S}=\Pi ^{\ast }(dI\circ S)$ on $B(N)$, that is,
$$
\widehat{S}(W_{1},W_{2})=dI\big(S(d\Pi (W_{1}),d\Pi (W_{2}))\big)
\tag{B8}
$$
for any $W_{1},W_{2}\in TB(N)$. Let $\{y^{a}\}_{1\leq a\leq K}$ be
the natural Euclidean coordinate system of $R^{K}$. Set
$u^{a}=y^{a}\circ u$, $\Pi ^{a}=y^{a}\circ \Pi $, and write
$\widehat{S}(\cdot ,\cdot )=\widehat{S}^{a}(\cdot ,\cdot
)\frac{\partial }{\partial y^{a}}$. Choose an orthonormal frame
field $\{e_{A}\}_{A=0,1,...,2m}$ around any point of $M$ such that
$span\{e_{A}\}_{1\leq A\leq 2m}=H(M)$. By definition of the second
fundamental form , we have
$$
\tau _{H}(u;\nabla ,D)=\bigtriangleup _{H}u^{a}\frac{\partial
}{\partial y^{a}},  \tag{B9}
$$
and
$$
\aligned tr_{g_{\theta }}\beta (I;\nabla ^{\widetilde{\theta }
},D)(df_{H},df_{H}) &=tr_{g_{\theta }}\beta (\Pi
;D,D)(du_{H},du_{H})\\
&=\sum_{A=1}^{2m}\frac{\partial ^{2}\Pi ^{a}}{\partial
y^{b}\partial y^{c}}
e_{A}(u^{b})e_{A}(u^{c})\frac{\partial }{\partial y^{a}}\\
&=\Pi _{bc}^{a}\langle \nabla _{H}u^{b},\nabla _{H}u^{c}\rangle
\frac{\partial }{\partial y^{a}}\endaligned  \tag{B10}
$$
where $\Pi _{bc}^{a}=\frac{\partial ^{2}\Pi ^{a}}{\partial
y^{b}\partial y^{c}}$. In terms of (B8), we derive that
$$
\aligned tr_{g_{\theta }}dI\big(S(df_{H},df_{H})\big)
&=tr_{g_{\theta }}\widehat{S}
^{a}(du_{H},du_{H})\frac{\partial }{\partial y^{a}} \\
&=\sum_{A=1}^{2m}\widehat{S}^{a}(\frac{\partial }{\partial
y^{b}},\frac{
\partial }{\partial y^{c}})e_{A}(u^{b})e_{A}(u^{c})\frac{\partial }{\partial
y^{a}} \\
&=\widehat{S}_{bc}^{a}\langle \nabla _{H}u^{b},\nabla
_{H}u^{c}\rangle \frac{\partial }{\partial y^{a}}\endaligned
\tag{B11}
$$
where $\widehat{S}_{bc}^{a}=\widehat{S}^{a}(\frac{\partial
}{\partial y^{b}}, \frac{\partial }{\partial y^{c}})$. It follows
from (B6), (B9), (B10) and (B11) that
$$
dI(\tau _{H}(f))=(\bigtriangleup _{H}u^{a}-\Pi _{bc}^{a}\langle
\nabla _{H}u^{b},\nabla _{H}u^{c}\rangle
-\widehat{S}_{bc}^{a}\langle \nabla _{H}u^{b},\nabla
_{H}u^{c}\rangle )\frac{\partial }{\partial y^{a}}. \tag{B12}
$$
In view of (B12), we obtain
\proclaim{Proposition B2} Let $M$,
$N$, $f$ and $u$ be as in Lemma B1. Then $ f $ is a special
$(H,\widetilde{H})$-harmonic map if and only if
$$
\bigtriangleup _{H}u^{a}-\Pi _{bc}^{a}\langle \nabla
_{H}u^{b},\nabla _{H}u^{c}\rangle -\widehat{S}_{bc}^{a}\langle
\nabla _{H}u^{b},\nabla _{H}u^{c}\rangle =0
$$
where $\Pi _{bc}^{a}=\frac{\partial ^{2}\Pi ^{a}}{\partial
y^{b}\partial y^{c}}$ and
$\widehat{S}_{bc}^{a}=\widehat{S}^{a}(\frac{\partial }{\partial
y^{b}},\frac{\partial }{\partial y^{c}})$ $(1\leq a,b,c\leq K)$
are smooth functions defined on $ B(N)\subset R^{K}$.
\endproclaim

In section 7, we study the existence problem of the special
$(H,\widetilde{H})$-harmonic map equation $\tau _{H}(f)=0$ by
solving the corresponding subelliptic heat flow (7.5), that is,
$$
\cases \frac{\partial f_t}{\partial t}&=\tau_H(f_t)\\
f|_{t=0}&=h \endcases
$$
for some map $h:M\rightarrow N$. Inspired by the above explicit
formulation, we will establish the fact that in order to solve
(7.5), it suffices to solve the following system
$$
\cases \frac{\partial u^{a}}{\partial t}&=\bigtriangleup
_{H}u^{a}-\Pi _{bc}^{a}\langle \nabla _{H}u^{b},\nabla
_{H}u^{c}\rangle -\widehat{S}
_{bc}^{a}\langle \nabla _{H}u^{b},\nabla _{H}u^{c}\rangle,  \\
u^{a}|_{t=0}&=h^{a}, \quad(1\leq a,b,c\leq K),\endcases\tag{B13}
$$
where $h^{a}=y^{a}\circ h$. Let us define a map $\rho
:B(N)\rightarrow R^{K}$ by
$$
\rho (y)=y-\Pi (y),\quad y\in B(N).
$$
Obviously, $\rho (y)$ is normal to $N$ and $\rho (y)=0$ if and
only if $y\in N$.

\proclaim {Lemma B3} Let $u(x,t)=(u^{a}(x,t))$ $((x,t)\in M\times
\lbrack 0,T_{0}))$ be a solution of (B13) with initial condition $
h=(h^{a}):M\rightarrow R^{K}$. Then the quantity
$$
\int_{M}|\rho (u(x,t)|^{2}dv_{\theta }
$$
is a nonincreasing function of $t$. In particular, if $h(M)\subset
N$, then $ u(x,t)\in N$ for all $(x,t)\in M\times \lbrack
0,T_{0})$.
\endproclaim
\demo{Proof}Since $\rho (y)=y-\Pi (y)$, we have
$$
\rho _{b}^{a}=\delta _{b}^{a}-\Pi _{b}^{a}  \tag{B14}
$$
and
$$
\rho _{bc}^{a}=-\Pi _{bc}^{a}  \tag{B15}
$$
where $\rho _{b}^{a}=\frac{\partial \rho ^{a}}{\partial y^{b}}$
and $\rho _{bc}^{a}=\frac{\partial ^{2}\rho ^{a}}{\partial
y^{b}\partial y^{c}}$. By applying the composition law ([EL]) to
the maps $u_{t}:(M,\nabla )\rightarrow (B(N),D)$ and $\rho
:(B(N),D)\rightarrow (R^{K},D)$, we have
$$
\bigtriangleup _{H}\rho (u)=d\rho (\bigtriangleup
_{H}u)+tr_{g_{\theta }}\beta (\rho ;D,D)(du_{H},du_{H}). \tag{B16}
$$
Using (B16), (B14), (B15) and (B13), we derive that
$$
\aligned (\bigtriangleup _{H}\rho (u))^{a}&=\rho
_{b}^{a}\bigtriangleup _{H}u^{b}+\rho _{bc}^{a}\langle \nabla
_{H}u^{b},\nabla
_{H}u^{c}\rangle\\
&=\bigtriangleup _{H}u^{a}-\Pi _{b}^{a}\bigtriangleup
_{H}u^{b}-\Pi _{bc}^{a}\langle \nabla _{H}u^{b},\nabla
_{H}u^{c}\rangle\\
&=\frac{\partial u^{a}}{\partial t}+\widehat{S}_{bc}^{a}\langle
\nabla _{H}u^{b},\nabla _{H}u^{c}\rangle -\Pi
_{b}^{a}\bigtriangleup _{H}u^{b}\\
&=\rho _{b}^{a}\frac{\partial u^{b}}{\partial t}+\Pi
_{b}^{a}(\frac{\partial u^{b}}{\partial t}-\bigtriangleup
_{H}u^{b})+\widehat{S}_{bc}^{a}\langle \nabla _{H}u^{b},\nabla
_{H}u^{c}\rangle\endaligned\tag{B17}
$$
Since $d\Pi (\frac{\partial u}{\partial t}-\bigtriangleup _{H}u)$
and $\widehat{S}(du_{H},du_{H})$ are tangent to $N$ and $\rho (u)$
is normal to $N$, we find from (B17) that
$$
\rho ^{a}(u)\left( \bigtriangleup _{H}\rho (u)\right) ^{a}=\rho
^{a}(u)\rho _{b}^{a}(u)\frac{\partial u^{b}}{\partial t}.\tag{B18}
$$
Using (B18), we deduce that
$$
\aligned \frac{\partial }{\partial t}\int_{M}(\rho
^{a}(u))^{2}dv_{\theta } &=2\int_{M}\rho ^{a}(u)\rho
_{b}^{a}(u)\frac{\partial u^{b}}{\partial t}
dv_{\theta } \\
&=2\int_{M}\rho ^{a}(u)\left( \bigtriangleup _{H}\rho (u)\right)
^{a}dv_{\theta } \\
&=-2\int_{M}|\nabla _{H}\rho (u)|^{2}dv_{\theta } \\
&\leq 0\endaligned
$$
which implies that $\int_{M}|\rho (u)|^{2}dv_{\theta }$ is
decreasing in $t$.\qed
\enddemo

In terms of (B12) and Lemma B3, we conclude that

\proclaim{Theorem B4} Let $h:M\rightarrow N\subset R^{K}$ be a
smooth map given by $h=(h^{1},...,h^{K})$ in the Euclidean
coordinates. If $u:M\times \lbrack 0,T_{0})\rightarrow N\subset
R^{K}$ \ is a solution of the following system
$$
\frac{\partial u^{a}}{\partial t}=\bigtriangleup _{H}u^{a}-\Pi
_{bc}^{a}\langle \nabla _{H}u^{b},\nabla _{H}u^{c}\rangle
-\widehat{S}_{bc}^{a}\langle \nabla _{H}u^{b},\nabla
_{H}u^{c}\rangle ,\quad 1\leq a\leq K,
$$
with initial condition $(u^{a}(x,0))=(h^{a}(x))$ for all $x\in M$,
then $u$ solves the subelliptic heat flow
$$
\frac{\partial u}{\partial t}=\tau _{H}(u)
$$
with initial condition $u(x,0)=h(x)$.
\endproclaim

{\bf Acknowledgments}: The author would like to thank Professors
Aziz El Kacimi and E. Macas-Virg\'{o}s for their useful
discussions on the homotopy invariance of basic cohomology of
Riemannian foliations. He would also like to thank Ping Cheng for
helpful conversations.

\vskip 0.5 true cm

\vskip 1 true cm

School of Mathematical Science

and

Laboratory of Mathematics for Nonlinear Science

Fudan University,

Shanghai 200433, P.R. China

\vskip 0.2 true cm yxdong\@fudan.edu.cn

\enddocument